\documentclass{amsart}

\usepackage[all]{xy}
\usepackage{latexsym}
\usepackage{amssymb}
\usepackage{amsmath}
\usepackage{amsthm}
\usepackage{amscd}
\usepackage{fancyhdr}
\usepackage[dvips]{graphicx}
\usepackage{a4wide}

\theoremstyle{plain}
\newtheorem{pro}{Proposition}
\newtheorem{thm}[pro]{Theorem}
\newtheorem{lem}[pro]{Lemma}
\newtheorem{cor}[pro]{Corollary}
\theoremstyle{definition}
\newtheorem*{dei}{Definition}
\newtheorem*{Rq}{\sc Remark}
\newtheorem*{Ex}{Example}
\newtheorem*{Exs}{Examples}

\newcommand{\Id}{\mathrm{Id}}
\newcommand{\KK}{k}
\newcommand{\bc}{\boxtimes}

\newcommand{\bt}{\boxtimes}

\newcommand{\R}{\mathcal{R}}
\newcommand{\End}{\mathrm{End}}

\newcommand{\As}{\mathcal{A}s}

\newcommand{\B}{\mathcal{B}}
\newcommand{\oC}{\bar{\mathcal{C}}}

\newcommand{\oPo}{\overline{\Po}}
\newcommand{\oR}{\overline{R}}
\newcommand{\uI}{\underline{I}}
\newcommand{\Bi}{\mathcal{B}i}
\newcommand{\BLi}{\mathcal{B}i\mathcal{L}ie}
\newcommand{\IBi}{\varepsilon\mathcal{B}i}

\newcommand{\NN}{\mathbb{N}}
\newcommand{\ZZ}{\mathbb{Z}}

\newcommand{\kj}{{\bar{k},\, \bar{\jmath}}}
\newcommand{\oi}{{\bar{\imath}}}
\newcommand{\oj}{{\bar{\jmath}}}
\newcommand{\ok}{{\bar{k}}}
\newcommand{\ol}{{\bar{l}}}
\newcommand{\coker}{\mathop{\mathrm{coker}}}

\newcommand{\Sy}{\mathbb{S}}
\newcommand{\CA}{\mathcal{A}}
\newcommand{\CB}{\mathcal{B}}
\newcommand{\C}{\mathcal{C}om}
\newcommand{\A}{\mathcal{A}s}
\newcommand{\D}{\mathcal{D}ias}
\newcommand{\Pe}{\mathcal{P}erm}
\newcommand{\Po}{\mathcal{P}}
\newcommand{\F}{\mathcal{F}}
\newcommand{\ac}{\scriptstyle \textrm{!`}}

\newcommand{\II}{I}
\newcommand{\Qo}{\mathcal{Q}}

\newcommand{\Li}{\mathcal{L}ie}
\newcommand{\Sc}{\mathbb{S}^c}
\newcommand{\Pli}{\mathcal{P}reLie}
\newcommand{\Dend}{\mathcal{D}end}
\newcommand{\Co}{\mathcal{C}}
\newcommand{\Trias}{\mathcal{T}rias}
\newcommand{\Tridend}{\mathcal{T}ri\mathcal{D}end}
\newcommand{\Comtri}{\mathcal{C}om\mathcal{T}rias}

\newcommand{\Qu}{\mathcal{Q}uad}

\newcommand{\ot}{\otimes}
\newcommand{\otk}{\otimes_\KK}
\newcommand{\oth}{\otimes_H}
\newcommand{\Hom}{\mathrm{Hom}}
\newcommand{\Der}{\mathrm{Der}}
\newcommand{\sgn}{\mathrm{sgn}}
\newcommand{\Zi}{\mathcal{Z}inb}
\newcommand{\Le}{\mathcal{L}eib}
\newcommand{\h}{\mathrm{hom}_\bullet}
\newcommand{\en}{\mathrm{end}}

\newcommand{\epi}{\twoheadrightarrow}
\newcommand{\mono}{\rightarrowtail}

\newcommand{\BlackSq}{\, \blacksquare \,}
\newcommand{\Sq}{\, \square \,}
\newcommand{\Y}{\vcenter{\xymatrix@M=0pt@R=6pt@C=6pt{
\ar@{-}[dr] &  &\ar@{-}[dl]  \\
 &\ar@{-}[d] &  \\  & &}}}
\newcommand{\YY}{\vcenter{\xymatrix@M=0pt@R=6pt@C=6pt{
\ar@{-}[dr] &  &\ar@2{-}[dl]  \\
 &\ar@2{-}[d] &  \\  & &}}}
\newcommand{\YYY}{\vcenter{\xymatrix@M=0pt@R=6pt@C=6pt{
\ar@{-}[dr] &  &\ar@3{-}[dl]  \\
 &\ar@3{-}[d] &  \\  & &}}}
\newcommand{\cop}{\vcenter{\xymatrix@M=0pt@R=6pt@C=6pt{
 & \ar@{-}[d] & \\
 &\ar@{-}[dr] \ar@{-}[dl] &  \\  & &}}}

\newcommand{\copL}{\xymatrix@M=0pt@R=6pt@C=6pt{
 & \ar@{-}[d] & \\
 &\ar@{-}[dr] \ar@{-}[dl] &  \\  & &\\  & &\\  & &}}
\newcommand{\YL}{\vcenter{\xymatrix@M=0pt@R=6pt@C=6pt{
\ar@{-}[dr] &  &\ar@{-}[dl]  \\
 &\ar@{-}[d] &   \\  & &\\  & &\\  & &}}}
\newcommand{\YYL}{\vcenter{\xymatrix@M=0pt@R=6pt@C=6pt{
\ar@{-}[dr] &  &\ar@2{-}[dl]  \\
 &\ar@2{-}[d] &  \\  & &\\  & &\\  & &}}}
\newcommand{\LYY}{\vcenter{\xymatrix@M=0pt@R=6pt@C=6pt{
\ar@{-}[dr] &  &\ar@2{-}[dl]  \\
 &\ar@2{-}[d] &  \\  & &\\  & &\\  & &}}}

\newcommand{\XX}{\vcenter{\xymatrix@M=0pt@R=6pt@C=6pt{\ar@{-}[ddrr]&&\ar@{-}[ddll] \\ && \\ &&   }}}

\newcommand{\Ta}{\vcenter{\xymatrix@M=0pt@R=6pt@C=6pt{ \ar@{-}[dddrrr] && \ar@{-}[dl] &&  \\
&&& \ar@{-}[dl]  &  \\ &&&&  \ar@{-}[dl]  \\&&&  \ar@{-}[d] &
\\&&&& }}}
\newcommand{\Tb}{\vcenter{\xymatrix@M=0pt@R=6pt@C=6pt{  & \ar@{-}[dr]&&\ar@{-}[dl] \\
\ar@{-}[dr]&&\ar@{-}[dl]& \\&\ar@{-}[dr]&&\ar@{-}[dl]
\\&&\ar@{-}[d]& \\&&& }}}
\newcommand{\Tc}{\vcenter{\xymatrix@M=0pt@R=6pt@C=6pt{   \ar@{-}[dr]&&\ar@{-}[dl]& \\
&\ar@{-}[dr]&& \ar@{-}[dl] \\\ar@{-}[dr]&&\ar@{-}[dl]&
\\&\ar@{-}[d]&& \\&&& }}}
\newcommand{\Td}{\vcenter{\xymatrix@M=0pt@R=6pt@C=6pt{ && \ar@{-}[dr]&&\ar@{-}[dddlll] \\
 &\ar@{-}[dr]&&& \\ \ar@{-}[dr]&&&& \\& \ar@{-}[d]&&& \\&&&&  }}}
\newcommand{\Te}{\vcenter{\xymatrix@R=3pt@C=3pt{\ar@{-}[drdr] &&\ar@{-}[dl]  *=0{}
\ar@{-}[dr]&& \ar@{-}[ddll] \\ &&& *=0{}& \\&& *=0{} \ar@{-}[d]&&
\\&&&& }}}
\newcommand{\TaC}{\vcenter{\xymatrix@M=0pt@R=6pt@C=6pt{ \ar@{-}[ddddddrrrrrr] && \ar@{-}[dl] && && \\
&&& \ar@{-}[dl]  &&&  \\ &&&&  \ar@{-}[dl]&&  \\&&& &&&
\\&&&&\ar@{-}[dl]&& \\&&&&&\ar@{-}[dl]&\\&&&&&& }}}
\newcommand{\TreeL}{\vcenter{\xymatrix@M=0pt@R=5pt@C=5pt{ \ar@{-}[dr] &
&\ar@{-}[dl] & &  \\
& \ar@{-}[dr] & &\ar@{-}[dl]  & \\
& &\ar@{-}[d] & & \\
& & \\ & & }}}
\newcommand{\TreeR}{\vcenter{\xymatrix@M=0pt@R=5pt@C=5pt{
 & &\ar@{-}[dr] & & \ar@{-}[dl]  \\
& \ar@{-}[dr] & &\ar@{-}[dl]  & \\
& &\ar@{-}[d] & & \\
& & \\ & & }}}
\newcommand{\Lepi}{\twoheadleftarrow}
\newcommand{\Lmono}{\leftarrowtail}
\newcommand{\mona}{\mathcal{M}\textrm{on}_\CA}
\newcommand{\coma}{\mathcal{C}\textrm{omon}_\CA}

\newcommand{\draftnote}[1]{}

\newenvironment{proo}{\begin{trivlist} \item{\sc {Proof.}}}
  {\hfill $\square$ \end{trivlist}}

\long\def\symbolfootnote[#1]#2{\begingroup%
\def\thefootnote{\fnsymbol{footnote}}\footnote[#1]{#2}\endgroup}

\title{\bf Manin products, Koszul duality, Loday algebras and Deligne conjecture}
\author{\bf Bruno Vallette}
\begin{document}

\maketitle

\begin{center}
\textit{Dedicated to Jean-Louis Loday, on the occasion of his
sixtieth birthday\footnote{The title of this paper can be read
``How to use Manin's products to prove Deligne's conjecture for
Loday algebras with Koszul property''}}
\end{center}

\symbolfootnote[0]{2000 \emph{Mathematics Subject Classification.} 18D50, 18D10, 20C30, 55P48. \\
\emph{Keywords and phrases.} Monoidal categories, $2$-monoidal
categories, Algebras, Operads, Properads, Manin's black and white
products, Koszul duality, Cohomology operations, Deligne's
conjecture}

\begin{abstract}
In this article we give a conceptual definition of Manin products
in any category endowed with two coherent monoidal products. This
construction can be applied to associative algebras, non-symmetric
operads, operads,  colored operads, and properads presented by
generators and relations. These two products, called black and
white, are dual to each other under Koszul duality functor. We
study their properties and compute several examples of black and
white products for operads. These products allow us to define
natural operations on the chain complex defining cohomology
theories. With these operations, we are able to prove that
Deligne's conjecture holds for a general class of operads and is
not specific to the case of associative algebras. Finally, we
prove generalized versions of a few conjectures raised by M.
Aguiar and J.-L. Loday related to the Koszul property of operads
defined by black products. These operads provide infinitely many
examples for this generalized Deligne's conjecture.
\end{abstract}

\section*{Introduction}

In his works on quantum groups and non-commutative geometry, Yu.~
I. Manin defined two products in the category of quadratic
algebras. An associative algebra $A$ is called \emph{quadratic} if
it is isomorphic to a quotient algebra of the form $A=T(V)/(R)$,
where $T(V)$ is the free algebra on $V$ and where $(R)$ is the
ideal generated by $R\subset V^{\ot 2}$. Let $A=T(V)/(R)$ and
$B=T(W)/(S)$ be two quadratic algebras. Any quadratic algebra
 generated by the tensor product $V\ot W$ is determined by a  subspace of $(V\ot W)^{\ot
2}$. Since $R\subset V^{\ot 2}$ and $S\subset W^{\ot 2}$, one has
to introduce the isomorphism $(23) \, : \, V \ot V \ot W \ot W \to
V \ot W \ot V \ot W$ defined by the permutation of the second and
third terms. The black and white products were defined by Manin as
follows
\begin{eqnarray*}
A \bullet B &:=& T(V\ot W)/\big( (23)(R\ot S) \big), \\
A \circ B &:=& T(V\ot W)/\big( (23)(R\ot W^{\ot 2} + V^{\ot 2}\ot
S ) \big).
\end{eqnarray*}
Since $(23)$ is an isomorphism, many properties of the algebras
$A$ and $B$  remain true for their black and white products. For
instance, the white product of two quadratic algebras is equal to
their degreewise tensor product $A \otimes B := \oplus_{n \ge 0}
\, A_{n} \otimes B_{n}$. Therefore, one can apply the method of J.
Backelin \cite{B} to prove that the white product of two Koszul
algebras is again a Koszul algebra.\\

Koszul duality theory is a homological algebra theory developed by
S. Priddy \cite{Priddy} in 1970 for quadratic algebras. To a
quadratic algebra $A=T(V)/(R)$ generated by a finite dimensional
vector space $V$, one can associate the \emph{Koszul dual algebra}
$A^!:=T(V^*)/(R^\perp)$. Under this finite dimensional hypothesis,
we have $(A \circ B)^!=A^! \bullet B^!$, that is black and white
constructions are dual to each other under Koszul duality functor.
The main result of Manin is the following adjunction in the
category of finitely generated quadratic algebras
$$\Hom_{Q.\, Alg}(A \bullet B^!,\, C) \cong \Hom_{Q.\, Alg}(A,\, B\circ C ). $$
 Using the general properties of internal
cohomomorphisms, Manin proved that $A \bullet A^!$ is a Hopf
algebra and was able to realize well known quantum groups as black
products of an algebra with its Koszul dual.  For more properties
of Manin's products for quadratic algebras,
we refer the reader to the book of A. Polishchuk and L. Positselski \cite{PP}.\\

Koszul duality theory was later generalized to binary quadratic
operads by V. Ginzburg and M. Kapranov \cite{GK} in 1994. This
generalization comes from the fact that an operad, like an
associative algebra, is a monoid in a monoidal category. A
quadratic operad $\Po=\F(V)/(R)$ is a quotient of a free operad by
an ideal generated by a sub-$\Sy$-module $R$ of $\F_{(2)}(V)$, the
part of weight $2$ of $\F(V)$. Let $\Po=\F(V)/(R)$ and
$\Qo=\F(W)/(S)$ be two quadratic operads. A quadratic operad
generated by the tensor product $V\ot W$ is determined by a
subspace of $\F_{(2)}(V\ot W)$. Since $R\subset\F_{(2)}(V)$ and
$S\subset\F_{(2)}(W)$, we need a map from $\F_{(2)}(V)\ot
\F_{(2)}(W) $ to $\F_{(2)}(V\ot W)$. In the binary case, Ginzburg
and Kapranov mentioned in \cite{GK2} two maps $\Psi \, : \,
\F_{(2)}(V)\ot \F_{(2)}(W) \to \F_{(2)}(V\ot W)$ and $\Phi \, : \,
\F_{(2)}(V\ot W) \to \F_{(2)}(V)\ot \F_{(2)}(W)$ and defined the
black and white products for binary quadratic operads as follows.
\begin{eqnarray*}
\Po \bullet \Qo &:=& \F(V\ot W)/\big( \Psi(R\ot S) \big), \\
\Po \circ \Qo &:=& \F(V\ot W)/\big( \Phi^{-1}(R\ot W^{\ot 2} +
V^{\ot 2}\ot S ) \big).
\end{eqnarray*}

When the operad $\Po=\F(V)/(R)$ is a binary quadratic operad
generated by a finite dimensional $\Sy$-module $V$, they defined a
\emph{Koszul dual operad} by the formula
$\Po^!:=\F(V^\vee)/(R^\perp)$, where $V^\vee(2):=V^*(2)\otimes
\sgn_{\Sy_2}$  is the dual representation twisted by the signature
representation. As in the case of algebras, they proved that $(\Po
\circ \Qo)^!=\Po^! \bullet \Qo^!$ and they showed the adjunction
$$\Hom_{Bin. Q.\, Op.}(\Po \bullet \Qo^!,\, \mathcal{R}) \cong
\Hom_{Bin. Q.\, Op.}(\Po,\, \Qo\circ \mathcal{R} ), $$ in the
category of finitely generated binary quadratic operads.\\

From the properties of black and white products for associative
algebras and binary quadratic operads, a few natural questions
arise. Where do the functors $\Psi$ and $\Phi$ conceptually come
from ? Is the black or white product of two binary Koszul operads
still a Koszul operad ? Can one do non-commutative geometry with
an operad of the form $\Po \bullet \Po^!$ ? One can also add : is
it possible to recover classical operads as black or white
products of more simple operads? Can black and white products help
to describe the natural operations acting on cohomological spaces
? The aim of this paper is to answer these questions.
\\

Let us recall that Koszul duality theory of associative algebras
and binary quadratic operads was extended to various other
monoidal categories in the last few years. The following diagram
shows these monoidal categories were Koszul duality holds. They
are represented by the name of their monoids.
$$\xymatrix@M=6pt{\textrm{Associative Algebras\ } \ar@{>->}[r] \ar@{>->}[dr] & \textrm{Non-symmetric Operads}
\ar[d]^{\Sigma}& & &   \\
 & \textrm{Operads\ \ }  \ar@{>->}[r] \ar@{>->}[d]& \textrm{Colored operads} \\
&  \textrm{Dioperads} \ar@{>->}[d]& \\
&  \textrm{Properads} & }$$ Koszul duality for dioperads was
proved by W.L. Gan in \cite{Gan}, it was proved by P. Van der Laan
in \cite{VdL} for colored operads and by the author for properads
in \cite{V1}. A \emph{properad} is an object slightly smaller than
a prop  which encodes faithfully a large variety of algebraic
structures like bialgebras or Lie bialgebras, for instance (see
Appendix A for more details). We would like to emphasize that the
Koszul dual that appears naturally, without finite dimensional
assumptions, is a comonoid (coalgebra, cooperad, coproperad, etc
...). See Section~\ref{Koszul Duality Pattern} for more details. \\

To answer the first question about the conceptual definition of
the functors $\Psi$ and $\Phi$, we introduce a new notion of
 category endowed with $2$ coherent monoidal products. We call it
  \emph{$2$-monoidal category} in Section~\ref{2-monoidal
category}. This definition generalizes previous notions given by
A. Joyal and R. Street in \cite{JS} in the framework of braided
tensor categories and by C. Balteanu, Z. Fiedorowicz, R.
Schw\"anzl and R. Vogt in \cite{BFSV} in the framework of iterated
monoidal category and iterated loop spaces. All the examples given
above are monoids in a $2$-monoidal category. In a $2$-monoidal
category, we define the functors $\Psi$ and $\Phi$ by universal
properties. This allows us to define white products for monoids
presented by generators and relations in
Section~\ref{ManinProducts}. Since the Koszul dual is a comonoid,
we define a black product for comonoids presented by generators
and relations. (This notion is introduced and detailed in
Appendix~\ref{Categorical algebra}).\\

The white product defined here coincides with the one of Yu. I.
Manin for quadratic algebras $A\circ B$, with the one of
Ginzburg-Kapranov for binary quadratic operads $\Po\circ \Qo$ and
with the one of R. Berger, M. Dubois-Violette and M. Wambst
\cite{BergerDuboisWambst} for $N$-homogenous algebras. Note that
the white product is defined without homogenous assumption.
Therefore, one could apply them to non-homogenous cases. In this
sense, it would be interesting to study the properties of the
white products of Artin-Schelter algebras \cite{ArtinSchelter,
Palmieri}. \\

Under finite dimensional assumptions, the twisted linear dual of
the Koszul dual cooperad gives the Koszul dual operad defined by
\cite{GK}. Using this relation, we define a black product for
operads in Section~\ref{Manin products for operads}. We do several
computations of black and white products and show that some
classical operads can be realized as products of simpler operads.
All these examples are products of Koszul operads and the result
is again a Koszul operad. This fact is not true in general and we
provide a counterexample in Section~\ref{Counterexample}. Whereas
this property holds for associative algebras, it is not true here
because the functors $\Phi$ and $\Psi$
are not isomorphisms. \\

We extend the adjunction of Manin and Ginzburg-Kapranov and prove
that $\Po \bullet \Po^!$ is a Hopf operad. Since operads are
non-linear generalizations of associative algebras, the notion of
Hopf operad can be seen as a non-linear generalization of
bialgebras. Hopf operads of the form $\Po \bullet \Po^!$ can
provide new examples of ``quantum groups'', in the philosophy of
\cite{M}. This adjunction also allows us to understand the
algebraic structures on tensor products or spaces of morphisms of
algebras. For instance, it gives a description of the structure of
cohomology spaces. \\

Non-symmetric operads are operads without the action of the
symmetric group. One can symmetrize a non-symmetric operad to get
an operad. (It corresponds to the functor $\Sigma$ in the diagram
above). The image of a non-symmetric operad under functor is
called a \emph{regular} operad. We define black $\BlackSq$ and
white $\Sq$ square products for regular operads as the image of
black and white products of non-symmetric operads in
Section~\ref{BW Square}. In the case of binary quadratic regular
operads, the black square product given here corresponds to the
one introduced by K. Ebrahimi-Fard and L. Guo in \cite{EFG} (see
also J.-L. Loday \cite{L2}). We prove the same kind of results for
regular operads and square products than
the ones for operads and Manin's products. \\

The adjunction for black and white square products allows us to
construct natural operations on the chain complex defining the
cohomology of an algebra over a non-symmetric (regular) operad.
The example of associative algebras is very classical. Since the
introduction of this (co)chain complex by Hochschild in 1945, it
has been extensively studied. M. Gerstenhaber proved in the
sixties that the cohomology of any associative algebra is endowed
with two coherent products : the commutative cup product and a Lie
bracket. This structure is now called a \emph{Gerstenhaber
algebra}. (Gerstenhaber also used this Lie bracket to study
deformations of associative algebras. This led to the work of
Kontsevich on deformation-quantization of Poisson manifold thirty
years latter). In homotopy theory, there is a topological operad,
formed by configurations of disks in the plane and called the
\emph{little disks operad}, whose action permits to recognize
two-fold loop spaces. In 1976, F. Cohen proved that the operad
defined by the homology of the little disks operad is equal to the
operad coding Gerstenhaber algebras. Therefore the Hochschild
cohomology space is an algebra over the homology of the little
disks operad. This surprising link between algebra, topology and
geometry led Deligne to formulate the conjecture that this
relation can be lifted on (co)chain complexes, that is the
singular chain complex of the little disks operad acts on the
Hochschild (co)chain complex of an associative algebra. This
conjecture was proved by several researchers using different
methods. In the present paper, we take a transversal approach. We
prove that Deligne's conjecture holds for a general class of
operads and is not specific to the case of associative algebras.
Using Manin's products, we construct operations on the chain
complex of any algebra over an operad of this class. (To be more
precise, finitely generated binary non-symmetric Koszul operads
form this class). Since these operations verify the same relations
than the ones on the Hochschild (co)chain complex, Deligne's
conjecture is then proved with the same methods. \\

Since the white square product is the Koszul dual of the black
square product, we can compute the Koszul duals of operads defined
by black square product. The first example is the operad
$\Qu=\Dend\BlackSq\Dend$ defined by M. Aguiar and J.-L. Loday in
\cite{AL}. Using the explicit description of its Koszul dual and
the method of partition posets of \cite{V2}, we prove that it is
Koszul over $\ZZ$, which answers a conjecture of Aguiar-Loday.
Actually, with the same method, we show that the families
$\Dend^{\BlackSq n}$, $\D^{\Sq n}$ and $\Tridend^{\BlackSq n}$,
$\Trias^{\Sq n}$ are Koszul over
$\ZZ$. These families provide infinitely many examples for which Deligne's conjecture hold over $\ZZ$. \\

Appendix A is a survey on the notions of operads and properads.
Appendix B yields a categorical approach of algebra with monoids
and
comonoids ((co)ideal, (co)modules).\\

Unless stated otherwise, we work over a field $\KK$ of
characteristic $0$.

\tableofcontents

\section{2-monoidal categories}\label{2-monoidal category}

In this section, we define the general framework of
\emph{2-monoidal category} verified by the examples studied
throughout the text. The notion of 2-monoidal category given here
is a lax and more general version of the one given by A. Joyal and
R. Street in \cite{JS} and the one given by C. Balteanu, Z.
Fiedorowicz, R. Schw\"anzl and R. Vogt in \cite{BFSV}.

\subsection{Monoidal category}

We recall briefly the definitions of \emph{monoidal category},
\emph{lax monoidal functor} and \emph{monoid} in order to settle
the notations for the next section. We refer to the book of S.
MacLane \cite{MacLane} Chapter VII and to the article of J.
B\'enabou \cite{Be} for full references on the subject.

\begin{dei}[Monoidal category]
A monoidal category $(\CA,\, \bt,\, I,\, \alpha,\, r,\,l)$ is a
category $\CA$ equipped with a bifunctor $\bt$ : $\CA \times \CA
\to \CA$ and a family of isomorphisms
$$\alpha_{A, B, C} \ : \ (A\bt B)\bt C \xrightarrow{\sim} A \bt (B \bt C),$$
for every $A$, $B$ and $C$ in $\CA$. These isomorphisms are
supposed to verify the pentagon axiom. For every object $A$ in
$\CA$, there exists two isomorphisms $l_a\, : \, I\bt A \to A$ and
$r_a\, : \, A\bt I \to A$ compatible  with $\alpha$.
\end{dei}

\begin{Ex}
Let $(\CA,\, \bt,\, I,\, \alpha,\, r,\,l)$ be a monoidal category.
The cartesian product $\CA\times \CA$ is a monoidal category where
the monoidal product $\bt^2$ is defined by $(A,\, B) \bt^2 (A',\,
B'):= (A\bt B,\, A'\bt B')$. The unit is $(I,\, I)$. The
associative isomorphisms are given by
$$\xymatrix{\big((A,\, B) \bt^2 (A',\, B')\big)\bt^2 (A'',\, B'') = \big( (A\bt A')\bt A'',\, (B\bt B')\bt B'' \big)
\ar[d]^{(\alpha_{A, A', A''},\, \alpha_{B, B', B''})} \\
(A,\, B) \bt^2 \big((A',\, B')\bt^2 (A'',\, B'')\big) = \big( A\bt
(A'\bt A''),\, B\bt (B'\bt B'') \big).}$$ The other isomorphisms
are $l_{(A,\, B)}:=(l_A,\, l_B)$ and  $r_{(A,\, B)}:=(r_A,\,
r_B)$.
\end{Ex}

\begin{dei}[Monoid]
Let $(\CA,\, \bt, \, I)$ be a monoidal category. A \emph{monoid}
$(M, \, \mu,\, u)$ is an object $M$ of $\CA$ endowed with two
morphisms : an associative \emph{product} $\mu \, :\, M\bt M \to
M$ and a \emph{unit} $u\, :\, I \to M$.
\end{dei}

\begin{dei}[Lax monoidal functor]
A \emph{lax monoidal functor} is a functor $F$ between two
monoidal categories $(\CA,\, \bt_\CA,\, I_\CA)\to (\CB,\,
\bt_\CB,\, I_\CB)$ such that there exists a map $\iota \, : \,
I_\CB \to F(I_\CA)$ and a natural transformation
$$\varphi_{A,\, A'} \ : \ F(A)\bt_\CB F(A') \to F(A \bt_\CA A'),$$
for every $A,A'$ in $\CA$. This natural transformation is supposed
to be compatible with the associativity and the units of the
monoidal categories :
\begin{itemize}
\item{\underline{Associativity condition} : } For every $A$, $A'$
and $A''$ in $\CA$, the following diagram is commutative

$$\xymatrix@C=80pt{  \big( F(A)\bt_\CB F(A') \big) \bt_\CB F(A'')  \ar[d]^{\varphi_{A, A'}\bt_\CB id
\ } \ar[r]^{\alpha^\CB_{F(A), F(A'), F(A'')}}&   F(A)\bt_\CB \big(
F(A')  \bt_\CB F(A'') \big)
 \ar[d]^{id \bt_\CB \varphi_{A', A''}} \\
 F(A \bt_\CA A') \bt_\CB F(A'')  \ar[d]^{\varphi_{A\bt_\CA A', A''}}  &
 F(A)\bt_\CB F(A'\bt_\CA  A'')   \ar[d]^{\varphi_{A, A'\bt_\CA A''}}  \\
 F\big( (A \bt_\CA A')\bt_\CA A'' \big)   \ar[r]^{F\left( \alpha^\CA_{A, A', A''}\right)} &
 F\big(A \bt_\CA (A' \bt_\CA A'')\big). }$$

\item{\underline{Unit condition} : } For every $A$ in $\CA$, the
following diagram is commutative
$$\xymatrix@R=30pt@C=50pt{ I_\CB \bt_\CB F(A) \ar[drr]_{l^\CB_{F(A)}}
 \ar[r]^{\iota \bt_\CB F(A)}  & F(I_\CA)\bt_\CB F(A) \ar[r]^{\varphi_{I_\CA, A}}
  & F(I_\CA \bt_\CA A )  \ar[d]^{F(l^\CA_A)}   \\
& & F(A).   } $$ The same statement holds on the right hand side.
\end{itemize}
\end{dei}

The purpose of the definition of lax monoidal functors is to
preserve monoids.

\begin{pro}[\cite{Be}]
\label{produitdemonoides} Let $F \, : \, (\CA,\, \bt_\CA,\,
I_\CA)\to (\CB,\, \bt_\CB,\, I_\CB)$ be a lax monoidal functor and
let $(M,\, \mu,\, u)$ be a monoid in $\CA$. The image of $M$ under
$F$ is a monoid in $\CB$. The product $\widetilde{\mu}$ is defined
by
$$\widetilde{\mu} \ : \ F(M)\bt_\CB F(M) \xrightarrow{\varphi_{M, M}} F(M\bt_\CA M) \xrightarrow{F(\mu)}  F(M). $$
And the unit $\widetilde{u}$ is defined by
$$\widetilde{u} \ : \ I_\CB \xrightarrow{\iota} F(I_\CA)
\xrightarrow{F(u)} F(M).$$
\end{pro}

\subsection{Definition of lax $2$-monoidal category}

Motivated by the examples treated in the sequel, we define a
general notion of category with two compatible monoidal products.

\begin{dei}[Lax 2-monoidal category]
A \emph{lax 2-monoidal category} is a category $\left(
\mathcal{A},\, \boxtimes,\, I, \, \ot,\, K\right)$, such that both
$\left( \mathcal{A},\, \boxtimes,\, I \right)$ and  $(\CA,\,
\ot,\, K)$ are monoidal categories and such that the bifunctor
$\ot \, : \, \CA \times \CA \to \CA$ is a lax monoidal functor
with respect to the monoidal products $\boxtimes^2$ and
$\boxtimes$.
\end{dei}

The last assumption of the definition describes the compatibility
between the two monoidal structures. The next proposition makes it
more explicit.

\begin{pro}\label{Def:Lax 2-mono}
A \emph{lax 2-monoidal category} is a category $\left(
\mathcal{A},\, \boxtimes,\, I, \, \ot,\, K\right)$, such that both
$\left( \mathcal{A},\, \boxtimes,\, I \right)$ and  $(\CA,\,
\ot,\, K)$ are monoidal categories. These two monoidal structures
are related by a natural transformation called the
\emph{interchange law}
$$(A\ot A')\boxtimes (B\ot B') \xrightarrow{\varphi_{A,A',B,B'}}(A\boxtimes B)\ot (A'\boxtimes B'),$$
where $A$, $A'$, $B$ and $B'$ are in $\CA$. This interchange law
is supposed to be compatible with the associativity of the first
monoidal product $\bt$, that is

$$\xymatrix@C=80pt{   \big( (A\ot A')\bt (B\ot B')
\big)\bt (C\ot C') \ar[d]^{\varphi_{A,A',B,B'}\bt id  }
\ar[r]^{\alpha^\bt_{A\ot A', B\ot B' , C\ot C'}}   &
 (A\ot A')\bt \big( (B\ot B') \bt (C\ot C') \big)  \ar[d]^{id\bt \varphi_{B,B',C,C'}} \\
\big( (A\bt B)\ot (A' \bt B') \big)\bt (C\ot C')
\ar[d]^{\varphi_{A\bt B, A'\bt B', C, C'}} &
(A\ot A')\bt \big( (B\bt C) \ot (B'\bt C') \big)  \ar[d]^{\varphi_{A, A', B\bt C, B'\bt C'}}  \\
 \big( (A\bt B)\bt C    \big) \ot  \big( (A'\bt B')\bt C'    \big)
 \ar[r]^{\alpha^{\bt}_{A, B, C}\ot \alpha^{\bt}_{A', B', C'}\ }   &
\big( A\bt (B\bt C) \big) \ot  \big(A'\bt (B'\bt C') \big), }$$

 where $\alpha^\bt_{A, B, C}$ is the associativity morphism for
the monoidal product $\bt$ : $(A\bt B)\bt C \xrightarrow{\sim} A \bt (B\bt C)$.\\

There exists a map  $\iota \, : \, I \to I\ot I$ such that for
every $A$ and $A'$ in $\CA$, the following diagram is commutative
$$\xymatrix@R=30pt@C=50pt{ I \bt (A\ot A') \ar[drr]_{l^\bt_{(A \ot A')}}
\ar[r]^{\iota \bt id}  & (I\ot I) \bt (A\ot A')
\ar[r]^{\varphi_{I, I
, A, A'}}  & (I\ot A) \bt (I\ot A')  \ar[d]^{l^\ot_A \bt l^\ot_{A'}}   \\
& & F(A).   } $$ The same statement holds on the right hand side.
\end{pro}

\begin{proo}
The proof is a straightforward application of the definition.
\end{proo}

\begin{pro}\label{tenseur de monoides}
Let $\left( \mathcal{A},\, \boxtimes,\, I, \, \ot,\, K\right)$ be
a lax $2$-monoidal category. Consider two $\bt$-monoids $M$ and
$N$ in $\CA$. Their $\ot$-product $M\ot N$ is a $\bt$-monoid.
\end{pro}

\begin{proo}
It is a direct corollary of Definition~\ref{Def:Lax 2-mono} and
Proposition~\ref{produitdemonoides}.
\end{proo}


Motivated by the example of braided monoidal categories, A. Joyal
and R. Street gave the first notion of a category endowed with two
compatible monoidal products in \cite{JS}. In their definition,
the monoidal categories are non-necessarily strict but the
interchange law is supposed to be a natural isomorphism. This last
condition forces the two monoidal products to be isomorphic.

In order to model $n$-fold loop spaces, C. Balteanu, Z.
Fiedorowicz, R. Schw\"anzl and R. Vogt introduced in \cite{BFSV}
the notions of \emph{$n$-fold monoidal category}. Their notion of
$2$-fold monoidal category is, in some sense, a lax version of the
one given by Joyal and Street since they do not assume the
interchange law to be an isomorphism. But they require the
monoidal structures to be strict and the two units are equal.

In the definition of a lax $2$-monoidal category, we do not ask
the monoidal structures to be strict. The two units need not be
isomorphic. And the interchange law is not an isomorphism.
Therefore, the notion given here is a lax version of the one of
Joyal-Street and the one of  Balteanu-Fiedorowicz-Schw\"anzl-Vogt.
The definition of lax $2$-monoidal category was suggested by our
natural examples, that we make explicit in Section~\ref{Examples
of 2-monoidal categories}.

\subsection{Definition of $2$-monoidal category}

Working in the opposite category, we get the dual notion of
\emph{colax $2$-monoidal category}. Finally, we call a
\emph{$2$-monoidal category} a category which is both lax and
colax $2$-monoidal.

\begin{dei}[Comonoid]
A \emph{comonoid} $C$ is a monoid in the opposite category. It is
endowed with two morphisms : a coassociative \emph{coproduct} $C
\to C\bt C$ and a \emph{counit} $C\to I$.
\end{dei}

\begin{dei}[Colax monoidal functor]
A \emph{colax monoidal functor} is a functor $F$ between two
monoidal categories $(\CA,\, \bt_\CA,\, I_\CA)\to (\CB,\,
\bt_\CB,\, I_\CB)$ such that there exists a map $I_\CB \gets
F(I_\CA)$ and a natural transformation
$$\psi_{A,\, A'} \ : \ F(A)\bt_\CB F(A') \gets F(A \bt_\CA A').$$
This natural transformation is supposed to be compatible with the
associativity and the units of the monoidal categories.
Explicitly, these compatibilities are given by the reversed
diagrams defining a lax monoidal functor.
\end{dei}

The purpose of the definition of colax monoidal functors is to
preserve comonoids.

\begin{pro}[Be]
\label{images de comonoides} Let $F \, : \, (\CA,\, \bt_\CA,\,
I_\CA)\to (\CB,\, \bt_\CB,\, I_\CB)$ be a colax monoidal functor
and let $C$ be a comonoid in $\CA$. The image of $C$ under $F$ is
a comonoid in $\CB$.
\end{pro}

\begin{dei}[Colax 2-monoidal category]\label{Def:coLax 2-mono}
A \emph{colax 2-monoidal category} is a category $\left(
\mathcal{A},\, \boxtimes,\, I, \, \ot,\, K\right)$, such that both
$\left( \mathcal{A},\, \boxtimes,\, I \right)$ and  $(\CA,\,
\ot,\, K)$ are monoidal categories and such that the bifunctor
$\ot \, : \, \CA \times \CA \to \CA$ is a colax monoidal functor.
\end{dei}

A category $\left( \mathcal{A},\, \boxtimes,\, I, \, \ot,\,
K\right)$ is a colax $2$-monoidal category if it is endowed with
natural transformations, called \emph{the interchange laws},
$$(A\ot A')\boxtimes (B\ot B') \xleftarrow{\psi_{A,A',B,B'}}(A\boxtimes B)\ot (A'\boxtimes B'),$$
verifying the same commutative diagram than the one defining a lax
$2$-monoidal category, with the maps $\varphi$ replaced by the
maps $\psi$.

\begin{pro}\label{Produit de Comonoides}
Let $\left( \mathcal{A},\, \boxtimes,\, I, \, \ot,\, K\right)$ be
a colax $2$-monoidal category. Consider two $\bt$-comonoids $M$
and $N$ in $\CA$. Their $\ot$-product $M\ot N$ is a
$\bt$-comonoid.
\end{pro}

\begin{proo}
It is a direct corollary of Definition~\ref{Def:coLax 2-mono} and
Proposition~\ref{images de comonoides}.
\end{proo}

\begin{dei}[2-monoidal category] A \emph{2-monoidal category} is a
category $\left( \mathcal{A},\, \boxtimes,\, I, \, \ot,\,
K\right)$, such that both $\left( \mathcal{A},\, \boxtimes,\, I
\right)$ and  $(\CA,\, \ot,\, K)$ are two monoidal categories and
such that the bifunctor $\ot \, : \, \CA \times \CA \to \CA$ is a
lax and colax monoidal functor.
\end{dei}

\begin{dei}[Strong 2-monoidal category] A \emph{strong 2-monoidal
category} is a $2$-monoidal category where the bifunctor $\ot \, :
\, \CA \times \CA \to \CA$ is a strong monoidal functor, that is
the interchange laws are isomorphisms.
\end{dei}

\subsection{Examples of $2$-monoidal categories}\label{Examples of 2-monoidal categories}

In the this section, we study the relation between the composition
product $\bc$  and the Hadamard product $\oth$
 in the category of $\Sy$-bimodules and in the subcategories of $\Sy$-modules and $\KK$-modules.
These notions are recalled in Appendix~\ref{Alg-Op-Properad}.

\begin{pro}\label{Prop:examples 2-mono}
 The categories $(\KK\textrm{-Mod},\, \otk,\, \KK)$,
$(\Sy\textrm{-Mod}, \circ,\, I)$ and $(\Sy\textrm{-biMod}, \bc,\,
I )$ endowed with the Hadamard tensor products are 2-monoidal
categories. The first one is a strong 2-monoidal category and it
is a full sub-2-monoidal category of the second one, which is a
full sub-2-monoidal category of the last one.
\end{pro}

\begin{proo}

\begin{itemize}
\item In the first category, the two monoidal products are equal ,
that is $\boxtimes=\ot=\ot_k$. The interchange laws are given by
the twisting isomorphism $(23) \, : \, V_1\ot V_2 \ot V_3 \ot V_4
\xrightarrow{\sim} V_1\ot V_3 \ot V_2 \ot V_4$.

\item In the category of $\Sy$-modules with $\boxtimes=\circ$,
$\ot=\ot_H$ the first interchange law $\varphi_{V,V',W,W'}$ map
comes from the well defined natural map
\begin{eqnarray*}
&(V \ot V')\circ (W \ot W')(n) :=&\\
 & \displaystyle \left( \bigoplus_{i_1+\cdots +i_l=n} (V\ot
V')(l)\ot_k \big( (W\ot
W')(i_1) \ot_k \cdots \ot_k (W\ot W')(i_l) \big)\ot_{\Sy_{i_1}\times \cdots \times
\Sy_{i_l}}k [\Sy_n]    \right)_{\Sy_l}& \\
 &\downarrow& \\
& \displaystyle \left( \bigoplus_{i_1+\cdots +i_l=n}
V(l)\ot_k \big(W(i_1) \ot_k \cdots \ot_k W(i_l)\big)\ot_{\Sy_{i_1}\times
\cdots \times \Sy_{i_l}} k [\Sy_n]  \right)_{\Sy_l} \bigotimes& \\
&\displaystyle \left( \bigoplus_{i_1+\cdots +i_l=n} V'(l)\ot_k
\big( W'(i_1) \ot_k \cdots \ot_k W'(i_l)\big) \ot_{\Sy_{i_1}\times
\cdots \times \Sy_{i_l}} k [\Sy_n] \right)_{\Sy_l}=\\
&(V\circ W) \ot (V'\circ W')(n).&
\end{eqnarray*}
The other map corresponds to the transpose of this one. It is well
defined on invariants instead of coinvariants. Since we work over
a field $\KK$ of characteristic $0$, we use the classical
isomorphism between invariants and coinvariants to fix this.

\item In the last case, which includes the two first, the
interchange law map is the direct generalization of the one
written above. Its explicit description is

\begin{eqnarray*}
& \big( V\ot V'\big) \bc\big( W\ot W' \big)(m,\, n) := &\\
& \bigoplus_{N\in \mathbb{N}^*} \left( \bigoplus_{\ol,\, \ok,\,
\oj,\, \oi} \KK[\mathbb{S}_m]
 \otimes_{\mathbb{S}_\ol}
\big( V\ot V'\big)(\ol,\, \ok)\otimes_{\mathbb{S}_{\ok}}
k[\Sc_\kj] \otimes_{\mathbb{S}_\oj} \big( W\ot W' \big)(\oj,\,
\oi) \otimes_{\mathbb{S}_\oi} k[\mathbb{S}_n]
\right)_{\Sy_b^{\textrm{op}}\times \Sy_a}&\\
 &\downarrow &\\
& \bigoplus_{N\in \mathbb{N}^*} \Big( \bigoplus_{\ol,\, \ok,\,
\oj,\, \oi} \big( \KK[\mathbb{S}_m]
 \otimes_{\mathbb{S}_\ol}
V (\ol,\, \ok)\otimes_{\mathbb{S}_{\ok}} k[\Sc_\kj]
\otimes_{\mathbb{S}_\oj}  W(\oj,\, \oi) \otimes_{\mathbb{S}_\oi}
k[\mathbb{S}_n] \big) \ot&\\& \big( \KK[\mathbb{S}_m]
 \otimes_{\mathbb{S}_\ol}
V'(\ol,\, \ok)\otimes_{\mathbb{S}_{\ok}} k[\Sc_\kj]
\otimes_{\mathbb{S}_\oj}  W'(\oj,\, \oi) \otimes_{\mathbb{S}_\oi}
k[\mathbb{S}_n]
\big)\Big)_{\Sy_b^{\textrm{op}}\times \Sy_a}=&\\
 &\downarrow &\\
& \bigoplus_{N\in \mathbb{N}^*} \left( \bigoplus_{\ol,\, \ok,\,
\oj,\, \oi} \KK[\mathbb{S}_m]
 \otimes_{\mathbb{S}_\ol}
V (\ol,\, \ok)\otimes_{\mathbb{S}_{\ok}} k[\Sc_\kj]
\otimes_{\mathbb{S}_\oj}  W(\oj,\, \oi) \otimes_{\mathbb{S}_\oi}
k[\mathbb{S}_n] \right)_{\Sy_b^{\textrm{op}}\times \Sy_a} \ot&\\&
\bigoplus_{N\in \mathbb{N}^*} \left( \bigoplus_{\ol,\, \ok,\,
\oj,\, \oi} \KK[\mathbb{S}_m]
 \otimes_{\mathbb{S}_\ol}
V'(\ol,\, \ok)\otimes_{\mathbb{S}_{\ok}} k[\Sc_\kj]
\otimes_{\mathbb{S}_\oj}  W'(\oj,\, \oi) \otimes_{\mathbb{S}_\oi}
k[\mathbb{S}_n]
\right)_{\Sy_b^{\textrm{op}}\times \Sy_a}=&\\
&\big( V \bc W\big) \otimes  \big( V' \bc W'\big) (m,\, n). &
\end{eqnarray*}
\end{itemize}
Note that the first map preserves the shape of the underlying
graph of the composition, whereas the second one does not.
Therefore, this interchange law map is injective but not an
isomorphism. The reverse natural transformation $\big( V \bc
W\big) \otimes  \big( V' \bc W'\big) (m,\, n) \epi  \big( V\ot
V'\big) \bc\big( W\ot W' \big)(m,\, n)$ is given by the projection
on pairs of composition  of $\big( V \bc W\big) \otimes  \big( V'
\bc W'\big)$ based on the same 2-levelled graph (see
\ref{properads}). To such pairs, it is straightforward to
associate an element of $\big( V\ot V'\big) \bc\big( W\ot W'
\big)$. This map is the transpose of the first one. It is the
composite of an epimorphism with an isomorphism, therefore it is
an epimorphism.
\end{proo}

\begin{Rq}
In the same way, we can also show that the underlying category of
non-symmetric operads, $\frac{1}{2}$-props \cite{MarklVoronov},
dioperads \cite{Gan}, colored operads  are $2$-monoidal
categories. We refer to \cite{V4} Section 5 and to \cite{MK}
Section 9 for surveys of these notions.
\end{Rq}

\subsection{Bimonoids}\label{Bimonoids}

In this section, we define the notion of \emph{bimonoid} that
generalizes the notion of \emph{bialgebra} in any lax
$2$-monoidal category.\\

Let $\left( \mathcal{A},  \boxtimes, I, \ot, K \right)$ be a lax
$2$-monoidal category. Proposition~\ref{produitdemonoides} shows
that the category of $\bt$-monoids, denoted by $\mona^\bt$, is a
monoidal category for the monoidal product $\ot$.

\begin{dei}[Bimonoid]
A \emph{bimonoid} is a comonoid in the monoidal category
$(\mona^\bt, \ot, K)$.
\end{dei}

\begin{Exs}

The examples of the categories $(\KK\textrm{-Mod},\, \otk,\,
\KK)$, $(\Sy\textrm{-Mod}, \circ,\, I)$ and $(\Sy\textrm{-biMod},
\bc,\, I )$ endowed with the Hadamard tensor products, give  the
following notions.

\begin{itemize}

\item In the case of $\KK$-modules, we find the classical notion
of \emph{bialgebras}.

\item In the case of $\Sy$-modules, we find the notion of
\emph{Hopf operads}. We refer the reader to the recent preprint of
M. Aguiar and S. Mahajan \cite{AguiarMahajan} for a study of
\emph{Hopf monoids} in the category of species which is a very
close notion.

\item In the case of $\Sy$-bimodules, this generalizes the notion
of Hopf operads to properads. We call them \emph{Hopf properads}.
\end{itemize}
\end{Exs}

When $\Po$ is a Hopf operad, the category of $\Po$-algebras is
stable under the tensor product(see \ref{Alg-Op-Properad}).

\begin{pro}
Let $\Po$ be a Hopf properad. The tensor product $A\ot B$ of two
$\Po$-gebra is again a $\Po$-gebra.
\end{pro}

\begin{proo}
The proof is straightforward.
\end{proo}

\section{Koszul duality pattern}\label{Koszul Duality Pattern}

We work in the abelian monoidal category $\big(
\textrm{dg-}\Sy\textrm{-biMod},\, \bc,\, I\big)$ of
dg-$\Sy$-bimodules (see Appendix~\ref{Alg-Op-Properad}). A monoid
in this category is called a (dg-)\emph{properad}. Since the
abelian monoidal categories of differential graded vector spaces
and dg-$\Sy$-modules are abelian monoidal subcategories of
dg-$\Sy$-bimodules, the sequel includes the cases of
(dg-)associative algebras and (dg-)operads. In the following of
the text, we will implicitly work in the differential graded
context without writing ``dg'', for the sake of simplicity. We use
a very general language since most of what follows can be
generalized to another examples (colored operads, non-symmetric
operads, for instance). Denote by $\F(V)$ the free properad
(monoid) on $V$ and by $\F^c(V)$ the cofree connected coproperad
(comonoid) on $V$ (see Appendix~\ref{Alg-Op-Properad} for more
details).\\

The Koszul dual coproperad is usually defined by the top homology
groups of the bar construction. The purpose of this section is to
prove that the construction of the Koszul dual can be described
with a pure categorical or algebraic point of view. This section
is a generalization of Section $2.4$ of E. Getzler and J.D.S.
Jones preprint \cite{GJ}.

\subsection{Quadratic (co)properads}
\label{Koszuldual}

Let $(V,\, R)$ be a \emph{quadratic data}, that is $R\subset
\F_{(2)}(V)$. Since the underlying $\Sy$-bimodule of the free
properad $\F(V)$ and the cofree connected coproperad $\F^c(V)$ are
isomorphic \footnote{This should also come from the fact that the
colored operad coding properads is Koszul-autodual.}, we consider
the following sequence in $\Sy$-biMod
$$R \mono \F_{(2)}(V) \mono  \F(V)\cong\F^c(V) \epi \F^c_{(2)}(V) \epi  \F^c_{(2)}(V)/R =:\oR.$$

A quadratic data will be written $(V, R)$ or equivalently $(V,
\oR)$. To such a sequence, we can naturally define a quotient
properad of $\F(V)$ and a subcoproperad of $\F^c(V)$ (see
Appendix~\ref{generatedby}).

\begin{dei}[Quadratic properad generated by $V$ and $R$]
The \emph{quadratic properad generated by $V$ and $R$} is the
quotient properad of $\F(V)$ by the ideal generated by
$R\mono\F(V)$. We denote it by $\Po(V,\, R)=\F(V)/(R)$.
\end{dei}

\begin{dei}[Quadratic coproperad generated by $V$ and $\oR$]
The \emph{quadratic coproperad generated by $V$ and $\oR$} is the
subcoproperad of $\F^c(V)$ generated by $\F^c(V)\epi\oR$. We
denote it by $\Co(V,\, \oR)$.
\end{dei}

For example, the quadratic coalgebra generated by $(V, \oR)$ is
equal to $$ \Co(V,\oR)=\KK \oplus V \oplus \bigoplus_{n\ge 2}
\bigcap_{i=0}^{n-2} V^{\otimes i} \otimes R \otimes  V^{\otimes
n-2-i}.$$

\begin{Rq}
We proved in \cite{V1} Corollary 7.5 that when a properad is
Koszul, it is necessarily quadratic. Therefore, there is no
restriction to study only the quadratic case.
\end{Rq}

\subsection{Definition of the Koszul dual revisited}

Koszul duality theory comes from homological algebra, when one
tries to find small resolutions (minimal models) of algebraic
structures (associative algebras, operads, properads, colored operads, for instance). \\

The Koszul dual cooperad of an operad $\Po$ is defined by the top
homology of the bar construction $\B(\Po)$ (see \cite{F} Section
$5$ and \cite{GJ} Section $2.4$). In \cite{V1} Section 7, we used
the same idea to define the Koszul dual coproperad of a properad.
The purpose of this section is to prove that the Koszul dual
coproperad is a quadratic coproperad and to prove the dual statement. \\

Let $\Po=\Po(V,\, R)$ be a quadratic properad. Recall from
\cite{V1} Section 4 that the \emph{bar construction} $\B(\Po)$ of
$\Po$ is the chain complex defined on $\F^c(s\oPo)$ by the unique
coderivation $\delta$ which extends the partial composition of
$\Po$. Dually, the \emph{cobar construction} of a coproperad is
the chain complex $\F(s^{-1}\oC)$, where the differential $d$ is
the unique derivation which extends the partial composition
coproduct
of $\Co$.  \\

When $\Po=\Po(V,\, R)$ is a quadratic properad, it is weight
graded. Denote this grading by $(\omega)$. In this case, the bar
construction of $\Po$ decomposes with respect to this grading. The
part of weight $(\omega)$ of $\B_\bullet(\Po)$ begins with
$$\B_\bullet(\Po)_{(\omega)} \ : \ 0 \to \F^c_{(\omega)}(s \Po_{(1)}) \xrightarrow{\delta} \cdots . $$
Let $\Po^{\ac}_{(\omega)}$ be its top homology group
$H_\omega\left( \B_\bullet(\Po)_{(\omega)} \right) $ and
$\Po^{\ac}:=\bigoplus_{(\omega)} \Po^{\ac}_{(\omega)}$. Using
$H_\omega\left( \B_\bullet(\Po)_{(\omega)} \right)=\ker \delta$,
we proved in \cite{V1} Proposition 7.2, that $\Po^{\ac}$ is a
subcoproperad of $\F^c(s \Po_{(1)})\cong \F^c(V)$.\\

Dually, let $\Co=\Co(V,\, \overline{R})$ be a quadratic
coproperad. It is a connected coproperad, that is weight graded
and such that $\Co_{(0)}=\KK$. Once again, its cobar construction
is the direct sum of subcomplexes indexed by the weight
$$\Omega_\bullet(\Co)_{(\omega)} \ : \ \cdots \xrightarrow{d}
 \F_{(\omega)}(s^{-1}
\Co_{(1)}) \to 0. $$ Define $\Co^{\ac}$ to be the top homology
groups of the cobar construction of $\Co$, that is
$\Co^{\ac}:=\bigoplus_{(\omega)}H_{-\omega}\left(
\Omega_\bullet(\Co)_{(\omega)} \right)$. Since $H_{-\omega}\left(
\Omega_\bullet(\Co)_{(\omega)} \right)=\coker d$, $\Co^{\ac}$ is a
quotient properad of $\F(s^{-1} \Co_{(1)}) \cong \F(V)$.

\begin{thm}
Let $(V,\, R)$ be a quadratic data. Denote by $s^2R$ the image of
$R$ in $F_{(2)}(sV)$ and by $s^{-2}\oR$ the quotient of
$\F^c_{(2)}(s^{-1}V)$ by $s^{-2}R$.

The Koszul dual coproperad of $\Po(V,R)$ is equal to
$\Po(V,R)^{\ac}=\Co(sV, s^2 \oR)$. Dually, the Koszul dual
properad of $\Co(V,\oR)$ is equal to $\Co^{\ac}:=\Po(s^{-1}V,\,
s^{-2}R)$. Therefore, we have   ${\Po^{\ac}}^{\ac}=\Po$ and
${\Co^{\ac}}^{\ac}=\Co$.
\end{thm}

\begin{proo}
The cobar construction of $\Co$ has the following form
$$\Omega_\bullet(\Co)_{(\omega)} \ : \ \cdots \to
\F_{(\omega)}(s^{-1} \Co_{(1)} +\underbrace{s^{-1} \Co_{(2)}}_{1}
) \xrightarrow{d} \F_{(\omega)}(s^{-1} \Co_{(1)}) \to 0, $$ where
$\F_{(\omega)}(s^{-1} \Co_{(1)} +\underbrace{s^{-1} \Co_{(2)}}_{1}
)$ stands for the sub-$\Sy$-bimodule of $\F_{(\omega)}(s^{-1}
\Co_{(1)} +s^{-1} \Co_{(2)})$ composed by graphs with $\omega-1$
vertices indexed by elements of $s^{-1}\Co_{(1)}$ and just one
vertex indexed by an element of $s^{-1}\Co_{(2)}$. The image of
$d$ is the kernel of the cokernel $\F_{(\omega)}(s^{-1}V)\epi
\Co^{\ac}_{(\omega)}$ of $d$. Since $\Co^{\ac}$ is a quotient
properad of $\F(s^{-1}V)$, $\textrm{Im}\, d \mono \F(s^{-1}V)$ is
an ideal monomorphism. From the shape of $\Omega_\bullet(\Co)$, we
see that the image of $d$ is made of graphs indexed by $s^{-1}V$
with at least one subgraph graph in $s^{-2} R$. Therefore, the
image of $d$ is equal to the  image of $\mu^2 \, : \,
\F(s^{-1}V)\bt (\F(s^{-1}V) + \underline{s^{-2}R})  \bt
\F(s^{-1}V)$, that is the ideal generated by $s^{-2}R$ by
Proposition~\ref{Ideal-Caracterisation} of
Appendix~\ref{Alg-Op-Properad}.

We dualize the arguments (in the opposite category) to get the
dual statement. The last assertion is easily verified.
\end{proo}

A properad is called a \emph{Koszul properad} when the homology of
its bar construction is concentrated in top dimension, that is
when $H_\bullet\big( \B(\Po)\big)=\Po^{\ac}$.

\subsection{Relation with the Koszul dual properad}
\label{dualproperad}

To an $\Sy$-bimodule $M$, we associate its linear dual $M^*:= \{
M(m,n)^*\}_{m,n}$. The linear dual ${}^*$ of a coproperad $(\Co,
\Delta)$ is always a properad : define the composition product by
the formula $\Co^* \bt \Co^* \to (\Co \bt \Co)^*
\xrightarrow{\Delta^t} \Co^*$. But we need a finite dimensional
assumption on the underlying $\Sy$-bimodule to have the dual
result. The main explanation for such a phenomenon is that there
exists a map $V^*\ot W^* \to (V \ot W)^*$, which is an isomorphism
when $V$ and $W$ are finite dimensional vector spaces.

\begin{dei}[Locally finite $\Sy$-bimodule]
An $\Sy$-bimodule $M$ is \emph{locally finite} if for every $m$
and $n$ in $\NN$, the dimension of the module $M(m,n)$ is finite
over $\KK$.
\end{dei}

\begin{pro}\label{Coproperad->properad}
When $V$ is a locally finite $\Sy$-bimodule, the linear dual of
the quadratic coproperad $\Co(V, \oR)$ generated by $V$ and $\oR$
is the quadratic properad $\F(V^*)/(R^\perp)$, where $R^\perp
\subset \F_{(2)}(V)^*\cong \F_{(2)}(V^*)$.
\end{pro}

\begin{proo}
The image under ${}^*$ of the terminal object (see
Appendix~\ref{generatedby})
$$\xymatrix@M=6pt@W=6pt@H=10pt{\F^c_{(2)}(V)/R    & & \\
\F^c(V)/\Co(V,\,R) \ar@{->>}[u]& \F^c(V) \ar@{->>}[l]
\ar@{->>}[ul]& \ar@{>->}[l] \ \Co(V,\, \oR) \ar@/^/[ull]_{0}}$$
gives the initial object of

$$\xymatrix@M=6pt@W=6pt@H=10pt{\left( \F^c_{(2)}(V)/R \right)^* \ar@/_/[drr]^{0}
\ar@{>->}[dr] \ar@{>->}[d]& & \\
\left(\F^c(V)/\Co(V,\,R)\right)^* \ar@{>->}[r] & \F^c(V)^*
\ar@{->>}[r] &  \Co(V,\,R)^*.}$$

Since $V$ is locally finite, we can identify $\left(
\F^c(V)\right)^*$ with the free properad on $V^*$ : $\F(V^*)$.
(The (co)free  (co)properad on $V$ is given by a direct sum of
particular tensor powers of $V$). Therefore, $\left(
\F^c_{(2)}(V)/R \right)^*$ is isomorphic to the orthogonal of $R$,
that is $R^\perp:=\{f\in \F_{(2)}(V)^*\cong \F_{(2)}(V^*) \, | \,
f_R = 0 \}$. We conclude by the uniqueness property of the initial
object.
\end{proo}

When $\Po(V,R)$ is a quadratic properad generated by a locally
finite $\Sy$-bimodule, we consider the linear dual of the Koszul
dual coproperad ${\Po^{\ac}}^*$. By
Proposition~\ref{Coproperad->properad}, we have
${\Po(V,R)^!}^*=\Po(s^{-1}V^*, s^{-2}R^\perp)$. In the case of
finitely generated associative algebra, it is the definition given
by S. Priddy \cite{Priddy}. In the case of binary quadratic
operads, V. Ginzburg and M. Kapranov (\cite{GK} Section $2$)
defined a twisted Koszul dual operad by the formula
$\Po^!:=\Po(V^\vee, R^\perp)$, where $M^\vee(n):= M^*(n)\otimes
\textrm{sgn}_{\Sy_n}$. The reason for this lies in Quillen
functors which are the bar and cobar constructions between
$\Po$-algebras and $\Po^{\ac}$-coalgebras (see \cite{QuillenRHT}
and \cite{GJ} Section $2$). The bar construction of a
$\Po$-algebra $A$ is the cofree $\Po^{\ac}$-coalgebra on the
suspension of $A$, that is $\Po^{\ac}(s A)$. In general, we have
$$\Po^{\ac}(s A)=\bigoplus_{n \ge 1} \Po^{\ac}(n)\otimes_{\Sy_n} (s A)^{\otimes n}
=\bigoplus_{n \ge 1} s^n \Po^{\ac}(n)\otimes \textrm{sgn}_{\Sy_n}
\otimes_{\Sy_n} A^{\otimes n}. $$ We define the \emph{suspension}
operad by $S(n):=s^{n-1}.k\ot \sgn_{\Sy_n}$, with the signature
action of the symmetric group. Actually, $S$ is equal to the
operad of endomorphisms of $s^{-1}k$, that is $S=\End(s^{-1}k)$.
We have $\Po^{\ac}(sA)=s(S\otimes\Po^{\ac})(A)$. Up to
suspensions, $\Po^{\ac}(sA)$ is the cofree ``$\Po^!$-coalgebra''
on $A$. The operad $\Po$ is Koszul if and only if $\Po^{\ac}$ is
Koszul, which is also equivalent to $\Po^!$ is Koszul.

\section{Manin products}
\label{ManinProducts} The aim of this section is to provide a
general and intrinsic framework for the definitions of Manin's
black and white products. We first give the conceptual definition
of Manin's white product of monoids in any lax $2$-monoidal
category. Then, we dualize the arguments to get the notion of
black product of comonoids in any colax $2$-monoidal
category.\\

We make explicit all the constructions in the category of
$\Sy$-bimodules. But they remain valid in general $2$-monoidal
categories with mild assumptions (existence of the free monoid,
cofree comonoid, for instance). These constructions also hold for
non-symmetric operads (see Section \ref{BW Square}) and colored
operads, for instance. We denote the vertical connected
composition product of $\Sy$-bimodules and the Hadamard horizontal
tensor $\ot_H$ by $\ot$, to lighten the notations.

\subsection{A canonical map between free
monoids}\label{Phi}

V. Ginzburg and M.M. Kapranov mentioned  in \cite{GK2} a morphism
of operads $\Phi \, : \, \F(V \ot W) \to \F(V)\ot \F(W)$ ``which
reflects the fact that the tensor product of an $\F(V)$-algebra
and an $\F(W)$-algebra is an $\F(V\ot W)$-algebra''. We describe
and extend this map $\Phi$ to a more general setting.

\begin{pro}
Let $(\CA,\, \boxtimes,\, I,\, \ot , K)$ be a lax $2$-monoidal
category such that $(\CA,\, \boxtimes,\, I)$ admits free monoids.
There exists a natural morphism of monoids $\Phi \, : \F(V\ot W)
\to \F(V) \ot \F(W)$.
\end{pro}

\begin{proo}
Let $V$ and $W$ be two objects in $\CA$. There is a natural map
$u_{\F(V)} \ot u_{\F(W)}\, : \, V\ot W \to \F(V)\ot \F(W)$. Using
Proposition~\ref{tenseur de monoides}, we know that $\F(V)\ot
\F(W)$ is a monoid for $\bt$. By the universal property of the
free monoid on $V\ot W$, there exists a unique morphism of monoids
$\Phi$ which factors the previous map
$$\xymatrix{V \ot W \ar[dr]_{u_v\ot u_W}  \ar[r] ^{u_{V\ot W}} &
  \F(V\ot W) \ar@{-->}[d]^{\exists ! \Phi} \\
 & \F(V)\ot \F(W).} $$
\end{proo}

\begin{Exs}
$ $

\begin{itemize}

\item When $\CA$ is the category of $k$-modules, the map $\Phi$ is
the direct sum of the isomorphisms $(V\ot W)^{\ot n} \cong V^{\ot
n} \ot W^{\ot n}$ induced by the twisting map.

\item In the category of $\Sy$-modules, the map $\Phi$ corresponds
to the injective morphism of operads $\F(V\ot W) \mono \F(V)\ot
\F(W)$ mentioned in  \cite{GK2}.

\item For $\Sy$-bimodules, the previous construction gives a
morphism of properads between the free properad $\F(V\ot W)$and
the Hadamard product $\F(V)\ot\F(W)$. Once again, this map is
always injective but not an isomorphism in general.
\end{itemize}
\end{Exs}

One remark before to conclude this section. The purpose of this
paragraph was to show that the definition of the map $\Phi$ is
canonical and does not depend on the bases of the modules involved
here. Now, if we choose a basis for the free operad, for instance,
we can make the map $\Phi$ more explicit. In this case, the image
of a tree $T$ with vertices indexed by elements of $V\ot W$ under
$\Phi$ is the tensor product of the same tree $T$ with vertices
indexed by the corresponding elements of $V$ with the tree $T$
whose vertices are indexed by the corresponding elements of $W$.

\subsection{Definition of the white product}

In this section we define the white product for every pair of
properads defined by generators and relations. When the two
properads are quadratic, the resulting white product is again
quadratic. Since an associative algebra is an operad and an operad
is a properad, this construction summarizes what can be found in
the literature. In the case of quadratic associative algebras, it
 corresponds to the original notions introduced
by Yu. I. Manin \cite{M} and in the case of binary quadratic
operads,
it corresponds to the definitions of V. Ginzburg and M. Kapranov \cite{GK, GK2}.\\

The properties of the morphism $\Phi$ lead directly to the
definition of the white product. Let $\Po$ and $\Qo$ be two
properads defined by generators and relations, $\Po=\F(V)/(R)$ and
$\Qo=\F(W)/(S)$. And denote the projections $ \pi_\Po \, : \,
\F(V) \epi \Po$ and $ \pi_\Qo \, : \, \F(W) \epi \Qo$.

Consider the following composite of morphisms of properads
$$\pi_\Po
\ot \pi_\Qo \, \circ \, \Phi \,  : \,  \F(V \ot W) \mono \F(V)\ot
\F(W) \epi \Po \ot \Qo.$$ Since it is a morphism of properads, its
kernel is an ideal of $\F(V\ot W)$. It is the ideal generated by
$\Phi^{-1}(R\ot \F(W) + \F(V)\ot S)$ in $\F(V\ot W)$.

\begin{dei}[White product]
Let $\Po=\F(V)/(R)$ and $\Qo=\F(W)/(S)$ be two properads defined
by generators and relations. The quotient properad
$$\Po \circ \Qo:= \F(V\ot W)/\big( \Phi^{-1}(R\ot \F(W) + \F(V)\ot
S)\big)$$ is called the \emph{white product of} $\Po$ and $\Qo$.
\end{dei}

The definition of the white product of two properads shows that
the morphism $\Phi$ factors through a natural morphism of
properads $\bar{\Phi} \, : \, \Po \circ \Qo \to \Po \ot \Qo$. In
the abelian category $\Sy$-bimodules, $\bar{\Phi}$ is the image of
$\pi_\Po \ot \pi_\Qo \, \circ \, \Phi$. Hence, it is a
monomorphism.
$$\xymatrix@M=6pt@W=6pt@H=10pt@C=40pt{\F(V\ot W)\ \ar@{>->}[r]^(0.45){\Phi} \ar@{->>}[drr]& \F(V) \ot \F(W)\
\ar@{->>}[r]^(0.6){\pi_\Po \ot \pi_\Qo} & \Po \ot \Qo    \\
& & \Po \circ \Qo \ar@{>->}[u]_{ \bar{\Phi}}} $$

Let $A$ be a $\Po$-gebra and $B$ a $\Qo$-gebra, since the tensor
$A\ot B$ is a $\Po \ot \Qo\,$-gebra, we get the following result.

\begin{pro}
\label{tensor-white} The tensor product $A\ot B$ is a gebra over
the white product $\Po \circ \Qo$.
\end{pro}

\begin{Ex}\label{Exemples Hom(C,A)}
Let $\Po$ and $\Qo$ be two operads. The tensor product of a
$\Po$-algebra with a $\Qo$-algebra is a $\Po\ot \Qo$-algebra. We
can partially dualize this statement. Let $C$ be a $\Qo$-coalgebra
and $A$ be a $\Po$-algebra, the space of morphism $\Hom_\KK(C,A)$
is $\Po \ot \Qo$-algebra (see \cite{BergerMoerdijk03} Proposition
1.1). It is also a $\Po \circ \Qo$-algebra by
Proposition~\ref{tensor-white}. As explained by G.~Barnich,
R.~Fulp, T.~Lada, and J.~Stasheff in
\cite{BarnichFulpLadaStasheff00}, when $C$ is a cocommutative
coalgebra, $\Hom_\KK(C,A)$ is always a $\Po$-algebra. This comes
from the fact that $\C$ is the unit object for $\ot$ and $\circ$.
Motivated by structures appearing in Lagrangian field theories in
physics, these authors studied the algebraic structures of
$\Hom_\KK(C,A)$ when $C$ is a coassociative coalgebra and $A$ a
Lie algebra or a Poisson algebra. Since $\Hom_\KK(C,A)$ is a
$\Po\circ \Qo$-algebra, Manin's white product for operads gives a
way to describe such structures.
\end{Ex}

The white product is a construction that preserves the grading of
the properads.

\begin{pro}
If $S\subset \bigoplus_{\omega =0 }^{N} \F_{(\omega)}(V)$ and
$R\subset \bigoplus_{\omega =0 }^{M} \F_{(\omega)}(W)$, the white
product of $\Po$ and $\Qo$ is a properad generated by $V\ot W$
with relations in $\bigoplus_{\omega =0 }^{\max(N,\, M)}
\F_{(\omega)}(V\ot W)$.

If $S$ and $R$ are homogenous of weight $N$, that is $S,\, R
\subset \F_{(N)}(V)$, then $\Po\circ \Qo$ is once again a properad
defined by homogenous relations of weight $N$.
\end{pro}

\begin{proo}
It comes from the definition of the morphism $\Phi$ which
preserves the grading.
\end{proo}

\begin{Exs} $ $

\begin{itemize}
\item Let $A$ and $B$ be two quadratic associative algebras. The
white product $A \circ B$ is equal to $T(V\ot W)/\big((23)(R\ot
W^{\ot 2} + V^{\ot 2} \ot S) \big)$, which is the definition given
by Manin in \cite{M,M2}. It is isomorphic via $\bar{\Phi}$ to the
Hadamard (or Segre) product $A \ot_H B:=\bigoplus_n A_{(n)} \ot_k
B_{(n)}$. This crucial property allowed J. Backelin to prove in
his thesis \cite{B} that the white product of two Koszul algebras
is a Koszul algebra.

\item An associative algebra $A=T(V)/(R)$ is \emph{$N$-homogenous}
if $R\subset V^{\ot N}$. R. Berger, M. Dubois-Violette and M.
Wambst generalized Manin's black and white products to
$N$-homogenous algebras in \cite{BergerDuboisWambst}. Berger and
Marconnet proved that the black product of two N-homogenous Koszul
algebra is still Koszul under some extra assumptions
(distributivity) in \cite{BergerMarconnet} Proposition 2.8. For
two $N$-homogenous algebras, the definition given above coincide
with their definition. Note that the definition given here can be
applied to non-homogenous algebras. The class of Artin-Schelter
algebras \cite{ArtinSchelter} provide interesting examples of
non-homogenous algebras. It would be interesting to study the
properties of the white product of such algebras, for instance the
ones of global dimension $4$ of \cite{Palmieri}.

\item When $\Po$ and $\Qo$ are binary quadratic operads, the
modules $\F_{(2)}(V)$ and $\F_{(2)}(W)$ are equal to $\F(V)(3)$
and $\F(W)(3)$. In that case, we get $R\ot \F(W) = R \ot \F(W)(3)$
and $\F(V)\ot S = \F(V)(3) \ot S$. This construction is the
original one described by Ginzburg and Kapranov in \cite{GK, GK2}.
Note that in this case, the white product is not, in general,
equal to the Hadamard product. (The morphism $\bar{\Phi}$ is not
an isomorphism in general). A direct consequence of this fact is
that the white product of two Koszul operads is not necessarily a
Koszul operad again. See section~\ref{Counterexample} for a
counterexample.
\end{itemize}

\end{Exs}

\subsection{The black product}
\label{BlackDefinition}

We dualize the arguments and work in the opposite category. This
gives the definition \emph{black product} of coproperads.

\begin{pro}
Let $(\CA,\, \boxtimes,\, I,\, \ot , K)$ be a colax $2$-monoidal
category such that $(\CA,\, \boxtimes,\, I)$ admits cofree
comonoids. There exists a natural morphism of comonoids $\Psi \, :
\F^c(V\ot W) \gets \F^c(V) \ot \F^c(W)$.
\end{pro}

\begin{dei}[Black product]
The black product of two coproperads $\Co(V, \oR)$ and $\Co(W,
\overline{S})$ is the image of the morphism of comonoids $\Psi
\circ (\iota\otimes \iota)$
$$\xymatrix@M=6pt@W=6pt@H=10pt@C=40pt{\F^c(V\ot W)\  & \ar@{->>}[l]_{\Psi} \F^c(V) \ot \F^c(W)\
 & \ar@{>->}[l]_{\iota \ot \iota}{\ \Co(V, \oR)\otimes
\Co(W,\overline{S})}   \ar@{->>}[d]^{\bar{\Psi}}  \\
& & \ar@{>->}[ull] \Co(V, \oR)\bullet \Co(W,\overline{S}). }$$ It
is equal to $\Co(V, \oR)\bullet \Co(W,\overline{S})=\Co(V\ot W,
\overline{\Psi(R\ot S)})$.
\end{dei}

Black and white constructions are dual to each other under linear
duality.

\begin{thm}
\label{BlackDualWhite} Let $(V,\, R)$ and $(W,\, S)$ be two
quadratic data, with $V$ and $W$ locally finite. We have the
following isomorphism of properads $\left(\Co(V,\, R) \bullet
\Co(W,\, S) \right)^* \cong \Co(V,\,R)^* \circ \Co(W,\,S)^*$.
\end{thm}

\begin{proo}
Since $\psi$ is the transpose of $\varphi$, we have
$\Psi={}^t\Phi_{V^*,W^*}$, up to isomorphism like $\big(\F(V^*)\ot
\F(W^*)\big)^*\cong \F(V) \ot \F(W)$. Therefore, we get $\Psi(R\ot
S)^\perp=\Phi^{-1}_{V^*, W^*}(R^\perp \ot \F(W^*) + \F(V^*)\ot
S^\perp)$. By Proposition~\ref{Coproperad->properad}, we have
\begin{eqnarray*}
\left(\Co(V,\, R) \bullet \Co(W,\, S) \right)^* &\cong& \Po(V^*\ot
W^*, \Psi(R\ot S)^\perp)\\
&\cong& \Po(V^*\ot W^*, \Phi^{-1}_{V^*, W^*}(R^\perp \ot \F(W^*) +
\F(V^*)\ot S^\perp))\\
&\cong& \Po(V^*, R^\perp) \circ \Po(W^*, S^\perp)\\
&\cong& \Co(V,\,R)^* \circ \Co(W,\,S)^*.
\end{eqnarray*}
\end{proo}

One of the main interest of the classical notions of black and
white products is that one gives the other via the Koszul dual
functor. In the next sections, we define a black product for
monoids (operad and non-symmetric operads). The translation of
Theorem~\ref{BlackDualWhite} in this framework will give the
relation with Koszul dual functor.

\section{Manin products for operads}\label{Manin products for operads}

In this section, we study Manin products for (symmetric) operads.
We first give a sufficient condition for the white product to be
equal to the Hadamard product. Then, we recall the bases used to
describe binary quadratic operads and their Koszul dual operad. We
refer the reader to \cite{GK}, \cite{L} and \cite{MSS} for
complete references. We make our constructions explicit for binary
quadratic operads in order to do computations. The linear dual
version of the black product for cooperads defines a product for
operads which corresponds to the definition of Ginzburg and
Kapranov, where we make the signs precise. We give an example of a
pair of Koszul operads such that their products is not Koszul.
This shows that black and white products for operads do not behave
like black and white products for associative algebras. Following
Yu. I. Manin, we prove that $\Po \bullet \Po^!$ is always a Hopf
operad. Finally, we describe the relation between unary operators
and black products.

\subsection{Relation between the Hadamard  product and the white product}
\label{Hadamard=White} We saw in the previous section that the
composite $(\pi_\Po \ot \pi_\Qo) \circ \Phi$ factors through its
image $\bar{\Phi} \, : \, \Po \circ \Qo \mono \Po\ot \Qo$.
Therefore, $\bar{\Phi}$ is an isomorphism if and only if the
composite $(\pi_\Po \ot \pi_\Qo) \circ \Phi$ is an epimorphism. We
shall give a sufficient
condition for this.\\

Consider the case of binary quadratic operads, that is quadratic
operads generated by binary operations ($V(n)=0$ for $n\neq 2$).
In this case, the free operad on $V$ is given by (non-planar)
binary trees with vertices labelled by operations of $V$. Denote
by $\mathbb{T}$ such a tree with $n-1$ vertices and the induced
\emph{label morphism} by $\mathcal{L}_\mathbb{T}^V \, :\, V^{\ot
(n-1)}\to \F(V)(n)$.

\begin{pro}
Let $\Po$ be a binary quadratic operad such that for every $n\ge
3$ and every binary tree $\mathbb{T}$ with $n-1$ vertices, the
composite $\pi_\Po \circ \mathcal{L}_\mathbb{T}^V \, : \, V^{\ot
(n-1)} \to \F(V)(n) \epi \Po(n)$ is surjective.

For every binary quadratic operad $\Qo$, the white product
$\Po\circ \Qo$ is equal to the Hadamard product $\Po \ot \Qo$.
\end{pro}

\begin{proo}
It is enough to prove that $(\pi_\Po \ot \pi_\Qo) \circ \Phi$ is
an epimorphism. Let $p\ot q$ be an elementary tensor of $\Po(n)\ot
\Qo(n)$, where $\Qo=\F(W)/(S)$. The element $q$ of $\Qo(n)$ can be
written $q=\sum_{i=1}^k \pi_\Qo \circ
\mathcal{L}_{\mathbb{T}_i}^W(w_1^i,\, \ldots ,\,  w_{n-1}^i)$,
with $\{\mathbb{T}_i\}$ a finite set of trees and $\{w^i_j\}$
elements of $W(2)$. By the assumption, there exists $v^i_1,\,
\ldots ,\, v^i_{n-1}$ in $V(2)$ such that $p=\pi_\Po\circ
\mathcal{L}_{\mathbb{T}_i}^V(v^i_1,\, \ldots ,\, v^i_{n-1})$, for
every $\mathbb{T}_i$. Therefore, we have $\displaystyle
p=\frac{1}{k} \sum_{i=1}^k \pi_\Po\circ
\mathcal{L}_{\mathbb{T}_i}^V(v^i_1 ,\, \ldots ,\, v^i_{n-1})$.
Finally, it shows that
$$p\ot q = (\pi_\Po \ot \pi_\Qo) \circ \Phi \left(
\frac{1}{k} \sum_{i=1}^k \mathcal{L}_{\mathbb{T}_i}^{V\ot
W}(v_1^i\ot w_1^i,\, \ldots ,\, v_{n-1}^i\ot w_{n-1}^i) \right)
.$$
\end{proo}

The condition of this proposition means that every operation of
$\Po$ can be written by any type of composition of generating
operations. In the next corollary, we show that the operads $\C$,
$\Pe$ and $\Comtri$ are examples of such operads. Recall briefly
that $\C$ is the operad for commutative algebras. The operad $\Pe$
was introduced by F. Chapoton in \cite{C} and $\Comtri$ was
defined in \cite{V2} Appendix A.

\begin{cor}\label{Hadamard=tensor}
For every binary quadratic operad $\Qo$, we have
\begin{itemize}
\item $\C \circ \Qo =\C \ot \Qo =\Qo$. The operad $\C$ is neutral
for the white product in the category of binary quadratic operads.

\item $\Pe \circ \Qo= \Pe \ot \Qo$ and $\Comtri \circ \Qo= \Comtri
\ot \Qo$
\end{itemize}
\end{cor}

\begin{proo}
The operad $\C$ is generated by $V(2)=k$ with trivial action of
$\Sy_2$ and the associativity relation. Hence, we have only one
commutative operation with arity $n$, that is $\C(n)=k$.
Therefore, for every tree $\mathbb{T}$, the morphism
$\mathcal{L}_{\mathbb{T}}$ is a surjection on $k$ and $\C\circ
\Qo=\C \ot \Qo=\Qo$.

The operad $\Pe$ corresponds to commutative operations with one
input emphasized (see \cite{V2} 4.2 and \cite{CV} 1.3.2). In arity
$n$, we have $n$ operations $\Pe(n)=k.e_1^n \oplus \cdots k.e_n^n$
where $e^n_i$ corresponds to the corolla with $n$ inputs such that
the $i$th input (or branch) is emphasized. The composition of
corollas gives a corolla where the leaf emphasized is the one with
a path to the root via emphasized branches.
$$\vcenter{\xymatrix@C=12pt@R=12pt{1 \ar@{-}[dr] & 2 \ar@{=}[d] & 3 \ar@{-}[dl] &
& 4 \ar@{-}[dr] & 5 \ar@{=}[d]  & 6 \ar@{-}[dl] \\
 &  *{} & & & & *{}  & \\
 & & &  *{}   \ar@{=}[rru]  \ar@{-}[llu] & & & \\
 & & &\ar@{-}[u] & & & }}\ \mapsto \ e^6_5 $$
Let $\mathbb{T}$ be a binary tree with $n-1$ vertices. To get
$e^n_i$, it is enough the look at the unique path from the $i$th
leaf to the root and index the vertices on this path with the
relevant operations.
$$\vcenter{\xymatrix@M=1pt@C=12pt@R=12pt{1 &2 \ar@{=}[dr] &3 &4  \ar@{-}[dl] \\
\ar@{-}[dr]& &e^2_1 \ar@{=}[dl] & \\& \ar@{=}[dr] e^2_2 & &
\ar@{-}[dl]\\& &
 e^2_1 \ar@{-}[d] &
\\& & & }}\  \mapsto \ e^4_2 $$
The operations of $\Comtri(n)$ are corollas with at least one leaf
emphasized and the proof is the same.
\end{proo}

This corollary shows that the Hadamard product of one operad $\C$,
$\Pe$ or $\Comtri$ with any other binary quadratic operad is again
a binary quadratic operad. For $\C$, the result is obvious. In the
particular case of $\Pe$ this result was proved directly by F.
Chapoton in \cite{C}. For every binary quadratic operad $\Qo$, he
constructed by hand a quadratic operad isomorphic to $\Pe \ot
\Qo$. This construction is actually the white product $\Pe \circ
\Qo$.

\begin{pro}
\label{PermWhiteAs} We have $\Pe \circ \A=\D$.
\end{pro}

\begin{proo}
Using the complete description of $\Pe$, $\A$ and $\D$, Chapoton
proved in \cite{C} that $\Pe \ot \A = \D$. Apply the previous
corollary to conclude.
\end{proo}


\subsection{Binary quadratic operad and Koszul dual operad}
\label{BinaryQuadraticOperad} The preceding section gives a method
for computing the white product for a particular class of operads.
When we cannot apply this method, we need the explicit form of the
products to compute them. In this section, we describe a
 basis for binary quadratic operads and their Koszul dual operads.\\

Recall that the free operad $\F(V)$ on $V$ is given by trees with
the vertices indexed by elements of $V$, with respect to the
action of the symmetric groups. When $V$ is an $\Sy_2$-module,
that is a module over the symmetric group $\Sy_2$, we have
$\F_{(2)}(V)=\F(V)(3)$, the part with $3$ inputs of the free
operad on $V$ which is isomorphic to
$$\F(V)(3)=\big(V\ot_{\Sy_2}(V\ot k \oplus k \ot V)\big) \ot_{\Sy_2} k[\Sy_3],$$
where the summand $V\ot(V \ot k)$ corresponds to the compositions
on the left $\vcenter{\xymatrix@M=0pt@R=5pt@C=5pt{ \ar@{-}[dr] &
&\ar@{-}[dl] & & \ar@{-}[dl]  \\
& \ar@{-}[dr] & &\ar@{-}[dl]  & \\
& &\ar@{-}[d] & & \\
& & \\ & & }} $ and the summand $V\ot(k \ot V)$ corresponds to the
compositions on the right $\vcenter{\xymatrix@M=0pt@R=5pt@C=5pt{
\ar@{-}[ddrr] & &\ar@{-}[dr] & & \ar@{-}[dl]  \\
&  & &\ar@{-}[dl]  & \\
& &\ar@{-}[d] & & \\
& & \\ & & }} $. Since the action of $\Sy_2$ maps one to the
other, we choose the one on the left and $\F(V)(3)$ is isomorphic
to the induced representation $\mathrm{Ind}_{\Sy_2\times
\Sy_1}^{\Sy_3} \left( V\ot (V\ot k)\right)$. Therefore, $\F(V)(3)$
can be identify with $3$ copies of $V\ot V$ represented by the
following types of tree :

$$\textrm{I :}  \vcenter{\xymatrix@R=10pt@C=10pt{x \ar@{-}[dr] & & \ar@{-}[dl] y  && \ar@{-}[ddll] z \\
& *+[F-,]{\nu} \ar@{-}[dr] &  && \\
&& *+[F-,]{\mu} \ar@{-}[d] && \\
&&&&}} \quad \textrm{II :} \vcenter{\xymatrix@R=10pt@C=10pt{y \ar@{-}[dr] & & \ar@{-}[dl] z  && \ar@{-}[ddll] x \\
& *+[F-,]{\nu} \ar@{-}[dr] &  && \\
&& *+[F-,]{\mu} \ar@{-}[d] && \\
&&&& }} \quad \textrm{III :}
\vcenter{\xymatrix@R=10pt@C=10pt{z \ar@{-}[dr] & & \ar@{-}[dl] x  && \ar@{-}[ddll] y \\
& *+[F-,]{\nu} \ar@{-}[dr] &  && \\
&& *+[F-,]{\mu} \ar@{-}[d] && \\
&&&& }}.$$ Denote them by $\mu\circ_\textrm{I}\nu$,
$\mu\circ_{\textrm{II}}\nu$ and $\mu\circ_{\textrm{III}}\nu$.

The action of the permutation $(12)$ is given by $\left(
\mu\circ_\textrm{I} \nu \right)^{(12)}= \mu\circ_\textrm{I}
\nu^{(12)}$, $\left( \mu\circ_\textrm{II} \nu \right)^{(12)}=
\mu\circ_\textrm{III} \nu^{(12)}$, $\left( \mu\circ_\textrm{III}
\nu \right)^{(12)}= \mu\circ_\textrm{II} \nu^{(12)}$  and the
action of $(132)$ is given by $\left( \mu\circ_\alpha \nu
\right)^{(123)}= \mu\circ_{(\alpha + \textrm{I})} \nu$.

\begin{Rq}
This basis is different from the one in \cite{GK} p. 228. The one
given here has nice symmetric properties with respect to the
action of $\Sy_3$ that we will use in \ref{Examples} to simplify
the computations.
\end{Rq}

The dual representation $V^*$ of an $\Sy_n$-module $V$ is the
vector space $V^*=\Hom(V,\, k)$ endowed with the following right
action of the symmetric group. For $f \, : \, V \to k$ and $\sigma
\in \Sy_n$, we have $\left( f^\sigma \right)
(x):=f(x^{\sigma^{-1}})$. We will need to twist the dual
representation by the signature, that is $V^\vee:=V^*\ot \sgn_{\Sy_n}$.\\

Let $V$ be an $\Sy$-module concentrated in arity $2$, that is an
$\Sy_2$-module. When $V$ is a finite dimensional $k$-vector space,
denote by $\mu,\, \nu, \eta,\, \zeta, \ldots $ one of its basis,
stable by the action of $\Sy_2$, and by $\mu^*,\, \nu^*, \eta^*,\,
\zeta^*, \ldots $ the dual basis. Therefore $\mu^\vee=\mu^*,\,
\nu^\vee=\nu^*, \eta^\vee=\eta^*,\, \zeta^\vee=\zeta^*, \ldots $
forms a basis of $V^\vee$ such that
$(\mu^\vee)^{(12)}=-(\mu^{(12)})^\vee$.We define the following
non-degenerate bilinear form
\begin{eqnarray*}
&& \F(V)(3)\ot \F(V^\vee)(3) \xrightarrow{<\,,\,>} k \\
&& <\mu\circ_\alpha \nu,\, \eta^\vee\circ_\beta
\zeta^\vee>:=\left\{\begin{array}{ll} 1 & \textrm{if}\
\alpha=\beta,
\,\mu=\eta \ \textrm{and} \ \nu=\zeta, \\
0 & \textrm{otherwise}. \end{array} \right.
\end{eqnarray*}

For a sub-$\Sy_3$-module $R$ of $\F(V)(3)$, we consider its
orthogonal $R^\perp:=\{\Omega \in\F(V^\vee)(3)\  |\ <\omega
,\,\Omega >=0, \ \forall\  \omega \in R \}$ for this bilinear
form.

Since the action of $\Sy_3$ on the bilinear form $<\, ,\,>$ is
given by the signature $<\omega^\sigma,\,
\Omega^\sigma>=\sgn(\sigma).<\omega,\, \Omega>$ we have that
$R^\perp$ is a sub-$\Sy_3$-module of $\F(V^\vee)(3)$. Note that
the non-degenerate bilinear form $<\, ,\,>$ defines an isomorphism
of $\Sy_3$-modules from $\F(V)(3)^\vee$ to $\F(V^\vee)(3)$.\\

Recall from~\ref{dualproperad} that under finite dimensional
assumptions, the Koszul dual operad of $\F(V)/(R)$ is
$\Po^!=\F(V^\vee)/(R^\perp)$. This bilinear form provides a method
for computing it. The canonical isomorphism $(V^\vee)^\vee \cong
V$ induces $(R^\perp)^\perp\cong R$ and
$\left(\Po^!\right)^!=\Po$.

\begin{Exs}
The operad $\C$ for commutative (associative) algebras $(A,\, *)$
is generated by the one dimensional $\Sy_2$-module $V:=k.*$ with
trivial action. Denote by $t_1=*\circ_\textrm{I}*$,
$t_2=*\circ_\textrm{II}*$ and $t_3=*\circ_\textrm{III}*$ the
elements of the basis of $\F(V)(3)$. The associativity relation is
the quadratic relation $t_1=t_2=t_3$. Therefore, the operad $\C$
has the following presentation $\C=\F(k.*)/(t_1-t_2,\, t_2-t_3)$.\\

The operad $\Li$ for Lie algebras $(L,\, [\,,\,])$ is generated by
the one dimensional $\Sy_2$-module $V':=k.{[\,,\,]}$ where the
action is given by the signature. If we denote by $t'_1$, $t'_2$
and $t'_3$ the elements of the basis of $\F(V')(3)$, the Jacobi
relation corresponds to $t'_1+t'_2+t'_3=0$ and the operad $\Li$ is
given by $\Li=\F(k.{[\,,\,]})/(t'_1+t'_2+t'_3)$.\\

Under the identification $V'\cong V^\vee$, we have
$\left((t_1-t_2).k \oplus (t_2-t_3).k
\right)^\perp=(t'_1+t'_2+t'_3).k$. Therefore we get $\C^!=\Li$
(and $\Li^!=\C$).
\end{Exs}

\subsection{Definition of the black product for
operads}\label{BlackproductOperadsSection}

Using the notions of the previous section, we define a black
product for binary quadratic operads.\\

The definition of the white product is based on the morphism
$\Phi$ (see \ref{Phi}). For binary quadratic operads, this
morphism  $\Phi\, : \, \F(V\ot W)(3) \to \F(V)(3)\ot\F(W)(3)$ is
the componentwise projection. For instance, for compositions of
type I, we have
$$  \vcenter{\xymatrix@R=8pt@C=8pt{x \ar@{-}[dr] & & \ar@{-}[dl] y  && \ar@{-}[ddll] z \\
& *+[F-,]{\nu_1\ot \nu_2} \ar@{-}[dr] &  && \\
&& *+[F-,]{\mu_1\ot \mu_2} \ar@{-}[d] && \\
&&&&}}\mapsto \vcenter{\xymatrix@R=8pt@C=8pt{x \ar@{-}[dr] & & \ar@{-}[dl] y  && \ar@{-}[ddll] z \\
& *+[F-,]{\nu_1} \ar@{-}[dr] &  && \\
&& *+[F-,]{\mu_1} \ar@{-}[d] && \\
&&&&}}\ot \vcenter{\xymatrix@R=8pt@C=8pt{x \ar@{-}[dr] & & \ar@{-}[dl] y  && \ar@{-}[ddll] z \\
& *+[F-,]{\nu_2} \ar@{-}[dr] &  && \\
&& *+[F-,]{\mu_2} \ar@{-}[d] && \\
&&&&}}.$$

We describe a general method that will be applied later in other cases. \\

When $V$ is finite dimensional, the Koszul dual of binary
quadratic operad $\F(V)/(R)$ can be defined by means of a
particular non-degenerate bilinear form on $\F(V)(3)\ot
\F(V^\vee)(3)$ (see \ref{BinaryQuadraticOperad}) denoted by
$<,>_{V}$. For the moment, we do not need its explicit
description. Since this bilinear form is non-degenerate, it
induces an isomorphism $\theta_V \, : \, \F(V)(3)
\xrightarrow{\simeq} \F(V^\vee)(3)^\vee$. Let $V$ and $W$ be two
finite dimensional $k$-modules. Define the morphism $\Psi$ by the
following commutative diagram
$$\xymatrix@C=40pt{\F(V)(3)\ot \F(W)(3) \ot k.\sgn_{\Sy_3} \ar[r]^(.55){\Psi}
\ar[d]^{\theta_{V} \ot \theta_{W}  \ot \textrm{sgn}} & \F(V\ot W\ot k.\sgn_{\Sy_2})(3) \\
\F(V^\vee)(3)^\vee \ot \F(W^\vee)(3)^\vee \ot k.\sgn_{\Sy_3}
\ar[d]^(.45){\simeq} & \ar[u]_{\theta^{-1}_{V\ot W \ot
\textrm{sgn}}}
\F((V\ot W\ot k.\sgn_{\Sy_2} )^\vee)(3)^\vee  \\
\big( \F(V^\vee)(3)\ot \F(W^\vee)(3) \big)^\vee
\ar[r]^(.57){{}^t\Phi_{V^\vee,\, W^\vee}} & \ar[u]_{\simeq}
\F(V^\vee \ot W^\vee)(3)^\vee,}$$ where $\simeq$ stands for the
natural isomorphism for the linear dual of a tensor product, since
the modules are finite dimensional. The morphism $\Psi$ defined
here is a twisted version of the one defined in
\ref{BlackDefinition}.

Recall that $\Phi_{V^\vee,\, W^\vee}$ is the morphism
$\F(V^\vee\ot W^\vee) \to \F(V^\vee)\ot \F(W^\vee)$.

\begin{lem}
\label{PhiPsiOrthogonal} Let $\Po=\F(V)/(R)$ and $\Qo=\F(W)/(S)$
be two binary quadratic operads such that $V$ and $W$ are finite
dimensional. The orthogonal of $\Psi(R\ot S)$ for $<\, ,\,
>_{V\ot W\ot \textrm{sgn}}$ is $\Phi_{V^\vee,\,
W^\vee}^{-1}(R^\perp \ot \F(W^\vee) + \F(V^\vee)\ot S^\perp)$.
\end{lem}

\begin{proo}
By definition of the transpose of $\Phi_{V^\vee,\, W^\vee}$, we
have
\begin{eqnarray*}
<\Psi(  r\ot s),\, X>_{V\ot W\ot \textrm{sgn}} &=&  <r\ot s,
\Phi_{V^\vee,\, W^\vee}(X)>_{(\F(V)\ot \F(V^\vee))\times
(\F(W)\ot \F(W^\vee))} \\
&=& \left( <r, - >_{V}.<s,- >_{W}\right) \circ \Phi_{V^\vee,\,
W^\vee}(X),
\end{eqnarray*}
for every $(r, s) \in R\times S$ and every $X\in \F\big((V \ot W
\ot k.\sgn_{\Sy_2})^\vee\big)$

Therefore, we have $\Psi(R\ot S)^\perp=$
\begin{eqnarray*}
&=& \left\{ X\in \F((V\ot W \ot k.\sgn_{\Sy_2})^\vee)(3) \ |\
\forall (r,\, s)\in R\times S \quad <\Psi(r\ot s),\, X>_{V\ot
W \ot k.\sgn_{\Sy_2}}=0\right\}\\
 &=& \left\{ X \in
\F(V^\vee \ot W^\vee)(3) \ | \ \forall (r,\, s)\in R\times S \quad
(<r,
->_{V'}.<s, ->_{W'}) \circ \Phi_{V^\vee, W^\vee }(X)=0\right\}\\
&=& \left\{ X \in \F(V^\vee \ot W^\vee)(3) \ | \ \Phi_{V^\vee,
W^\vee }(X) \in
R^\perp \ot\F(W^\vee) + \F(V^\vee)\ot S^\perp  \right\} \\
&=& \Phi_{V^\vee, W^\vee}^{-1}(R^\perp \ot\F(W^\vee) +
\F(V^\vee)\ot S^\perp ).
\end{eqnarray*}
\end{proo}

\begin{dei}[Black product for operads]
Let $\Po=\F(V)/(R)$ and $\Qo=\F(W)/(S)$ be two binary quadratic
operads with finite dimensional generating spaces. Define their
black product by the formula
\begin{eqnarray*}
\Po\bullet \Qo &=&  \F(V\ot W \ot  k.\sgn_{\Sy_2})/\big( \Psi(R\ot
S ) \big).
\end{eqnarray*}
\end{dei}

\begin{pro}
\label{BlackOperad} For binary quadratic operads generated by
finite dimensional $\Sy_2$-modules, this definition of black
product verifies $\big(\Po \bullet \Qo\big)^! = \Po^! \circ \Qo^!
$ and corresponds to the one of Ginzburg and Kapranov \cite{GK2}.
\end{pro}

\begin{proo}
It is a direct corollary of Lemma~\ref{PhiPsiOrthogonal}.
\end{proo}

Since $\C$ is the neutral element for $\circ$, we have that $\Li$
is the neutral element for $\bullet$.

\subsection{Examples}\label{Examples}

We make explicit some computations of black and white products.
For the definitions of the various operads encountered in this
section, we
refer the reader to \cite{L3}. \\

In order to compute black and white products for operads where the
space of generators $V$ is equal to $k[\Sy_2]=\mu.k\oplus \mu'.k$,
with $\mu.(12)=\mu'$, we will adopt the following convention.
Denote by $v_1,\ldots ,\, v_{12}$ the $12$ elements of $\F(V)(3)$.
$$\begin{array}{|c|c|c|c|c|c|}
\hline \textbf{1} & \mu\circ_\textrm{I}\mu \ \leftrightarrow \
(xy)z  & \textbf{5} &  \mu\circ_\textrm{III}\mu \ \leftrightarrow
\ (zx)y &
\textbf{9} &  \mu\circ_\textrm{II}\mu \ \leftrightarrow \ (yz)x \\
 \hline \textbf{2} & \mu'\circ_\textrm{II}\mu \
\leftrightarrow \ x(yz)  & \textbf{6} & \mu'\circ_\textrm{I}\mu \
\leftrightarrow \ z(xy) &
\textbf{10} &  \mu'\circ_\textrm{III}\mu \ \leftrightarrow \ y(zx) \\
 \hline
\textbf{3} &  \mu'\circ_\textrm{II}\mu' \ \leftrightarrow \ x(zy)
& \textbf{7} & \mu'\circ_\textrm{I}\mu' \ \leftrightarrow \ z(yx)
&
\textbf{11} & \mu'\circ_\textrm{III}\mu' \ \leftrightarrow \ y(xz) \\
  \hline
\textbf{4} &  \mu\circ_\textrm{III}\mu' \ \leftrightarrow \ (xz)y
& \textbf{8} & \quad \mu\circ_\textrm{II}\mu' \ \leftrightarrow \
(zy)x &
\textbf{12} &  \mu\circ_\textrm{I}\mu' \ \leftrightarrow \ (yx)z \\
 \hline
\end{array} $$

\begin{figure}
\centering
\includegraphics[scale=0.6]{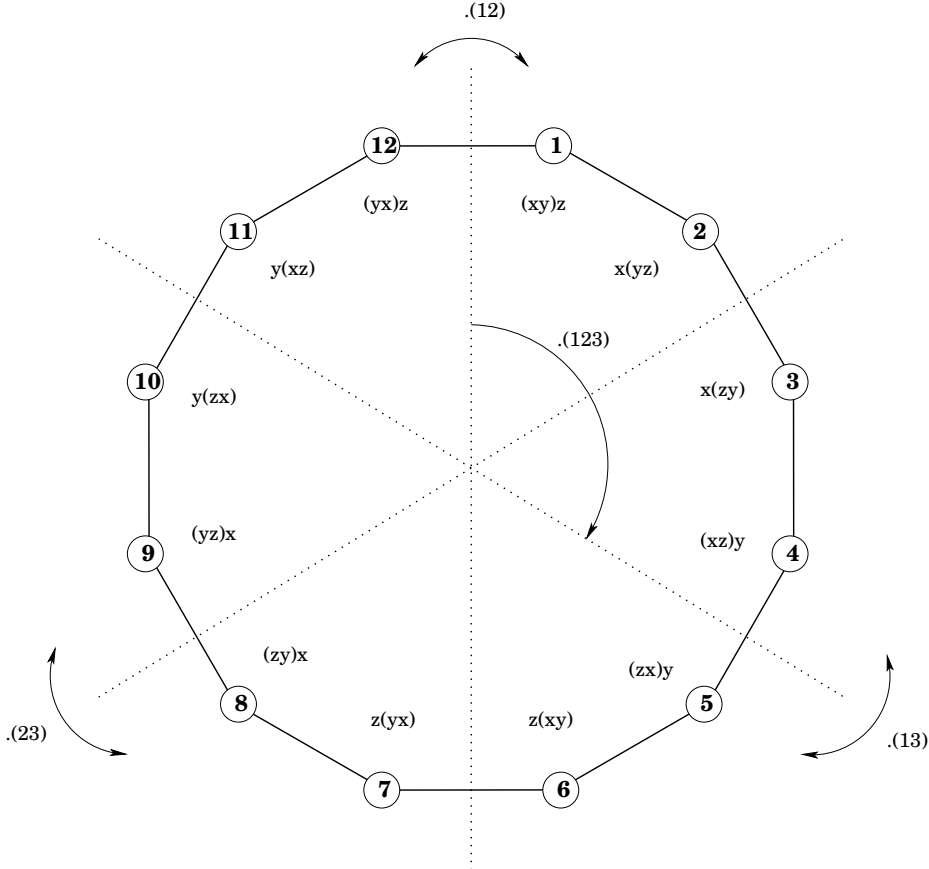}
\caption{The permutoassociahedron} \label{permutoassoc}
\end{figure}

This labelling corresponds to the labelling of the
permutoassociahedron \cite{K}. Figure~\ref{permutoassoc}
represents it with the action of the symmetric group $\Sy_3$.\\

An associative algebra is a vector space with a binary associative
operation, that is $\mu(\mu(a,\, b),\, c)=\mu(a,\, \mu(b,\, c))$.
With these notations, the relations of associativity of the operad
$\A$ become $v_i-v_{i+1}$, for $i=1,\, 3,\, 5,\, 7,\, 9,\, 11$. A
\emph{(right) preLie algebra} is a vector space with a binary
operation such that its associator is right symmetric, that is
$\mu(\mu(a,\, b),\, c)-\mu(a,\, \mu(b,\, c))=\mu(\mu(a,\, c),\,
b)-\mu(a,\, \mu(c,\, b))$. This relation corresponds to
$v_{i}-v_{i+1}+v_{i+2}-v_{i+3}$ for $i=1,\, 5,\, 9$ with our
conventions. The operation of a \emph{$\Pe$-algebra} verifies
$\mu(\mu(a,\, b),\, c)=\mu(a,\, \mu(b,\, c))=\mu(a,\, \mu(c,\,
b))$ which gives here $v_{i}=v_{i+1}=v_{i+2}=v_{i+3}$ for $i=1,\,
5,\, 9$.
Note that $\Pli$ is the Koszul dual of $\Pe$ and vice versa (\emph{cf.} \cite{CL}). \\

We now give an example of computation.

\begin{thm}
\label{Computations}
We have $\Pli\bullet\C=\Zi$,
$\Pli\bullet\A=\Dend$ and $\Pe\circ\Li=\Le$.
\end{thm}

\begin{proo} Denote by $\nu$ the commutative generating
operation of $\C$ and by $w_1,\, w_5, \, w_9$ the related elements
of $\F(\nu.k)(3)$. We write the associativity relation of $\nu$ :
$w_1-w_5=0$ and $w_5-w_9=0$. We have
\begin{eqnarray}
\Psi\left((v_1-v_2+v_3-v_4)\ot(w_1-w_5)\right) &=& \Psi\left(v_1\ot w_1 + v_4\ot w_5\right)\\
\Psi\left((v_1-v_2+v_3-v_4)\ot(w_5-w_9)\right) &=& \Psi\left((v_2-v_3)\ot w_9 - v_4\ot w_5\right)\\
\Psi\left((v_5-v_6+v_7-v_8)\ot(w_1-w_5)\right) &=& \Psi\left((v_7-v_6)\ot w_1 - v_5\ot w_5\right)\\
\Psi\left((v_5-v_6+v_7-v_8)\ot(w_5-w_9)\right) &=& \Psi\left(v_5\ot w_5 + v_8\ot w_9\right)\\
\Psi\left((v_9-v_{10}+v_{11}-v_{12})\ot(w_1-w_5)\right) &=& \Psi\left(-v_{12}\ot w_1 + (v_{10} - v_{11})\ot w_5 \right)\\
\Psi\left((v_9-v_{10}+v_{11}-v_{12})\ot(w_5-w_9)\right) &=&
\Psi\left((v_{11}-v_{10})\ot w_5 - v_9\ot w_9\right)
\end{eqnarray}
The action of $(132)$ sends $(1)$ to $(4)$, $(3)$ to $(6)$ and
$(5)$ to $(2)$. The image of $(1)$ under $(13)$ is $(3)$.
Therefore, we only need to make $(1)$ and $(2)$ explicit. If we
identify $(\mu.k\oplus \mu'.k)\ot \nu.k \ot k.\sgn_{\Sy_2}$ with
$\gamma.k\oplus \gamma'.k$ via the isomorphism of $\Sy_2$-modules
\begin{eqnarray*}
\mu\ot \nu\ot 1 &\mapsto& \gamma \\
\mu'\ot \nu\ot 1 &\mapsto& -\gamma',
\end{eqnarray*}
the morphism $\Psi$ becomes
\begin{eqnarray*}
\Psi((\mu\circ_\textrm{I} \mu) \ot (\nu\circ_\textrm{I}
\nu))=\Psi(v_1\ot w_1)&=& \gamma \circ_\textrm{I} \gamma =z_1
\quad \textrm{and} \\ \Psi((\mu'\circ_\textrm{II} \mu) \ot
(\nu\circ_\textrm{II} \nu))=\Psi(v_2\ot w_1)&=& - \gamma'
\circ_\textrm{I} \gamma=-z_2 . \end{eqnarray*} The image of the
other elements are obtained from these two by the action of
$\Sy_3$. For instance, we have $\Psi(v_3\ot w_1)=-z_3$,
$\Psi(v_4\ot w_1)=z_4$ and $\Psi(v_5\ot w_5)=z_5$.

We get
\begin{eqnarray*}
\Psi\left(v_1\ot w_1 + v_4\ot w_5\right)&=&\gamma
\circ_\textrm{I}\gamma - \gamma \circ_\textrm{III}\gamma' \\
\Psi\left((v_2-v_3)\ot w_9 - v_4\ot w_5\right)&=& -\gamma'
\circ_\textrm{II}\gamma -\gamma' \circ_\textrm{II}\gamma' +\gamma
\circ_\textrm{III}\gamma'.
\end{eqnarray*}
Finally, if we represent the operation $\gamma(x,\,y)$ by $x\star
y$, we have
\begin{eqnarray*}
(x\star y)\star z &=& (x\star z)\star y\\
(x\star z)\star y &=& x\star (z\star y)+x\star (y\star z),
\end{eqnarray*}
where we recognize the axioms of a Zinbiel algebra (\emph{cf.} \cite{L3}).\\

The two other identities are obtained by Koszul duality using
Proposition~\ref{BlackOperad}. From Proposition~\ref{PermWhiteAs}
$\Pe\, \circ\, \A =\D$, we get
 $\Pli \bullet \A =\big( \Pe \circ \A\big)^!
=\big(\D\big)^!=\Dend$. The last equality $\Pe\circ\Li=\Le$ is the
Koszul dual of the first one $\Pli\bullet\C=\Zi$.
\end{proo}

Jean-Louis Loday defined the operad $\Dend$ by two operations such
that their sum is an associative product (see \cite{L3}). In the
same way, he defined the operad $\Zi$ with one product such that
its symmetrized version is a associative (and commutative)
product. This process is often called a \emph{splitting of
associativity}. Proposition~\ref{Computations} shows that we can
interpret the operation $\Pli \, \bullet -$ as a natural way of
splitting the associativity. \\

A commutative algebra is an associative algebra. Therefore, we
have a morphism of operads $\A \to \C$. Since a commutative
algebra is a $\Pe$-algebra and a $\Pe$-algebra is an associative
algebra, the previous morphism factors through $\A \to \Pe \to\C$.
Similarly, a Zinbiel algebra is a dendriform algebra $\Dend \to
\Zi$. We can factor this morphism by a new operad $\Pli \, \bullet
\, \Pe$ using the functor $\Pli\, \bullet -$
$$\xymatrix{\A \ar[r]&
\Pe\ar[r] & \C & \ar@/^/[d]^{\ \Pli\, \bullet -} \\
\Dend \ar[r]& \Pli\, \bullet\,  \Pe \ar[r]& \Zi. &} $$

We describe this new type of algebra.

\begin{thm}
An algebra over the operad $\Pli \bullet \Pe$ is a dendriform
algebra such that the two operations $\prec$ and $\succ$ verify
the two extra relations
\begin{eqnarray*}
&& x\prec (y\prec z) + x \prec (y\succ z) =  x\prec (z\prec y) + x
\prec (z\succ y)\\
&& x\succ(y\prec z)= x\succ (z \succ y).
\end{eqnarray*}
Using the notation $x*y:=x\prec y + x \succ y$, we sum up the $5$
relations of a $\Pli\bullet \Pe$-algebra by
$$\left\{
\begin{array}{l}
 (x\prec y)\prec z = x\prec (y*z)\\
 (x\succ y)\prec z = x\succ (y\prec z)\\
 (x*y)\succ z= x\succ(y\succ z)\\
 x\prec (y*z) = x\prec (z*y)\\
 x\succ(y\prec z)= x\succ (z \succ y).
\end{array}\right.$$
\end{thm}

\draftnote{lui trouver un nom (right) Perm-Dend ...}

A $\Pe$-algebra is an associative algrebra which is symmetric on
the right. A $\Pli \bullet \Pe$-algebra is a dendriform algebra
with right-symmetric relations.

\begin{proo}
Denote by $\omega$ the generating operation of the operad $\Pe$
and by $w_1,\dots , \, w_{12}$ the related elements of
$\F(\omega.k \oplus \omega'.k)(3)$.  The space of relations
$\Psi(R\ot S)$ is generated by the elements
$\Psi\left((v_i-v_{i+1}+v_{i+2}-v_{i+3})\ot (w_j-w_{j+1})
\right)$, for $i\in \{1,\, 5,\, 9 \}$ and $j\in \{1,\, 2,\, 3,\,
5,\, 6,\, 7,\, 9,\, 10,\, 11 \}$. Reduce the computations using
the action of $\Sy_3$ (the symmetries can be seen on the
permutoassociahedron), it remains $5$ relations among which $3$
correspond to the following ones
\begin{eqnarray}
\Psi\left((v_1-v_2+v_3-v_4)\ot(w_1-w_2)\right) &=& \Psi\left(v_1\ot w_1 + (v_2 - v_3)\ot w_2\right)\\
\Psi\left((v_1-v_2+v_3-v_4)\ot(w_5-w_6)\right) &=& \Psi\left(-v_4\ot w_5    -v_1\ot w_6\right)\\
\Psi\left((v_1-v_2+v_3-v_4)\ot(w_7-w_8)\right) &=&
\Psi\left(v_1\ot w_7 + (v_2 - v_3)\ot w_8\right).
\end{eqnarray}

Identify the representation $(\mu.k\oplus \mu'.k)\ot (\omega.k
\oplus \omega'.k) \ot k.\sgn_{\Sy_2}$ with the two copies of
$k[\Sy_2]$ : $\alpha.k[\Sy_2]\oplus \beta.k[\Sy_2]= \alpha.k\oplus
\alpha'.k \oplus \beta.k\oplus \beta'.k$ via the isomorphism of
$\Sy_2$-modules
\begin{eqnarray*}
\mu\ot \omega \ot 1 \mapsto \alpha \quad &\textrm{and}& \quad  \mu'\ot \omega \ot 1 \mapsto -\beta, \\
\mu'\ot \omega' \ot 1 \mapsto -\alpha' \quad &\textrm{and}& \quad
\mu\ot \omega' \ot 1 \mapsto \beta'.
\end{eqnarray*}
The morphism $\Psi$ becomes
\begin{eqnarray*}
\Psi((\mu\circ_\textrm{I} \mu) \ot (\omega\circ_\textrm{I}
\omega))=\Psi(v_1\ot w_1)&=& \alpha \circ_\textrm{I} \alpha,
\\ \Psi((\mu'\circ_\textrm{II} \mu) \ot (\omega' \circ_\textrm{II}
\omega))=\Psi(v_2\ot w_2)&=& -\alpha'
\circ_\textrm{II} \alpha \quad \textrm{and} \\
\Psi((\mu'\circ_\textrm{II} \mu') \ot (\omega' \circ_\textrm{II}
\omega))=\Psi(v_3\ot w_2)&=& - \alpha' \circ_\textrm{II} \beta,
\end{eqnarray*}
for instance. Hence, the relations $(7)$, $(8)$ and $(9)$ are
\begin{eqnarray*}
&& \alpha \circ_\textrm{I} \alpha -\alpha' \circ_\textrm{II}
\alpha -
\alpha' \circ_\textrm{II} \beta \\
&& \alpha \circ_\textrm{III} \beta - \beta' \circ_\textrm{I} \alpha \\
&& \beta' \circ_\textrm{I} \beta' -\beta \circ_\textrm{II} \beta'
- \beta \circ_\textrm{II} \alpha'
\end{eqnarray*}
If we represent the operation $\alpha(x,\, y)$ by $x\prec y$ and
$\beta(x,\, y)$ by $x\succ y$, these $3$ relations become
\begin{eqnarray*}
&& (x\prec y)\prec z = x\prec (y\prec z) + x \prec (y\succ z) \\
&& (z\succ x)\prec y = z\succ (x\prec y)\\
&& z\succ(y\succ x)= (z\succ y)\succ x  + (z\prec y )\succ x,
\end{eqnarray*}
which are the axioms defining dendriform algebras \cite{L3}.

The two other relations are
\begin{eqnarray}
&&\Psi \left((v_1-v_2+v_3-v_4)\ot (w_2-w_3) \right)=
-\alpha'\circ_\textrm{II} \alpha +\alpha'\circ_\textrm{II} \beta'
-\alpha'\circ_\textrm{II} \beta+\alpha'\circ_\textrm{II} \alpha',\\
&&\Psi\left((v_1-v_2+v_3-v_4)\ot (w_6-w_7) \right)=
\beta'\circ_\textrm{I} \alpha-\beta'\circ_\textrm{I} \beta'.
\end{eqnarray}
And they give after identification
\begin{eqnarray*}
&&x\prec (y\prec z) + x \prec (y\succ z) =  x\prec (z\prec y) + x
\prec (z\succ y)\\
&& x\succ(y\prec z)= x\succ (z \succ y).
\end{eqnarray*}
\end{proo}

A $\Pli \bullet \Pe$-algebra is a $\Pe$-algebra with splitting of
the associativity relation.

\begin{pro}
Let $(A,\, \prec,\, \succ)$ be a $\Pli\bullet \Pe$-algebra. With
the operation $* := \prec + \succ$, the vector space $(A,\, *)$
becomes a $\Pe$-algebra.
\end{pro}

\begin{proo}
Consider the sum of the relations.
\end{proo}

Since a $\Pli$-algebra gives a $\Li$-algebra by
anti-symmetrization of the product, we have a morphism of operads
$\Li \xrightarrow{\lambda} \Pli$. Taking the black product of this
morphism with an operad $\Po$, we get a morphism of the form
$\Po=\Li\bullet \Po \xrightarrow{\lambda \bullet \Po} \Pli \bullet
\Po$. A $\Pli\bullet \Po$-algebra has twice more generating
operations then $\Po$ and this morphism corresponds to take the
sum of them. Denote it by $\lambda \bullet \Po=+$. In the previous
cases, we had
$$\xymatrix{\A \ar[r] \ar[d]^{+}&
\Pe\ar[r]  \ar[d]^{+}& \C  \ar[d]^{+}\\
\Dend \ar[r]& \Pli\bullet \Pe \ar[r]& \Zi.} $$

Therefore, the black product with $\Pli$ is a general splitting of
the relations.\\

One interesting property of the black and white products is to
recover classical operads and morphisms between them by means of
products from simpler operads. We have the dual diagram of operads
$$\xymatrix{\A  &
\Pli \ar[l] & \ar[l] \Li  \\
\D \ar[u] & \Pli\circ \Pe \ar[l] \ar[u] & \ar[l] \ar[u] \Le.} $$

The operad $\Pli$ allows to factor the map $\A \leftarrow \Li$.
The notion of $\Pli$-algebra is important and has application in
deformation theory and differential geometry for instance (see
\cite{CL}). The second row $\D \leftarrow \Le$ was introduced by
J.-L. Loday with a view toward applications in algebraic K-theory
(see the introduction of \cite{L3}). The operad $\D$ appears
naturally when one tries to build a bicomplex in algebraic
K-theory with the same form then the one in cyclic homology (the
additive counterpart of algebraic K-theory). Since the operad
$\Pli \circ \Pe$ factors the map $\Le \to \D$, we expect the
operad $\Pli
\circ \Pe$ to appear in these fields in the future.\\

Recall from \cite{CL}, that a basis for $\Pli(n)$ is given by the
set of rooted trees with $n$ vertices labelled by $\{1, \ldots,
n\}$. From Corollary~\ref{Hadamard=tensor}, we have $\Pli \circ
\Pe=\Pli \ot \Pe$. Therefore, a basis for $\Pli \circ \Pe(n)$ is
provided by the set of rooted trees with $n$ vertices labelled by
$\{1, \ldots, n\}$ with one vertex emphasized. We leave to the
reader to describe the composition map of this operad. (Use the
composition of $\Pli$ based on rooted trees given in \cite{CL}
with the fact that only the insertion of a tree in an emphasized
vertex keeps a vertex emphasized).

\subsection{A counterexample}\label{Counterexample}

In this section, we show that the category of Koszul operads is
not stable by white and black products. We exhibit a pair of
Koszul operads whose black product is not Koszul. \\

Consider the \emph{nilpotent} operad $\mathcal{N}$ defined by a
generating skew-symmetric binary operation such that every
composition of it vanishes.

\begin{lem}
The operad $\Pli\bullet \mathcal{N}$ is equal to the quadratic
operad generated by a binary operation $\diamond$ with the
following relations : $(x\diamond y)\diamond z =0$ and $ x
\diamond (y \diamond z)=x \diamond (z \diamond y)$, for every $x,\
y,\, z$.
\end{lem}

\begin{proo}
We use the same notations $v_i$ for the space $R$ of relations of
the operad $\Pli$. The space $S$ of relations of the nilpotent
operad is generated by $w_1$, $w_5$ and $w_9$. By symmetry of the
relations, we only have to compute the three terms
\begin{eqnarray}
&&\Psi \left((v_1-v_2+v_3-v_4)\ot w_1 \right)=\Psi(v_1\ot w_1)=
\diamond \circ_\textrm{I} \diamond,\\
&&\Psi\left((v_1-v_2+v_3-v_4)\ot w_5 \right)=
\Psi((-v_2+v_3)\ot w_5)= -\diamond'\circ_\textrm{II} \diamond + \diamond'\circ_\textrm{II} \diamond', \\
&&\Psi\left((v_1-v_2+v_3-v_4)\ot w_9 \right)= \Psi(-v_4\ot w_9)=
\diamond\circ_\textrm{III} \diamond.
\end{eqnarray}
They correspond to
$$(x\diamond y)\diamond z =0, \quad x
\diamond (y \diamond z)=x \diamond (z \diamond y) \quad
\textrm{and} \quad (x\diamond z)\diamond y =0.$$
\end{proo}

\begin{thm}
The operad $\Pli\bullet \mathcal{N}$ is not Koszul.
\end{thm}

\begin{proo}
Because of its relations, the operad $\Pli\bullet \mathcal{N}$ has
no operations in arity $n$ for $n$ greater than $4$, that is
$\left( \Pli\bullet \mathcal{N}\right)(n)=0$ for $n\ge 4$. Recall
that the Poincar\'e series of an operad $\Po$ is defined by
$\displaystyle f_\Po(x):=\sum_{n\ge 1} \frac{\textrm{dim}
(\Po(n))}{n!} x^n$ (see \cite{GK} Section $3$ or \cite{L3}
Appendix B.5.c.). When an operad $\Po$ is Koszul, its Poincar\'e
series and the Poincar\'e series of its dual verify the equation
$f_{\Po^!}(-f_{\Po}(-x))=x$ (\cite{GK} Formula $(3.3.2)$). The
Poincar\'e series $-f_{\Pli\bullet \mathcal{N}}(-x)$ is $x - x^2
+\frac{1}{2}x^3$. Its inverse for the composition is
$$x + x^2 + \frac{3}{2}x^3  + \frac{5}{2}x^4+ \frac{17}{4}x^5+
7x^6+ \frac{21}{2}x^7+ \frac{99}{8}x^8+ \frac{55}{16}x^9
\underline{-\frac{715}{16}x^{10}} +\cdots.$$ Since the
$10^{\textrm{th}}$ coefficient is negative, this series does not
correspond to the Poincar\'e series of an operad. Therefore the
operad $\Pli\bullet \mathcal{N}$ is not Koszul.
\end{proo}

The operad $\Pli$ is Koszul (see \cite{CL} for a proof in
characteristic $0$ and \cite{CV} for a more general one). Any
nilpotent operad is Koszul (the Koszul dual is a free operad,
which is Koszul). So the operad $\Pli\bullet \mathcal{N}$ is the
black product of two Koszul operads which is not a Koszul operad.
This result comes from the fact the morphism $\Psi$ (and the
morphism $\Phi$) is not an isomorphism in general. The morphism
$\Psi$ is a projection and kills part of the relations. Therefore,
the coherence between the relations, expressed by the Koszul
property, does not hold anymore.

\subsection{Adjunction}\label{SymMonCat}In this section, we
generalize the main result of \cite{M2}  about the adjunction
between the black and the white products to $k$-ary quadratic
operads.\\

Let $k$ be an integer greater than $2$. Consider the category of
\emph{$k$-ary quadratic operads}, that is quadratic operads
generated by a finite dimensional $\Sy$-module concentrated in
arity $k$. A morphism between two $k$-ary quadratic operads
$\F(V)/(R)$ and $\F(W)/(S)$ is a morphism induced by a map of
$\Sy_k$-modules $V(k) \to W(k)$. Denote this category by
\emph{$k$.q-Op}.

One can generalize the basis and the non-degenerate bilinear form
of \ref{BinaryQuadraticOperad} for the binary case to the $k$-ary
case. Then Lemma~\ref{PhiPsiOrthogonal} and
Proposition~\ref{BlackOperad} also hold in $k$.q-Op, which defines
black products in this category. Recall from V. Gnedbaye
\cite{Gnedbaye} the notion of \emph{$k$-Lie algebra}, that is a
module endowed with a $k$-ary antisymmetric bracket satisfying a
generalized Jacobi relation. We denote the associated operad by
$\Li^{<k>}$. Gnedbaye proved that $\Li^{<k>}$ is the Koszul dual
operad of $\C^{<k>}$ (denoted $stAs^{<k>}$ in \cite{Gnedbaye}),
where a $\C^{<k>}$-algebra is module equipped with a $k$-ary
commutative and totally associative operation.

\begin{pro}
The black and white products endow the category of $k$-ary
quadratic operads with a structure of symmetric monoidal category,
where the operad $\Li^{<k>}$ is the unit object for $\bullet$ and
the operad $\C^{<k>}$ is the unit object for $\circ$.
\end{pro}

\begin{proo}
The same arguments as in \ref{Hadamard=White} show that for a
$k$-ary quadratic operad $\Po$, we have $\C^{<k>} \circ
\Po=\C^{<k>} \otimes \Po=\Po$. If $n \in (k-1).\NN +1$,
$\C^{<k>}(n)=\KK$, otherwise $\C^{<k>}(n)=\Po(n)=0$. The rest of
the proof is straightforward.
\end{proo}

\begin{thm}\label{Adjunction}
There is a natural isomorphism $\Hom_{k.\textrm{q.Op}}(\Po\bullet
\Qo,\, \mathcal{R})\cong \Hom_{k.\textrm{q.Op}}(\Po,\, \Qo^!\circ
\mathcal{R})$. Hence, the tensor category of $k$-ary quadratic
operads with the black product $\bullet$ is endowed with an
internal $\Hom$ object denoted $\h(\Qo,\, \mathcal{R}):=\Qo^!
\circ \mathcal{R}$. Dually,
 $\mathop{cohom}(\Po,\, \Qo):=\Po
\bullet \Qo^!$ defines an internal $\mathrm{coHom}$ object in
$(k.\textrm{q.Op}, \circ, \Co_k)$.
\end{thm}

\begin{proo}
Let $\Po=\F(V)/(R)$, $\Qo=\F(W)/(S)$ and $\mathcal{R}=\F(X)/(T)$
be three $k$-ary quadratic operads. There is a one-to-one
bijection between maps $f \, :\, V \ot W \ot \sgn \to X $ and maps
$g\, :\, V \to W^\vee \ot X$. It remains to show that $\F(f)\big(
\Psi(R\ot S)\big) \subset T$ is equivalent to $\F(g)(R)\subset
\Phi^{-1}(S^\perp\ot \F(X) + \F(W^\vee)\ot T)$. By
Lemma~\ref{PhiPsiOrthogonal}, we have
\begin{eqnarray*}
<\F(g)(R), \Phi^{-1}(S^\perp\ot \F(X) + \F(W^\vee)\ot
T)^\perp>_{W^\vee\ot X} &=& <\F(g)(R), \Psi(S\ot
T^\perp)>_{W^\vee\ot X}\\
&=& <\F(f)\left( \Psi(R\ot S)\right), T^\perp>_{X},
\end{eqnarray*}
which concludes the proof.
\end{proo}

For another point of view on this type of adjunction and
$\textrm{coHom}$ objects in another operadic setting, we refer the
reader to D. Borisov and Yu.I. Manin \cite{BorisovManin}.

\begin{cor}
Let $\Po$ be a $k$-ary quadratic operad. The operad
$\en(\Po):=\Po\bullet \Po^!$ is a comonoid in $(k.\textrm{q.Op},
\circ, \Co_k)$.
\end{cor}

\begin{proo}
The proof comes from general methods of $\textrm{coHom}$ objects.
\end{proo}

Composing $\Delta$ with  $\bar{\Phi}\, :\, \en(\Po) \circ \en(\Po)
\to \en(\Po) \ot \en(\Po)$, we get that $\Po^! \bullet \Po$ is a
comonoid for the tensor product, that is a Hopf operad (see
\ref{Bimonoids}).

\begin{thm}
For every $k$-ary quadratic operad $\Po$, the operad
$\en(\Po)=\Po^!\bullet \Po$ is a Hopf operad.
\end{thm}

The first example is $\C=\C \bullet \Li$. Other examples are
$\Pli \bullet \Pe$, $\Zi \bullet \Le$.\\

In \cite{M}, Yu. I. Manin proved the equivalent theorem for
quadratic algebras. This allowed him to realize quantum groups as
black products of an algebra with its Koszul dual algebra. In this
spirit, the previous theorem gives a method to get new ``quantum
groups'', that is Hopf operads. \\

The tensor product of a Lie algebra with a commutative algebra is
again a Lie algebra (Courant algebras for instance). This result
can be widely generalized. Let $\Po$ be a $k$-ary quadratic
operad. For any $\Po^!$-algebra $A$ and any $\Po$-algebra $B$,
their tensor product $A\ot B$ is a $\Li^{<k>}$ algebra (see
\cite{L3} Appendix B.5.a. for a proof in the binary case and see
\cite{GW} Theorem $2.3$ for a proof in the ternary case). In the
language of operads, it means that there exists a morphism of
operads $\Li^{<k>} \xrightarrow{l} \Po^!\ot \Po$. In the
particular case of $\Po=\mathcal{L}eib$ and
$\Po^!=\mathcal{Z}inb$, J.-L. Loday and I. Dokas refined this
result and proved in \cite{DokasLoday} that the previous map
factors through $\Pli$. We now give a conceptual proof of the
existence of the map from $\Li^{<k>}$ to $\Po^!\ot \Po$ and show
that is always comes from a composite with the white product.

\begin{pro}\label{universalLiek}
For every $k$-ary quadratic operad $\Po$, there is a canonical
morphism of operads $\Li^{<k>} \xrightarrow{i} \Po^! \circ \Po$,
defined by the commutative diagram
$$\xymatrix{\Li^{<k>} \ar[rr]^{l} \ar[dr]_{i} & &\Po^!\ot \Po \\
&\Po^! \circ \Po \ar[ur]_{\bar{\Phi}}.  &} $$
\end{pro}

\begin{proo}
Apply Theorem~\ref{Adjunction} to the triple of operads
$\Li^{<k>}$, $\Po$ and $\Po$. We get a natural isomorphism
$\Hom_{k.\textrm{q.Op}}(\Li^{<k>} \bullet \Po,\, \Po)\cong
\Hom_{k.\textrm{q.Op}}(\Li^{<k>},\, \Po^!\circ \Po)$. Since
$\Li^{<k>}$ is the unit object for $\bullet$, we have
$\Hom_{k.\textrm{q.Op}}(\Po,\, \Po)\cong
\Hom_{k.\textrm{q.Op}}(\Li^{<k>},\, \Po^!\circ \Po)$. Define
$\Li^{<k>} \xrightarrow{i} \Po^! \circ \Po$ to be the image of the
identity of $\Po$ under this isomorphism.
\end{proo}

\subsection{Cohomology operations}\label{Cohomology Operations}

In this section, we recall the definition of the intrinsic Lie
bracket on the chain complex defining the cohomology
theories for algebras over a Koszul operad. We use the previous section to define
 another Lie bracket on the same
space. Because of the symmetries, this operation vanishes on cohomology. \\

 Let $(\Po, \mu^\Po)$ and $(\Qo,\,
\mu^\Qo)$ be two augmented dg-operads and let $\rho \ : \ \Po \to
\Qo$ be a morphism of augmented dg-operads. This morphism makes
$\Qo$ a module over $\Po$. Denote by $\mu_{(1,1)}^\Po$ the partial
composition of $\Po$, that is the composition of two non-trivial
operations of $\Po$.

\begin{dei}[Derivation]
A homogenous morphism $\partial\, :\, \Po \to \Qo$ is a
\emph{homogenous derivation of $\rho$}  if
$$\partial\circ \mu^\Po_{(1,\,
1)}=\mu^\Qo_{(1,\, 1)}\circ (\partial\otimes \rho) +
\mu^\Qo_{(1,\, 1)}\circ (\rho\otimes \partial).$$
 This formula, applied to elements $p_1 \ot
p_2$ of $\Po\ot \Po$, where $p_1$ and $p_2$ are homogenous
elements of $\Po$, gives
\begin{eqnarray*}
\partial\circ\mu^\Po(p_1\ot p_2) &=&
\mu^\Qo\big(\partial(p_1)\ot\rho(p_2)\big) +
(-1)^{|\partial||p_1|}\mu^\Qo\big(\rho(p_1)\ot
\partial(p_2)\big).
\end{eqnarray*}
A \emph{derivation} is a sum of homogenous derivations. The set of
homogenous derivations with respect to $\rho$ of degree $n$ is
denoted $\Der^n_{\rho}(\Po,\, \Qo)$ and the set of derivations is
denoted $\Der_{\rho}^\bullet(\Po,\, \Qo)$ or simply $\Der(\Po,\,
\Qo)$ when the morphism $\rho$ is obvious.
\end{dei}

We recall the definition of the cohomology of $\Po$-algebras, when
$\Po$ is a Koszul operad. Let $A$ be a $\Po$-algebra, that is
there is a morphism of operads $\Po \xrightarrow{\phi} \End(A)$.
Denote by
$\Omega(\Po^{\ac})=(\F(s^{-1}\oPo^{\ac}),\bar{\partial)}$ the
cobar construction of $\Po^{\ac}$, where the differential
$\bar{\partial}$ is the unique derivation which extends the
partial coproduct of the Koszul dual cooperad $\Po^{\ac}$. Since
$\Po$ is a Koszul operad, $\Omega(\Po^{\ac})$ is a quasi-free
resolution of $\Po$.

$$\xymatrix{\Omega(\Po^{\ac}) \ar[r]^{\sim} \ar[dr]_{\rho}&
\Po \ar[d]^{\phi} \\ & \End(A).}$$

\begin{lem}
Let $(\R,\, \bar{\partial})  \xrightarrow{\varepsilon} \Po$ be a
resolution of $\Po$ and let $f$ be an homogenous derivation of
degree $n$ in $\Der^n_{\rho} \left( \R ,\, \End(A)\right)$. One
has $f\circ \bar{\partial} \in \Der^{n-1}_{\rho} \left( \R,\,
\End(A) \right)$.
\end{lem}

\begin{proo}
The degree of $f \circ \bar{\partial}$ is $n-1$. It remains to
show that $f \circ \bar{\partial}$ is a derivation. Since $(\R,
\bar{\partial})$ is a dg-operad, we have
\begin{eqnarray*}
f\circ \bar{\partial}\circ  \mu^\R  &=&   f \circ \mu^\R\circ (
\bar{\partial}\ot \Id  +\Id \ot \bar{\partial})\\
&=&  \mu^{\End(A)} \circ (f\ot \rho + \rho \ot f) \circ
 ( \bar{\partial}\ot \Id +\Id \ot \bar{\partial})\\&=&
\mu^{\End(A)} \circ \big((f\circ \bar{\partial}) \ot \rho + \rho
\ot (f \circ \bar{\partial}))\big).
\end{eqnarray*}
Since $\R$ is concentrated in non-negative degree and $A$ is
concentrated in degree $0$, the composite $\rho\circ
\bar{\partial} =\phi\circ \varepsilon \circ \bar{\partial}$ is
null.
\end{proo}

The deformation theory of the map $\Po \xrightarrow{\phi} \End(A)$
is studied via the following cochain complex defined by M. Markl
in \cite{Markl96a}. The cohomology of a $\Po$-algebra $A$ is
defined on the space of derivations of $\rho$ (see also
\cite{Quillen}, \cite{Markl96b} and \cite{KontsevichSoibelman}).

\begin{dei}
The \emph{cohomology of a $\Po$-algebra $A$} is defined by the
(deformation) chain complex
$$C^\bullet_\Po(A):=\big( \Der^\bullet_\rho\left(\Omega(\Po^{\ac}) ,\,
\End(A)\right), \partial),$$ where the differential $\partial$ is
the pullback by $\bar{\partial}$, that is $\partial(f):= f \circ
\bar{\partial}$.
\end{dei}

Since $\Omega(\Po^{\ac})$ is a free operad, we have
\begin{eqnarray*} \Der^\bullet_\rho\left(\Omega(\Po^{\ac}),\,
\End(A)\right) \cong \Hom_\Sy^{\bullet}(\F(s^{-1}\oPo^{\ac}),
\End(A))\cong \Hom_\Sy^{\bullet-1}(\oPo^{\ac}, \End(A))\cong
\Hom^{\bullet-1}_\KK(\oPo^{\ac}(A), A),
\end{eqnarray*}
where $\Hom_\Sy(M,N)$ denotes the set of $\Sy$-equivariant maps between the $\Sy$-modules $M$ and $N$.\\

As in the paper of M. Kontsevich and Y. Soibelman
\cite{KontsevichSoibelman}, we can consider the augmented chain
complex $\Hom^{\bullet}_\Sy(\Po^{\ac}, \End(A))\cong
\Hom^{\bullet}_\KK(\Po^{\ac}(A), A)$. Up to a shift of degree, the
last space corresponds to the Hochschild (co)chain complex
 for associative algebras, Harrison cohomology of
commutative algebras and Chevalley-Eilenberg for Lie algebras.
Notice that in the literature, this cohomology is called the
cohomology of $A$ with coefficient in $A$. Since this chain
complex is defined to control the deformation of the morphism
$\Phi$, that is the structure of $\Po$-algebra on $A$, we call it
the cohomology $\Po$ with coefficient in $A$ or simply the
cohomology of $A$, once the
operad is chosen. \\

In these three cases, the chain complex is a dg-Lie algebra whose
bracket is often called the \emph{intrinsic bracket} (see J.
Stasheff \cite{Stasheff93}). The space $\Hom_\KK(\Po^{\ac},
\End(A))$  of morphisms from a dg-cooperad to a dg-operad is an
$\Sy$-module with the action by conjugation, that is
$(f.\sigma)(p):=\big(f(p.\sigma^{-1})\big).\sigma$. Moreover, it
is a dg-operad, called the \emph{convolution operad} in
\cite{BergerMoerdijk03} Section 1. On the direct sum of the
$\Sy_n$-modules of an operad, one can define a preLie product
$\star$ whose anti-symmetrization gives a Lie bracket. When the
operad is the convolution operad $\Hom^\bullet_\KK(\Po^{\ac},
\End(A))$, the preLie product is a degree $0$ operation given by
$$f\star g \ := \ \Po^{\ac} \xrightarrow{\Delta'} \Po^{\ac} \ot \Po^{\ac}
\xrightarrow{f\ot g} \End(A) \ot \End(A) \xrightarrow{\mu_A}
\End(A)$$ where $\Delta'$ is the partial coproduct of the cooperad
$\Po^{\ac}$. The intrinsic Lie bracket is defined by $[f,g]:=f
\star g - (-1)^{|f||g|}g \star f$. The space of $\Sy$-equivariant
morphisms $\Hom_\Sy(\Po^{\ac}, \End(A))$ is equal to the space of
invariants $\Hom_\KK(\Po^{\ac}, \End(A))^\Sy$ with respect to the
action by conjugation. It is a subspace of the convolution operad
$\Hom_\KK(\Po^{\ac}, \End(A))$ stable under the preLie product
$\star$. (See for instance \cite{VdL} for a proof of this in the
coinvariant context. Since we work over a field $\KK$ of
characteristic $0$, the isomorphism between invariants and
coinvariants allows us to conclude.) The induced Lie bracket on
$C_\Po^\bullet(A)=\Hom_\Sy(\Po^{\ac}, \End(A))$ defines
an intrinsic Lie bracket on cohomology. (We refer the reader to \cite{MerkulovVallette07}
a complete study of the deformation complex). \\

When $\Po=\A$, it is exactly the structure defined by M.
Gerstenhaber in \cite{Gerstenhaber63} and when $\Po=\Li$ it is the
Lie bracket of Nijenhuis and Richardson, which controls the formal
deformations of $\Po$-algebra structure (see D. Balavoine
\cite{Balavoine} Section $4$).\\

Let $A$ be a $\Po$-algebra and $C$ be a $\Po^{\ac}$-coalgebra, we
have by Proposition~\ref{universalLiek} that $\Hom_\KK(C,A)$ is
naturally endowed with a structure of $\Po^!\circ \Po$-algebra and
$\Li^{<k>}$-algebra (see also Section~\ref{Exemples Hom(C,A)}).
Applied to $C=\Po^{\ac}(A)$, this result gives that the chain
complex $C^\bullet_\Po(A)$ is a $\Po\circ \Po^!$-algebra and a
$\Li^{<k>}$-algebra. In the binary case, it means that
$\Hom^\bullet_\KK(\Po^{\ac}(A), A)$ is equipped with another Lie
bracket $\{\, ,\}$ of degree $-1$. Let $\alpha$ be a morphism of
degree $-1$ defined as follows
$$\alpha \ : \ \Po^{\ac} \epi \Po^{\ac}_{(1)}=\Po^{\ac}(2)=s \Po(2)
\to \Po(2) \mono \Po \xrightarrow{\Phi} \End(A).$$

It is a a \emph{twisting cochain}, that is $\alpha$ is solution to
the Maurer-Cartan equation $\alpha \star \alpha=0$, when $\Po$ and
$A$ are concentrated in degree $0$ (see Section $2.3$ of
Getzler-Jones \cite{GJ}). The Lie bracket $\{f , g \}$ is equal to
$$\Po^{\ac} \xrightarrow{\Delta} \Po^{\ac} \circ \Po^{\ac} \epi \Po^{\ac}(2)\ot_{\Sy_2}{\Po^{\ac}}^{\ot 2}
\xrightarrow{\alpha \ot (f\ot g + (-1)^{|f||g|}g\ot f)} \End(A)
\ot \End(A) \xrightarrow{\mu_A} \End(A).$$

Note that in the binary case, the latter Lie bracket $\{ \, ,\}$
is not equal to the intrinsic Lie bracket $[\, ,]$. For instance,
there is a shift of degree between the two.

\begin{lem}\label{Formula}
For every $f$ and $g$ in $C^\bullet_\Po(A)$, we have
\begin{eqnarray*}
\partial(f)=[f, \alpha] \quad \textrm{and} \quad
 \{ f,g\}=\partial f \star g + (-1)^{|f|}f \star \partial g -
\partial (f\star g).
\end{eqnarray*}
\end{lem}

\begin{proo}
The proof is straightforward and left to the reader.
\end{proo}

Equipped with the intrinsic Lie bracket $[\, , ]$,
$C^\bullet_\Po(A)$ becomes a dg-Lie algebra. The second formula
shows that the Lie bracket $\{\, , \}$ vanishes on cohomology.
This result and formula can be explained as follows. The preLie
product comes from the partial composition of the operad. The
general composition product of an operad defines symmetric braces.
Since the partial composition of the operad generates the global
one, the preLie product generates the symmetric braces. (See also
J.-M. Oudom and D. Guin \cite{OudomGuin04} and \cite{LadaMarkl05}
for a proof of this result). Therefore, we cannot expect to have
other products than the intrinsic Lie bracket in general. In
particular examples, it would be interesting to see if the
structure of $\Po \circ \Po^!$-algebra induces a non-trivial
structure on cohomology. We will see in \ref{AdjonctionReguliere}
how to refine this study when the operad is not symmetric
(regular).






\section{Black and white square-products for regular
operads}\label{BW Square}

K. Ebrahimi-Fard and L. Guo in \cite{EFG} and J.-L. Loday in
\cite{L2} defined and used an analog of Manin's black product for
regular operads that they called the \emph{black square product}.
In this section, we give the conceptual definitions of Manin's
black and white square products for regular operads. They are not
equal to the black and white ``circle''-products in the category
of operads. Actually, they come from the black and white products
in the category of non-symmetric operads.

\subsection{Definitions of non-symmetric and regular operads}

Recall that a \emph{non-symmetric operad} is an operad without the
actions of the symmetric groups. From a non-symmetric operad $\{
\Po'_n\}_{n\in \NN^*}$, we can associate an $\Sy$-module by the
collection of the free $\Sy_n$-modules on $\Po(n):=\Po'_n \ot_k
k[\Sy_n]$. The composition product for the operad $\Po$ is defined
from the non-symmetric one. Such an operad is called a
\emph{regular operad}. Denote $\Sigma$ this functor from
non-symmetric operads to operads. Therefore, the category of
regular operads is the image of $\Sigma$ and is equivalent to the
category of non-symmetric operads. Denote by $U$ the inverse
functor :
$$\xymatrix{ \textrm{Non-symmetric Operads} \ \ar@_{->}@<-1ex>[r]_(0.57){\Sigma} &
\ar@_{->}@<-1ex>[l]_(0.42){U} \ \textrm{Regular Operads}. }
$$

Let $\Po=\F(V)/(R)$ be a binary quadratic regular operad. In that
case, we have that $V$ and $R$ are regular modules, that is $V=V'
\ot_k k[\Sy_2]$ and $R=R' \ot_k k[\Sy_3]$. The non-symmetric
operad $\Po'=U(\Po)$ is once again binary and quadratic. It is
given by $\Po'=F(V')/(R')$.

\subsection{Definitions of black and white square-products}

A non-symmetric operad is a monoid in the category of non-negative
graded modules with a non-symmetric version of $\circ$ (see
Appendix~\ref{Alg-Op-Properad}). Under the Hadamard product, this
category forms a $2$-monoidal category. Hence, we can apply
arguments of section~\ref{ManinProducts} and consider the morphism
$\Phi$ and the induced white product for non-symmetric operads.
 From two binary quadratic regular operads $\Po=\F(V)/(R)$ and
$\Qo=\F(W)/(S)$, we study the associated white product
$$U(\Po) \circ U(\Qo) := \F(V'\ot W')/\big(
\Phi^{-1}(R'\ot \F(W') + \F(V')\ot S')\big).$$ The idea is now to
come back to the category of regular operads using the functor
$\Sigma$.

\begin{dei}[White square-product]
The \emph{white square-product} of two binary quadratic regular
operads $\Po$ and $\Qo$ is defined by the following formula
$$\Po \Sq \Qo := \Sigma(U(\Po)\circ U(\Qo)).$$
\end{dei}

 More explicitly, the white square-product of $\Po$ and $\Qo$ is
equal to $\Po\Sq \Qo=\F(V'\ot W' \ot k[\Sy_2])/\big( (
\Phi^{-1}(R'\ot \F(W') + \F(V')\ot S'))\ot k[\Sy_3]   \big)$.

Note that the definition given above does not correspond to
Definition $3.1$ of K. Ebrahimi-Fard and L. Guo in \cite{EFG} (See
Remark below.)

\begin{pro}
\label{TensorVSSquare} Let $A$ be a $\Po$-algebra and $B$ a
$\Qo$-algebra, their tensor product $A\ot B$ is an algebra over
the white square-product $\Po \Sq \Qo$.
\end{pro}

\begin{proo}
The proof is the same than Proposition~\ref{tensor-white}.
\end{proo}

Let $V$ be an $\Sy_2$-module. The part $\F(V)(3)$ with $3$ inputs
of the free operad on $V$ is isomorphic to
$$\F(V)(3)=\big(V\ot_{\Sy_2}(V\ot k \oplus k \ot V)\big) \ot_{\Sy_2} k[\Sy_3],$$
where the summand $V\ot(V \ot k)$ corresponds to the compositions
on the left $\vcenter{\xymatrix@M=0pt@R=5pt@C=5pt{ \ar@{-}[dr] &
&\ar@{-}[dl] & &  \\
& \ar@{-}[dr] & &\ar@{-}[dl]  & \\
& &\ar@{-}[d] & & \\
& & \\ & & }} $ and the summand $V\ot(k \ot V)$ corresponds to the
compositions on the right $\vcenter{\xymatrix@M=0pt@R=5pt@C=5pt{
 & &\ar@{-}[dr] & & \ar@{-}[dl]  \\
& \ar@{-}[dr] & &\ar@{-}[dl]  & \\
& &\ar@{-}[d] & & \\
& & \\ & & }} $. When $V$ is a sum of regular representations
$V=V'\ot k[\Sy_2]$, we have $\F(V)(3)=\big( V'\ot (V'\ot k) \oplus
V' \ot(k \ot V') \big)\ot k[\Sy_3]$. Therefore, $\F(V)(3)$ can be
identify with $2$ copies of $V'\ot V'$ represented by the
following types of tree $\vcenter{\xymatrix@M=0pt@R=5pt@C=5pt{
\ar@{-}[dr] &
&\ar@{-}[dl] & &  \\
& \ar@{-}[dr] & &\ar@{-}[dl]  & \\
& &\ar@{-}[d] & & \\
& & \\ & & }}$ and $\vcenter{\xymatrix@M=0pt@R=5pt@C=5pt{
 & &\ar@{-}[dr] & & \ar@{-}[dl]  \\
& \ar@{-}[dr] & &\ar@{-}[dl]  & \\
& &\ar@{-}[d] & & \\
& & \\ & & }}$. These two copies correspond to the part of arity
$3$ of the free non-symmetric operad on $V'$. We denote the first
composition based on the pattern $\TreeL$ by
$\mu \circ_1 \nu$ and the second one based on $\TreeR$ by $\mu \circ_2 \nu$, where $\mu$ is below $\nu$.\\

In the appendix B of \cite{L3}, J.-L. Loday described the
non-degenerate bilinear form $<\,,\,>$ for regular operads. It
comes from the following one for non-symmetric operads.
\begin{eqnarray*}
\F(V')(3)\ot \F(V'^*)(3) \xrightarrow{<\, , \, >}  k \\
<\mu \circ_1 \nu,\, \zeta \circ_1 \xi> := +\zeta(\mu).\xi(\nu) \\
<\mu \circ_2 \nu,\, \zeta \circ_2 \xi> := -\zeta(\mu).\xi(\nu),
\end{eqnarray*}
the other products being null.\\

We define the black product of binary non-symmetric operad like in
\ref{BlackproductOperadsSection} (the non-degenerate bilinear form
is given below). Applying the same ideas, we have the analog of
Lemma~\ref{PhiPsiOrthogonal} and Proposition~\ref{BlackOperad}.

\begin{lem}
\label{Form-Blacksquare} Let $\Po=\F(V'\ot k[\Sy_2])/(R'\ot
k[\Sy_3])$ and $\Qo=\F(W'\ot k[\Sy_2])/(S'\ot k[\Sy_3])$ be two
regular operads such that the $V'$ and $W'$ are finite
dimensional. The orthogonal of $\Psi(R'\ot S')$ for $<\, ,\, >$ is
$\left( \Phi_*^{-1}(R'^\perp \ot \F(W'^*) + \F(V'^*)\ot
S'^\perp)\right)$.
\end{lem}

\begin{dei}[Black square product]
Let $\Po'=\F(V')/(R')$ and $\Qo'=\F(W')/(S')$ be two binary
quadratic non-symmetric operads with finite dimensional generating
spaces. Define their black product by the formula
\begin{eqnarray*}
\Po'\bullet \Qo' &=&  \F(V'\ot W')/\big( \Psi(R'\ot S' ) \big).
\end{eqnarray*}
The \emph{black square product} of two binary quadratic regular
operads is defined by
$$\Po \BlackSq \Qo := \Sigma (U(\Po) \bullet U(\Qo)). $$
\end{dei}

\begin{pro}
\label{BlackSquareOperad} For binary quadratic regular operads
generated by finite dimensional modules, this definition of black
product verifies $\big(\Po \BlackSq \Qo\big)^! = \Po^! \Sq \Qo^!
$.
\end{pro}

Finally, we can use the particular form of the bilinear product
$<\, ,\, >$ to make explicit the morphism $\Psi$ and show that the
black square-product defined here corresponds to the one of
\cite{L2} and \cite{EFG}.

\begin{pro}
\label{BlackSquare} Under the same hypotheses, let $r\ot s$ be an
elementary tensor of $R'\ot S'$. Denote $r=r_1+r_2$, where $r_1$
is the part of $r$ corresponding to the compositions of the form
$\TreeL$ and $r_2$ is the part of $r$ corresponding to the
compositions of the form $\TreeR$. In the same way, write
$s=s_1+s_2$. The image of $r\ot s$ under $\Psi$ is $\Psi(r\ot
s)=\Phi^{-1}(r_1\ot s_1)-\Phi^{-1}(r_2\ot s_2)$.
\end{pro}

\begin{proo}
Note that $r_1\ot s_1$ and $r_2\ot s_2$ belong to $\textrm{Im}\,
\Phi$. For $X \in \F(V'^* \ot W'^*)(3)$, denote the image of $X$
under $\Phi_*$ by $\Phi_*(X)=\sum \Phi_{V'^*}(X)\ot
\Phi_{W'^*}(X)$. More precisely, we decompose the image of $X$
under $\Phi_*$ with the two types of compositions
$\Phi_*(X)=\Phi_*(X_1+X_2)=\sum_1 \Phi_{V'^*}(X_1)\ot
\Phi_{W'^*}(X_1) + \sum_2 \Phi_{V'^*}(X_2)\ot \Phi_{W'^*}(X_2)$.
We have
\begin{eqnarray*}
& & <\Phi^{-1}(r_1\ot s_1)-\Phi^{-1}(r_2\ot s_2),\,
\F({\Xi})(X)>_{V'\ot W'} \\ &=& <\Phi^{-1}(r_1\ot s_1),\,
\F({\Xi})(X)_1>_{V'\ot W'}-<\Phi^{-1}(r_2\ot
s_2),\, \F({\Xi})(X)_2>_{V'\ot W'} \\
&=& \sum_1 <r_1,\, \Phi_{V'^*}(X_1)>_{V'}.<s_1,\,
\Phi_{W'^*}(X_1)>_{W'} + \sum_2 <r_2,\,
\Phi_{V'^*}(X_2)>_{V'}.<s_2,\,
\Phi_{W'^*}(X_2)>_{W'}=\\
&=& \sum_1 <r,\, \Phi_{V'^*}(X_1)>_{V'}.<s,\,
\Phi_{W'^*}(X_1)>_{W'} + \sum_2 <r,\, \Phi_{V'^*}(X_2)>_{V'}.<s,\,
\Phi_{W'^*}(X_2)>_{W'}=\\
&=& \sum <r,\, \Phi_{V'^*}(X)>_{V'}.<s,\, \Phi_{W'^*}(X)>_{W'}\\
&=&  \left(\sum <r,\, ->_{V'}.<s,\, ->_{W'}\right) \circ
\Phi_*(X).
 \end{eqnarray*}
\end{proo}

\begin{cor}
The black-square product defined here is equal to the one defined
in \cite{EFG} and in \cite{L2}.
\end{cor}

\begin{Rq} The white square-product is equal to $\Po\Sq
\Qo:=\F(V'\ot W' \ot k[\Sy_2])/\big( ( \Phi^{-1}(R'\ot \F(W') +
\F(V')\ot S'))\ot k[\Sy_3]   \big)$ and the black square-product
to $\Po\, \blacksquare\, \Qo:=\F(V'\ot W' \ot k[\Sy_2])/\big(
(\Psi(R'\ot S'))\ot k[\Sy_3]   \big)$. The definition of $\square$
proposed by K. Ebrahimi-Fard and L. Guo in \cite{EFG} corresponds
to $\Psi(R'\ot \F(W') + \F(V')\ot S')$ instead of $\Phi^{-1}(R'\ot
\F(W') + \F(V')\ot S')$. We have $\Phi^{-1}(R'\ot \F(W') +
\F(V')\ot S')\subset \Psi(R'\ot \F(W') + \F(V')\ot S')$. But the
second module can be slightly bigger than the first one (see the
example of $\D \Sq \D$ on page $309$ of \cite{EFG}). This explains
why the white square-product defined in \cite{EFG} is not the
Koszul dual of the black square-product.
\end{Rq}

With the explicit form of the black square-product, we get the
following property which is Proposition $2.4$ of \cite{L2}.

\begin{pro}[\cite{L2}]
For two binary quadratic regular operads $\Po$ and $\Qo$, there
exists a canonical epimorphism $\Po \BlackSq \Qo \epi \Po \Sq
\Qo$.
\end{pro}

\begin{proo}
We have to show that $\Phi \circ \Psi \, (R' \ot S')\subset R'\ot
\F(W')(3)+ \F(V')(3)\ot S'$. Let $r\ot s$ be an elemental tensor
of $R'\ot S'$. Denote $r\ot s=(r_1+r_2) \ot (s_1 +s_2)$. From
Proposition~\ref{BlackSquare}, we get $\Phi \circ \Psi(r\ot
s)=r_1\ot s_1 - r_2\ot s_2= (r_1+r_2)\ot s_1 - r_2\ot(s_1+s_2)\in
R'\ot \F(W')(3)+ \F(V')(3)\ot S'$.
\end{proo}

The proposition means that  any $\Po \Sq \Qo$-algebra is a $\Po
\BlackSq \Qo$-algebra. This result together with
Proposition~\ref{TensorVSSquare}, gives the following corollary
(Proposition 3.3 of \cite{EFG}).

\begin{cor}[\cite{EFG}]
For any $\Po$-algebra $A$ and $\Qo$-algebra $B$, their tensor
$A\ot B$ is a $\Po \BlackSq \Qo$-algebra.
\end{cor}

\begin{Rq}
The operads $\chi^+$ and $\chi^-$ discovered by J.-L. Loday in
\cite{L} factors this projection
$$\xymatrix{ \Dend \BlackSq \D \ar@{->>}[r] \ar@{<->}@/_2pc/[rr]_{!}
& \chi^{\pm} \ar@(r,u)_{!} \ar@{->>}[r] & \Dend \Sq \D.}$$ Going
from the left to the right, there is one more relation each time.
(The dimensions of the spaces of relations is $15$, $16$ and $17$
respectively).
\end{Rq}

\draftnote{est-ce-que $\chi$ correspond a $\Psi( R'\ot \F(W')(3)+
\F(V')(3)\ot S')$ ? Est-ce-que cette construction est autoduale ?
Est-ce-que l'on peut toujours factoriser $\BlackSq \epi \Sq$ par
elle ?}

\subsection{Adjunction}\label{AdjonctionReguliere}

We can apply the same methods as in Section~\ref{SymMonCat} to
prove the same kind of adjunction for black and white square
products for regular operads. Consider the category of
\emph{$k$-ary quadratic regular operads} denoted by
\emph{$k$.q-Reg}. One can extend black and white square products
in this category. Recall from \cite{Gnedbaye} that a \emph{totally
associative $k$-ary algebra} is a module equipped with a regular
$k$-ary operation such that all the quadratic compositions are
equal. Denote the corresponding operad by $T\As^{<k>}$. Dually, a
\emph{partially associative $k$-ary algebra} is a module equipped
with a regular $k$-ary operation such that the sum of all
quadratic compositions is zero. Denote the corresponding operad by
$P\As^{<k>}$. Gnedbaye proved that these two operads are Koszul
dual to each other.

\begin{pro}
The black and white square products endow the category of $k$-ary
quadratic regular operads with a structure of symmetric monoidal
category, where the operad $P\As^{<k>}$ is the unit object for
$\BlackSq$ and the operad $T\As^{<k>}$ is the unit object for
$\Sq$.
\end{pro}

\begin{thm}\label{Adjunction}
There is a natural isomorphism
$\Hom_{k.\textrm{q.Reg}}(\Po\BlackSq \Qo,\, \mathcal{R})\cong
\Hom_{k.\textrm{q.Reg}}(\Po,\, \Qo^!\Sq \mathcal{R})$.
\end{thm}

\begin{pro}\label{Operations}
For every $k$-ary quadratic regular operad $\Po$, there is a
canonical morphism of operads $P\As^{<k>} \xrightarrow{i} \Po^!
\Sq \Po$, defined by the commutative diagram
$$\xymatrix{P\As^{<k>} \ar[rr]^{l} \ar[dr]_{i} & &\Po^!\ot \Po \\
&\Po^! \Sq \Po \ar[ur]_{\bar{\Phi}}.  &} $$
\end{pro}

This proposition can be seen as a refinement of
Proposition~\ref{universalLiek}. When $\Po$ is a $k$-ary regular
quadratic operad, the map $\Li^{<k>} \to \Po^! \ot \Po$ factors
through $P\As^{<k>}$, where the morphism $\Li^{<k>} \to
P\As^{<k>}$ is induced by the anti-symmetrization of the $k$-ary
partially associative product as in the binary case.

\subsection{Non-symmetric cohomology operations}

In this section, we refine the arguments of
Section~\ref{Cohomology Operations} for non-symmetric (regular)
operads. This gives non-vanishing natural operations on the
deformation chain complex of any algebras over such operads. More precisely, we prove that,
under some assumptions, the (co)chain complex defining the cohomology of algebras is a multiplicative operad.\\

Recall from \cite{GerstenhaberVoronov} that an \emph{operad with
multiplication} is a non-symmetric operad $\Po$ endowed with a
morphism $\A \to \Po$. Let $\Po$ be a finitely generated binary
non-symmetric Koszul operad. Following Section~\ref{Cohomology
Operations}, the chain complex defining the cohomology of a
$\Po$-algebra $A$ is equal to
$C^\bullet_\Po(A)=\Hom^{\bullet}_\KK(\Po^{\ac}, \End(A))$ which is
a non-symmetric (convolution) operad. By
Proposition~\ref{Operations}, there is a morphism of operads $\A
\to \Po^! \ot \Po$. Since $\Po^{\ac}={\Po^!}^*$, we have
$$\A \to \Po^! \ot \Po \cong \Hom_\KK(\Po^{\ac}, \Po) \xrightarrow{\Phi_*} \Hom_\KK(\Po^{\ac}, \End(A)).$$
These results form the following proposition.

\begin{pro}\label{multiplicative Operad}
For every finitely generated binary non-symmetric Koszul operad
$\Po$ and every $\Po$-algebra $A$, the chain complex defining its
cohomology $C^\bullet_\Po(A)$ is an operad with multiplication.
\end{pro}

The multiplication $\As \to \Po$ of an operad $\Po$ allows us to
define a canonical cosimplicial structure on it (see
\cite{McClureSmith02} Section $3$) and then a differential map $d$
by alternate summation (see \cite{GerstenhaberVoronov} Formula
$(5)$). Denote by $m$ the image of the associative operation. The
face maps $d^i\,:\, \Po(n) \to \Po(n+1)$ are defined by
$$ d^i(p) := \left\{ \begin{array}{lll}
m \circ_2 p & \textrm{if} &  i=0\\
p \circ_i m & \textrm{if} &  0<i<n+1\\
m \circ_1 p & \textrm{if} &  i=n+1.
\end{array}              \right.   $$

The differential $d$ is equal to $d (f) := m \star f - (-1)^{|f|}f
\star m=[m,f]$.

\begin{lem}
With the same assumptions, the differential $\partial(f)$ on
$C^\bullet_\Po(A)$ is equal to $(-1)^{|f|}d(f)$. Hence, the chain
complex $C^\bullet_\Po(A)$ is always cosimplicial.
\end{lem}

\begin{proo}
The image of the associative operation in $\Hom_\KK(\Po^{\ac},
\End(A))$ is the map $\alpha \, :\, \Po^{\ac}(2) \to \Hom(A^{\ot
2}, A)$ defined in Section~\ref{Cohomology Operations}. We prove
in Lemma~\ref{Formula} that the differential on $C^\bullet_\Po(A)$
is equal to $\partial(f)=[f, \alpha]=(-1)^{|f|}d(f)$.
\end{proo}

Therefore, the chain complex $C_\Po^\bullet(A)$ is endowed with
two types of operations : \emph{braces} operations induced by the
non-symmetric operadic structure and an associative operation
called the \emph{cup product} coming from the properties of
Manin's products. In \cite{GerstenhaberVoronov}, M. Gerstenhaber
and A.A. Voronov defined the notion of \emph{homotopy $G$-algebra}
which gives the compatibility between these types of operations.
Their purpose was to describe the operations acting of the chain
complex of Hochschild cohomology of an associative algebra.
Actually, the structure of homotopy $G$-algebra on the deformation
chain complex and the structure of Gerstenhaber algebra on
cohomology is universal among finitely generated binary
non-symmetric Koszul operads.

\begin{cor}
For every finitely generated binary non-symmetric Koszul operad
$\Po$ and every $\Po$-algebra $A$, the chain complex
$C_\Po^\bullet(A)$ is a homotopy $G$-algebra and the cohomology
space $H^\bullet_\Po(A)$ is a Gerstenhaber algebra.
\end{cor}

\begin{proo}
Apply Theorem $3$ of \cite{GerstenhaberVoronov} which asserts that
any multiplicative operad induces a homotopy $G$-algebra on the
direct sum of its components. To prove the second part, apply the
computations of the proof of Corollary $5$ of
\cite{GerstenhaberVoronov}.
\end{proo}

\subsection{Generalized Deligne's conjecture}

Finally, we extend and prove Deligne's conjecture to any algebra
over a finitely generated binary non-symmetric Koszul
operad, which includes the original case of associative algebras.\\


The little disk operad $\mathcal{D}_2$ is a topological operad
defined by configurations of disks on the plane. In 1976, F. Cohen
showed that the homology operad $H_\bullet(\mathcal{D}_2)$ is
equal to the operad coding Gerstenhaber algebras \cite{FredCohen}.
This led P. Deligne to make the following wish "I would like the
complex computing Hochschild cohomology to be an algebra over [the
singular chain operad of the little disks] or a suitable version
of it" in \cite{Deligne93}. By suitable version of it, he meant
another operad homotopically equivalent to $\mathcal{D}_2$. This
conjecture can be seen as a lifting on the level of chain
complexes of the result of F. Cohen. In 1999, J.E. McClure and
J.H. Smith gave a prove of this conjecture in the following way.
First, they construct a topological operad $\mathcal{C}$ whose
chain version acts on any multiplicative operad. Then, they show
that this operad is equivalent to the little disks operad. This
proof with Proposition~\ref{multiplicative Operad} shows that
Deligne conjecture can be generalized to any finitely generated
binary non-symmetric Koszul operads and is not specific to the
case of associative algebras.

\begin{thm}\label{DeligneConjecture}
For every finitely generated binary non-symmetric Koszul operad
$\Po$ and every $\Po$-algebra $A$, the chain complex
$C^\bullet_\Po(A)$ is an algebra over an operad equivalent to the
singular chains of the little disks operad.
\end{thm}

\begin{proo}
Since $C^\bullet_\Po(A)$ is an operad with multiplication, the
operad $\mathcal{C}$ of \cite{McClureSmith02} acts on it. And this
operad is weakly equivalent to the little disks operad by Theorem
$3.3$ of \cite{McClureSmith02}.
\end{proo}

Notice that the non-symmetric case is very different from the
symmetric one. The Lie bracket $\{\,  , \}$ described in
Section~\ref{Cohomology Operations} vanishes on cohomology. When
the algebra is modelled by a non-symmetric operad, this Lie
bracket is the symmetrization of an associative operation, the
cup-product, which
is not necessarily trivial on cohomology.\\

In \cite{Markl-Operations}, M. Markl defined the notion of
\emph{natural operations} on cohomology and asked a few questions
and conjectures about the operad $\mathcal{B}_\Po$ generated by
these operations. Here we have proved that, for binary
non-symmetric Koszul operads, the Gerstenhaber operad imbeds into
$\mathcal{B}_\Po$. For more precise statements depending on the
operad $\Po$, one has to work with $\Po^!\circ \Po$. In the
symmetric case, operations of $\Po^!\circ \Po$ could give
non-trivial operations in cohomology.

\begin{Rq}
In \cite{Yau}, D. Yau proved this generalized Deligne's conjecture
for a few operads found by J.-L. Loday. His method is based on a
notion of \emph{pre-operadic system} which ensures that
$C^\bullet_\Po(A)=\Hom^{\bullet}_\KK(\Po^{\ac}, \End(A))$ is an
operad. Actually this notion comes from the axioms of a basis for
the Koszul dual cooperad. The cohomology space of an algebra over
any non-symmetric Koszul operad is always a non-symmetric operad
(convolution operad from the Koszul dual cooperad to the
endomorphism operad). Then, the author shows, case by case, that
the cohomology of the operad is multiplicative. In fact, the
adjunction of Manin's products for non-symmetric binary operads
always provides a morphism $\A \to C^\bullet_\Po(A)$.
\end{Rq}

\subsection{The operad $\Qu$ and its Koszul dual}

In this section, we study the example of black square-product
$\Qu=\Dend \BlackSq \Dend$ introduced by M. Aguiar and J.-L. Loday
in \cite{AL}. We prove that the Koszul dual of $\Qu$
is the operad $\Qu^!=\Pe \ot \D=\Pe \circ \D$. \\

The operad $\Dend$ is a split of one associative product $\star$
into two products $\prec$ and $\succ$, $\star=\prec + \succ$. The
operad $\Qu$ was defined by M. Aguiar and J.-L. Loday in \cite{AL}
as a split of an associative product $\star$ into four products
$\nearrow$, $\searrow$, $\swarrow$ and $\nwarrow$, that is
$\star=\nearrow+\searrow+\swarrow+\nwarrow$. It was proved in
\cite{EFG} that this operad $\Qu$ is equal to the black
square-product $\Dend \BlackSq \Dend$. Therefore one can interpret
the splitting of associativity with the black square-product with
$\Dend$. At the end of their paper, M. Aguiar and J.-L. Loday
raised one question ``what is the Koszul dual of the operad $\Qu$
?'' and two conjectures. The first conjecture deals with the
dimensions of the $\Sy_n$-modules $\Qu(n)$ and the second one is
that the operad $\Qu$ is Koszul. In the rest of this section, we
answer
these questions.\\

The previous section give a direct answer to the first question.

\begin{pro}
The Koszul dual of $\Qu$ is equal to $\Qu^!=\D \Sq \D$.
\end{pro}

\begin{proo}
Since the Koszul dual of $\Dend$ is the operad $\Dend^!=\D$
(\cite{L3} Proposition $8.3$), we have
$$\Qu^!=(\Dend \BlackSq \Dend)^!=\Dend^! \Sq
\Dend^!=\D \Sq \D,$$ from Proposition~\ref{BlackSquareOperad}.
\end{proo}

It remains to use the explicit form of the white square-product to
describe $\D \Sq \D$.

\begin{thm}
\label{DualQuad} The operad $\Qu^!=\D \Sq \D$ is isomorphic to
$\Pe \ot \D=\Pe\ot \Pe \ot \A$.
\end{thm}

\begin{proo}
Denote the basis of $\Dend'(2)$ by $\prec.k \, \oplus \succ.k$ and
its dual basis, the basis of $\D'(2)$, by $\dashv.k \, \oplus
\vdash.k$. The induced basis of $\Qu'(2)=(\Dend \BlackSq
\Dend)'(2)$ is $\{ \prec \ot \prec,\, \prec\ot \succ,\,
\succ\ot\prec,\, \succ\ot \succ \}$ and the induced basis of
${{\Qu}^!}{}'(2)=(\D \Sq \D)'(2)$ is $\{ \dashv \ot \dashv,\,
\dashv\ot \vdash,\, \vdash \ot\dashv,\, \vdash\ot \vdash \}$ . The
relations of $\D$ are easy to remember. Represent the operation
$\dashv$ by the tree $\vcenter{\xymatrix@M=0pt@R=6pt@C=6pt{
\ar@2{-}[dr] &  &\ar@{-}[dl]  \\
 &\ar@{-}[d] &  \\  & &}}$ and the operation $\vdash$ by the tree
 $\vcenter{\xymatrix@M=0pt@R=6pt@C=6pt{
\ar@{-}[dr] &  &\ar@2{-}[dl]  \\
 &\ar@{-}[d] &  \\  & &}}$. Any element of $\F(\dashv.k
\, \oplus \vdash.k)(3)$ can be seen as a tree with exactly one
path from one leaf to the root. For example, the composition
$\dashv \circ_1 \vdash$ corresponds to the tree
$\vcenter{\xymatrix@M=0pt@R=6pt@C=6pt{  \ar@{-}[dr] &
&\ar@2{-}[dl] & &   \\
&\ar@2{-}[dr]  &   & \ar@{-}[dl] & \\
& &\ar@{-}[d] & & \\
& & }}$. To get the relations of $\D$, identify the trees with
paths from the same leaf. For instance, we have $\dashv \circ_1
\vdash\, =\, \vdash \circ_2 \dashv$, which corresponds to
$\vcenter{\xymatrix@M=0pt@R=6pt@C=6pt{ \ar@{-}[dr] &
&\ar@2{-}[dl] & &   \\
&\ar@2{-}[dr]  &   & \ar@{-}[dl] & \\
& &\ar@{-}[d] & & \\
& & }}=\vcenter{\xymatrix@M=0pt@R=6pt@C=6pt{ &
&\ar@2{-}[dr] & & \ar@{-}[dl]  \\
& \ar@{-}[dr] & &\ar@2{-}[dl]  & \\
& &\ar@{-}[d] & & \\
& & }} $. The relations of $\D$ are
$$ \left\{ \begin{array}{lc}
\dashv \circ_1 \dashv\, =\, \dashv \circ_2 \dashv =\dashv \circ_2
\vdash & (L)\\
\dashv \circ_1 \vdash\, =\, \vdash \circ_2 \dashv & (M)\\
\vdash \circ_1 \dashv\, =\, \vdash \circ_1 \vdash = \vdash \circ_2
\vdash & (R),
\end{array}
\right. $$ where the first line corresponds to the Left leaf, the
second one to the Middle leaf and the last one to the Right leaf.
For simplicity, denote these compositions and relations by
$$ \left\{ \begin{array}{l}
 L_1=L'_2=L''_2 \\
M_1=M_2\\
R'_1=R''_1=R_2.
\end{array}
\right. $$ The operad $\D$ is equal to $\F(V'\ot k[\Sy_2])/(R'\ot
k[\Sy_3])$. One can see that the following relations are elements
of $(R\ot \F(V')+\F(V')\ot R)\cap \textrm{Im}\, \Phi$
$$ \left\{ \begin{array}{ll}
(L,L) & L_1\ot L_1=L'_2\ot L'_2=L''_2\ot L'_2=L'_2\ot
L''_2=L''_2\ot L''_2 \\
(L,M) & L_1 \ot M_1 = L'_2\ot M_2 =  L''_2\ot M_2 \\
(L,R) & L_1 \ot R'_1 = L_1 \ot R''_1 = L'_2\ot R_2 = L''_2\ot R_2
\\
(M,M) & M_1\ot M_1 = M_2 \ot M_2
\end{array} \right. $$
The other ones are obtained by the symmetries
$\vcenter{\xymatrix@M=0pt@R=6pt@C=6pt{ \ar@{-}[dr] &
&\ar@{-}[dl] & &   \\
&\ar@{-}[dr]  &   & \ar@{-}[dl] & \\
& &\ar@{-}[d] & & \\
& & }} \leftrightarrow \vcenter{\xymatrix@M=0pt@R=6pt@C=6pt{ &
&\ar@{-}[dr] & & \ar@{-}[dl]  \\
& \ar@{-}[dr] & &\ar@{-}[dl]  & \\
& &\ar@{-}[d] & & \\
& & }} $ and $a\ot b \leftrightarrow b\ot a$. We get $23$ linearly
independent elements in $(R\ot \F(V')+\F(V')\ot R)\cap
\textrm{Im}\, \Phi$. Since the dimension of $\Qu'(3)$ is
$23=32-9$, we know that these elements form a basis of $(R\ot
\F(V')+\F(V')\ot R)\cap \textrm{Im}\, \Phi$. Hence, they give the
relations defining $\Qu^!$.

Interpret these relations in the same way as the ones of $\D$. An
element of $\F\big((\dashv.k \, \oplus \vdash.k) \ot (\dashv.k \,
\oplus \vdash.k)\big)(3)$ can be seen as a tree with two kind of
paths, one given by the left side of $\ot$ and the second one by
the right side of $\ot$. For instance, the tree
$\vcenter{\xymatrix@M=0pt@R=6pt@C=6pt{ \ar@{.}[dr] &
&\ar@2{-}[dl] & &   \\
&\ar@2{-}[dr]  &   & \ar@{.}[dl] & \\
& &\ar@{-}[d] & & \\
& & }}$ represents $(\dashv \ot \vdash) \circ_1 (\vdash \ot
\dashv)$, where the left side corresponds to $=$ and the right
side to $\cdots$. This produces two indexes for the leaves. With
this identification, the relations of $\Qu^!$ mean that any
elements written with trees such that the same leaves are indexed
by the same ``colors'' are equal. Therefore, a basis for
$\Qu^!(n)$ is given by planar corollas with $n$ leaves indexed by
two colors. The composition of such trees is a corollas and to
know which leaf is indexed by which color, follow the path of the
same color. As a consequence, we have $\D \Sq \D =\Pe \ot \D$. (A
basis for the operad $\D$ is given by corollas with one leaf
emphasized. Tensoring with $\Pe$ induces another independent index
of the leaves. And the compositions are the same.)
\end{proo}

\begin{Rq}
More generally, we get the duals of the operad coding
octo-algebras of P. Leroux \cite{Le} and its follow-up. Since
$\mathcal{O}cto=\Dend ^{\BlackSq 3}$, we get $\mathcal{O}cto^!=\Pe
^{\ot 3} \ot \A$ and, for any $n\in \NN$, $\left( \Dend^{\BlackSq
n}\right)^!=\D^{\Sq n}= \Pe^{\ot n}\ot \A$.
\end{Rq}

\begin{cor}
\label{dimQuad!} The dimensions of the components of the Koszul
dual of $\Qu$ are equal to $$\textrm{dim}\, (\Qu^!(n))=n^2.n!.$$
\end{cor}

\begin{proo}
We have $\textrm{dim} (\Qu^!(n))=\textrm{dim} (\Pe(n) \ot_k
\D(n))=n^2.n!$.
\end{proo}

Proposition~\ref{PermWhiteAs} gives that $\D \Sq \D = \Pe \ot \D=
\Pe \circ \D=\Pe \ot \Pe \ot \A=\Pe \circ \Pe \circ \A$. We have
$\D=\Sigma(\Pe)$ and $\D \Sq \D = \Sigma\big( \Pe \circ \Pe \big)=
\Sigma \big( \Pe^{\ot 2} \big)$. By duality, we get another form
for $\Qu$.

\begin{cor}
We have $\Qu=\Pe \bullet \Pe \bullet \A$
\end{cor}

\subsection{Koszulity of $\Qu$ and other operads defined by square products}
Aguiar and Loday made in \cite{AL} the conjecture that the operad
$\Qu$ is Koszul. We show this statement using poset's method of
\cite{V2}. More generally, we prove that the operads of the form
 $\Dend^{\BlackSq n}$ and $\D^{\Sq n}$ are Koszul. P. Leroux introduced in
\cite{Leroux04} the operad $\mathcal{E}nnea=\Tridend \BlackSq
\Tridend$. We prove the same results for the family
$\Tridend^{\BlackSq n}$ and $\Trias^{\Sq n}$. All these families
provide infinitely many examples of the
generalized Deligne's conjecture proved in Theorem~\ref{DeligneConjecture}.\\

In order to study the homological properties of the algebras over
an operad, it is crucial to prove that the operad is Koszul. We
refer the reader to the paper of B. Fresse \cite{F} or to the book
of M. Markl, S. Shnider and J. Stasheff \cite{MSS} for a full
treatment of the subject. Since an operad is Koszul if and only if
its Koszul dual is Koszul, we work with the simplest one to prove
that the pair is a pair of Koszul operads. In the case of the
operads $\Qu$ and $\Qu^!=\Pe \ot \D$, we will prove that the
Koszul property holds for the last one, $\Pe \ot \D$.\\

Let $\Po$ be an algebraic operad coming from an operad in the
category of Sets. For instance, it is the case when the relations
defining the operads only involve equalities between two terms and
no linear combination. The operads $\A$, $\C$, $\Pe$ and $\D$ are
of this type. In \cite{V2}, we defined a family of partition type
posets associated to such an operad $\Po$ and proved that the
operad is Koszul over $\ZZ$ and over any field $k$
if and only each maximal intervals of the posets are Cohen-Macaulay.\\

We saw in the proof of Theorem~\ref{DualQuad} that $\Qu^!=\Pe \ot
\D = \Pe \ot \Pe \ot \As$. Therefore, $\Qu^!$ is a set operad with
basis $\{(i,\, j,\, \sigma) \ | \ 1\leq i,\,j \leq n, \, \sigma\in
\Sy_n \}$. The partitions associated to $\Qu^!$ are of the form
$\big( \sigma(1),\ldots,\, \overline{\sigma(i)},\ldots,\,
\underline{\sigma(j)} ,\ldots ,\, \sigma(n) \big)$, where $1\leq
i,\,j \leq n $ and $\sigma\in \Sy_n$. The order between the
$\Qu^!$-partitions is given by the refinement of partitions with
respect to the two indexes. For instance, we have
$\{(\overline{3},\, \underline{1},\, 5),\,
(\underline{\overline{2}},\, 4) \}\leq \{(3,\, \underline{1},\, 5,
\, \overline{2},\, 4) \}$.

\begin{lem}
For each $n\in \NN$, the maximal intervals of the poset $\Pi_{\Pe
\ot \D}(n)$ associated to the operad $\Pe \ot \D$ are totally
semi-modular.
\end{lem}

\begin{proo}
The proof is the same than Lemma $1.10$, $1.15$ and $2.6$ of
\cite{CV}.
\end{proo}

\begin{thm}
The operad $\Qu$ is Koszul over $\ZZ$.
\end{thm}

\begin{proo}
The maximal intervals of the posets $\Pi_{\Pe \ot \D}(n)$ are
totally semi-modular. Therefore they are Cohen-Macaulay over $\ZZ$
by \cite{Baclawski, Farmer} . One can see that the operad $\Pe \ot
\D$ is a basic set operad (see \cite{V2} Page 6). Then we can
apply Theorem $9$ of \cite{V2}.
\end{proo}

\begin{cor}
The dimensions of the homogenous components of $\Qu$ are equal to
$$\textrm{dim}
(\Qu(n))=(n-1)!\sum_{j=n}^{2n-1}\binom{3n}{n+1+j}.\binom{j-1}{j-n}.$$
\end{cor}

\begin{proo}
When an operad $\Po$ is Koszul, there are relations between the
dimensions of $\Po(n)$ and the dimensions of $\Po^!(n)$ (see
\cite{GK} Theorem $(3.3.2)$ or \cite{L3} Appendix B.5.c.). Use
these relations with Corollary~\ref{dimQuad!} to conclude.
\end{proo}

More generally, we have seen that, for every $n\in \NN$,
$\left(\Dend^{\BlackSq n}\right)^!=\Pe^{\ot n}\ot \A$, which is a
basic set operad. The related partitions have the same form than
the ones for $\Qu^!$ but with $n$ types of indices instead of $2$.

\begin{thm}\label{Dendn-Koszul}
For every $n$, the operad $\Dend^{\BlackSq n}$ is Koszul over
$\ZZ$.
\end{thm}

\begin{proo}
Apply the same arguments.
\end{proo}

J.-L. Loday and M. Ronco introduced in \cite{LodayRonco04} the
pair of Koszul dual operads $\Trias$ and $\Tridend$. A
$\Trias$-algebra is a $\D$-algebra with an extra operation. In
\cite{V2}, we defined a commutative analogue of $\Trias$ which we
denoted by $\Comtri$. The $\Comtri$-partitions are partitions with
at least one element of each block emphasized. The
$\Trias$-partitions are ordered partitions with at least one
element in the block emphasized. Using the same ideas than before,
we have to following results. The operad $\Trias^{\Sq n}\cong
\Comtri^{\ot n}\ot \A$. The maximal intervals of $\Pi_{\Trias^{\Sq
n}}$ are totally semi-modular.

\begin{thm}
For every $n$, the operads $\Trias^{\Sq n}$ and
$\Tridend^{\BlackSq n}$ are Koszul over $\ZZ$.
\end{thm}

Recall from \cite{EFG} Proposition $3.5$ that $\Tridend \BlackSq
\Tridend$ is isomorphic to the $\mathcal{E}nnea$ operad defined by
P. Leroux in \cite{Leroux04}. The previous theorem gives that the
$\mathcal{E}nnea$ operad is Koszul over $\ZZ$.\\

From this result, we get four infinite families of operads for
which Deligne's conjecture holds.

\begin{cor}
Let $\Po$ be an operad of the form $\Dend^{\BlackSq n}$, $\D^{\Sq
n}$, $\Tridend^{\BlackSq n}$ or $\Trias^{\Sq n}$. Then for any
$\Po$-algebra $A$, the chain complex $C^\bullet_\Po(A)$ is an
algebra over an operad equivalent to the singular chains of the
little disks operad.
\end{cor}

\begin{proo}
By the previous theorems, these operads are finitely generated
binary non-symmetric and Koszul. Then apply
Theorem~\ref{DeligneConjecture}.
\end{proo}

Notice that the poset's method of \cite{V2} allowed us to prove
that these operads are Koszul over $\ZZ$. Since the proof of
Deligne's conjecture \cite{McClureSmith02} also works over the
ring of integers, this last corollary holds over $\ZZ$.

\appendix

\section{Associative algebras, operads and properads}
\label{Alg-Op-Properad}

This appendix is a short survey on the notions of associative
algebras, operads and properads which are the main examples of
$2$-monoidal categories treated in this text. For a complete
treatment of the subject, we refer the reader to \cite{V1}.

\subsection{Associative algebras}

Associative algebras, operads and properads are monoids in some
monoidal categories.\\

Let $k$ be the ground field and let $(\KK\textrm{-Mod}, \otk,\,
\KK)$ be the monoidal category of $\KK$-modules equipped with the
tensor product over $\KK$.
\begin{dei}[Associative algebra]
A monoid $(A,\, \mu,\, \eta)$ in $(\KK\textrm{-Mod}, \otk,\, \KK)$
is an \emph{associative algebra}. The product $\mu \, : \, A\otk A
\xrightarrow{\mu} A $ is associative and $\KK \xrightarrow{\eta}
A$ is the unit of $A$.
\end{dei}

\begin{figure}[h]
$$\xymatrix@=14pt{ \ar[d] \\ *+[F-,]{a_1} \ar[d] \\ \vdots \ar[d]\\ *+[F-,]{a_l} \ar[d] \\ {} }$$
\caption{Product of $a_1,\, \ldots,\, a_l$.} \label{Fig1}
\end{figure}
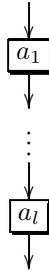

The product of elements $a_1,\ldots ,\, a_l$ of $A$ can be
represented by an indexed branch, see Figure~\ref{Fig1}

\begin{Ex}
Let $M$ be a $\KK$-module. Denote by $\End(M):=\Hom_\KK(M,\,M)$
the space of endomorphisms of $M$. With the composition of
endomorphisms, $\End(M)$ is an associative algebra.
\end{Ex}

An  \emph{associative coalgebra} is a comonoid in
$(\KK\textrm{-Mod}, \otk,\, \KK)$, that is a monoid in the
opposite category.

\subsection{Operads}

An \emph{$\Sy$-module} is a collection $\{\Po(n) \}_{n\in \NN^*}$
of right modules over the symmetric group $\Sy_n$. In the category
of $\Sy$-modules, one defines a monoidal product by the following
formula

$$\Po \circ \Qo (n) := \bigoplus_{1\leqslant l \leqslant n} \left(
\bigoplus_{i_1+\cdots+i_l=n} \Po(l) \otimes_k \big( \Qo(i_1)
\otimes_k \cdots \otimes_k \Qo(i_l)\big) \otimes_{\Sy_{i_1}\times
\cdots \times \Sy_{i_l}} k[\Sy_n]  \right)_{\Sy_l},$$

where the coinvariants are taken with respect to the action of the
symmetric group $\Sy_l$ given by $(p\otk q_1 \ldots q_l\otk
\sigma)^\nu := p^\nu \otk q_{\nu(1)} \ldots q_{\nu(l)} \otk
\bar{\nu}^{-1}.\sigma$ for $p\in \Po(l)$, $q_j\in \Qo(i_j)$,
$\sigma\in \Sy_n$ and $\nu\in \Sy_l$, such that $\bar{\nu}$ is the
induced block permutation.\\

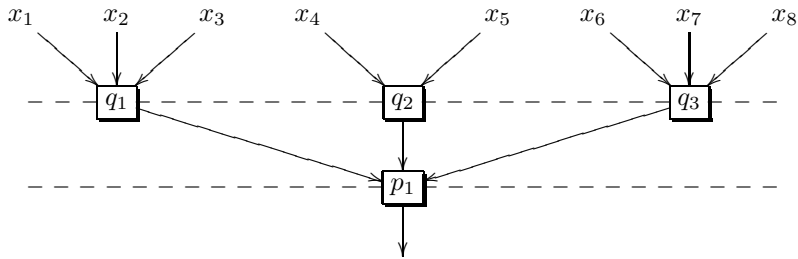
\begin{figure}[h]
$$\xymatrix@R=20pt@C=20pt{ x_1 \ar[dr] & x_2 \ar[d] & x_3 \ar[dl]& x_4\ar[dr]& &x_5\ar[dl]
& x_6\ar[dr] & x_7\ar[d] & x_   8\ar[dl] \\
\ar@{--}[r] & *+[F-,]{q_1}\ar@{--}[rrr] \ar[drrr]& & &
*+[F-,]{q_2}\ar@{--}[rrr] \ar[d]
 & & & *+[F-,]{q_3} \ar@{--}[r]
\ar[dlll]& \\
\ar@{--}[rrrr] & & & & *+[F-,]{p_1} \ar[d]\ar@{--}[rrrr] & & & & \\
 & & & & & & & &  }$$
 \caption{The monoidal product $\Po \circ \Qo$.}
\label{Fig2}
\end{figure}

The notion of $\Sy$-module is used to model the multi-linear
operations acting on some algebras. The monoidal product $\circ$
reflects the compositions of operations and can be represented by
2-levelled trees whose vertices are indexed by the elements of
$\Po$ and $\Qo$, see Figure~\ref{Fig2}.
The unit of this monoidal category is given by the $\Sy$-module
$I=(k,\, 0,\, 0,\, \ldots )$, which corresponds to the identity
operation represented by $|\,$.

\begin{dei}[Operad]
A monoid $(\Po, \mu, \eta)$ in $(\Sy\textrm{-Mod}, \circ,\, I )$
is called an \emph{operad}. The associative product $\mu \, : \,
\Po \circ \Po \to \Po$ is called the \emph{composition product}
and $\eta \, :\, \II \to \Po$ is the \emph{unit}.
\end{dei}

\begin{Ex}
Let $M$ be a $\KK$-module and consider $\End(M):=\bigoplus_{n\in
\NN^*} \Hom_\KK(M^{\ot n},\, M)$. The permutation of the inputs of
a morphism in $\Hom_\KK(M^{\ot n},\, M)$ makes $\End(M)$ into an
$\Sy$-module. With the natural composition of morphisms, $\End(M)$
is an operad called the \emph{endomorphism operad}.
\end{Ex}

A \emph{cooperad} is a comonoid in $(\Sy\textrm{-Mod}, \circ,\, I
)$.\\

To every $\KK$-module $V$, one can associate an $\Sy$-module
$\widetilde{V}:=(V,\,  0 ,\,0,\, \ldots)$ concentrated in arity
$1$. This defines an embedding of $\KK$-Mod into $\Sy$-Mod. One
can check that this embedding is compatible with the monoidal
products, that is $\widetilde{V\ot W}= \widetilde{V} \circ
\widetilde{W}$. Therefore, $(\KK\textrm{-Mod}, \otk,\, \KK)$ is a
full monoidal subcategory of $(\Sy\textrm{-Mod},\, \circ,\, I)$.
 Thus every associative algebra is an operad.\\

We can forget the action of the symmetric groups and work in the
category of $\NN^*$-graded vector spaces. This category is endowed
with a monoidal product, a non-symmetric analog of the previous
one. We still denote it by $\circ$  :
$$\Po \circ \Qo (n) := \bigoplus_{1\leqslant l \leqslant n} \left(
\bigoplus_{i_1+\cdots+i_l=n} \Po(l) \otimes_k \big( \Qo(i_1)
\otimes_k \cdots \otimes_k \Qo(i_l)\big) \right).$$

\begin{dei}[Non-symmetric operad]
 A monoid in
this monoidal category is called a \emph{non-symmetric operad}.
\end{dei}

We can also define a notion of operads with colors indexing the
inputs and the outputs. The composition of such operations have to
fit with the colors. Such operads are called
\emph{colored operads} (cf. \cite{VdL, BM}).\\

\subsection{Properads}

\label{properads} 
\begin{figure}[h!]
$$ \xymatrix{
x_1 \ar[dr] &x_2 \ar[d] &x_ 3\ar[dl] &x_4 \ar[dr] & &x_5 \ar[dl] \\
\ar@{--}[r]& *+[F-,]{q_1} \ar@{-}[dl] \ar@{-}[drr] \ar@{--}[rrr]&
& & *+[F-,]{q_2}
\ar@{-}[dr]  \ar@{-}[dll]|(0.75) \hole  \ar@{--}[r]& \\
 *=0{} \ar[dr]& &*=0{}\ar[dl] &*=0{}\ar[dr]& &*=0{} \ar[dl]\\
 \ar@{--}[r] & *+[F-,]{p_1} \ar[d]
\ar@{--}[rrr]& & & *+[F-,]{p_2} \ar[dl] \ar[d] \ar[dr]
\ar@{--}[r] & \\
 & y_1 & &y_2  &y_3 &y_4 } $$
\caption{Composition of operations with multiple inputs and
multiple outputs.} \label{Fig3}
\end{figure}
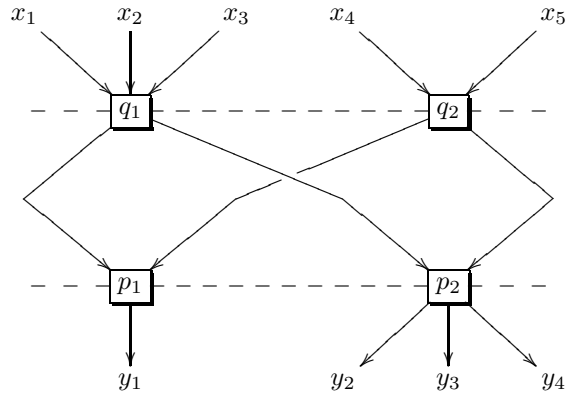

We are going to pursue this generalization.
Elements of an associative algebra can be seen as operations with
one input and one output (see Figure~\ref{Fig1}). Elements of an
operad represent operations with multiple inputs but one output.
To model operations with multiple inputs and multiple outputs, one
uses the notion of \emph{$\Sy$-bimodule}. An \emph{$\Sy$-bimodule}
is a collection $\{\Po(m, n) \}_{m, \, n\in \NN^*}$ of modules
over the symmetric groups $\Sy_n$ on the right and $\Sy_m$ on the
left. In this category, we define a monoidal product based on the
composition of operations indexing the vertices of a 2-levelled
directed connected
graph, see Figure~\ref{Fig3}.\\

Let $a$ and $b$ be the number of vertices on the first and on the
second level respectively. Let $N$ be the number of edges between
the first and the second level. To an $a$-tuple of integers
$\oi:=(i_1,\, \ldots ,\, i_a)$, we associate
$|\oi|:=i_1+\cdots+i_a$. Given two $a$-tuples $\oi$ and $\oj$, we
denote by $\Qo(\oj,\, \oi)$ the tensor product $\Qo(j_1,\,
i_1)\otimes \cdots \otimes \Qo(j_a,\, i_a)$ and we denote by
$\Sy_{\oi}$ the image of the direct product of the groups
$\Sy_{i_1}\times\cdots \times \Sy_{i_n}$ in $\Sy_{|\oi|}$.

\begin{dei}[Connected permutations]
Let $N$ be an integer. Let $\ok=(k_1,\, \ldots,\, k_b)$ be a
$b$-tuple and $\oj=(j_1,\, \ldots ,\, j_a)$ be a $a$-tuple such
that $|\ok| = k_1+\cdots+k_b =|\oj|=j_1+\cdots+ j_a= N$.

A \emph{$(\ok,\, \oj)$-connected permutation} $\sigma$ is a
permutation of $\Sy_N$ such that the graph of a geometric
representation of $\sigma$ is connected if one gathers the inputs
labelled by $j_1+\cdots+ j_i +1, \,
 \ldots,\, j_1+\cdots +j_i + j_{i+1}$, for $0\leq i \leq a-1$, and
 the outputs labelled by $k_1+\cdots+ k_i +1, \,
 \ldots,\, k_1+\cdots +k_i + k_{i+1}$, for $0\leq i \leq b-1$.
The set of $(\ok,\, \oj)$-connected permutations is denoted by
$\Sc_\kj \,$.
\end{dei}

\begin{Ex}
Consider the permutation $(1324)$ in $\Sy_4$ and take  $\ok=(2,
2)$ and $\oj=(2, 2)$. If one links the inputs $1$, $2$ and $3$,
$4$ and the outputs $1$, $2$ and $3$, $4$, it gives the following
connected graph
$$\xymatrix{ *+[o][F-]{1}  \ar@{-}[r] \ar@{-}[d] & *+[o][F-]{2}
\ar@{-}[dr]& *+[o][F-]{3}  \ar@{-}[dl] |\hole \ar@{-}[r]
& *+[o][F-]{4}  \ar@{-}[d] \\
 *+[o][F-]{1}  \ar@{-}[r] &*+[o][F-]{2}  & *+[o][F-]{3}  \ar@{-}[r] & *+[o][F-]{4}} $$
Therefore, the permutation $(1324)$ is $((2, 2),\, (2,
2))$-connected.
\end{Ex}

Let $\Po$ and $\Qo$ be two $\Sy$-bimodules, their monoidal product
is given by the formula
$$\mathcal{P}\boxtimes_c \mathcal{Q}(m,\, n) :=
 \bigoplus_{N\in \mathbb{N}^*} \left( \bigoplus_{\ol,\, \ok,\, \oj,\, \oi} \KK[\mathbb{S}_m]
 \otimes_{\mathbb{S}_\ol}
\mathcal{P}(\ol,\, \ok)\otimes_{\mathbb{S}_{\ok}} k[\Sc_\kj]
\otimes_{\mathbb{S}_\oj} \mathcal{Q}(\oj,\, \oi)
\otimes_{\mathbb{S}_\oi} k[\mathbb{S}_n]
\right)_{\Sy_b^{\textrm{op}}\times \Sy_a} ,$$ where the second
direct sum runs over the $b$-tuples $\ol$, $\ok$ and the
$a$-tuples $\oj$, $\oi$ such that  $|\ol|=m$, $|\ok|=|\oj|=N$,
$|\oi|=n$ and where the coinvariants correspond to the following
action of $\Sy_b^{\textrm{op}}\times \Sy_a$ :
\begin{eqnarray*}
&&\theta \otimes p_1\otimes \cdots \otimes p_b \otimes \sigma
\otimes q_1 \otimes \cdots \otimes q_a \otimes\omega
  \sim \\
&&\theta \,\tau^{-1}_\ol \otimes p_{\tau^{-1}(1)}\otimes \cdots
\otimes p_{\tau^{-1}(b)} \otimes
 \tau_\ok\,
 \sigma \, \nu_\oj \otimes q_{\nu(1)} \otimes \cdots \otimes q_{\nu(a)}
 \otimes \nu^{-1}_{\oi} \,
 \omega,
\end{eqnarray*}
for $\theta \in \Sy_m$, $\omega \in \Sy_n$, $\sigma \in \Sc_\kj$
and for $\tau \in \mathbb{S}_b$ with $\tau_\ok$ the corresponding
block permutation, $\nu \in \mathbb{S}_a$ and $\nu_\oj$ the
corresponding block permutation. The unit $\II$ for this monoidal
product is given by $$\left\{
\begin{array}{l}
I(1,\, 1):=\KK, \quad \textrm{and}   \\
I(m,\,n):=0 \quad \textrm{otherwise}.
\end{array} \right.$$

We denote by $(\Sy\textrm{-biMod}, \bc,\, I)$ this monoidal
category.

\begin{Rq}
We need to restrict compositions to connected graphs and connected
permutations in order to get a monoidal category (see Proposition
1.6 of \cite{V1}).
\end{Rq}

\begin{dei}[Properad]
A \emph{properad} is a monoid in the monoidal category
$(\Sy\textrm{-biMod}, \bc,\ I)$.
\end{dei}

\begin{Ex}
Let $M$ be a $\KK$-module and consider $\End(M):=\bigoplus_{m, \,
n\in \NN^*} \Hom_\KK(M^{\ot n},\, M^{\ot m})$. The permutation of
the inputs and the outputs of an element of $\Hom_\KK(M^{\ot n},\,
M^{\ot m})$ makes $\End(M)$ into an $\Sy$-bimodule. Once again,
$\End(M)$, endowed with the natural (connected) composition of
morphisms, is a properad.
\end{Ex}

A comonoid in $(\Sy\textrm{-biMod}, \bc,\ I)$ is called a
\emph{coproperad}.\\

To an $\Sy$-module $V$, we associate an $\Sy$-bimodule
$\widetilde{V}$ defined by
$$\left\{ \begin{array}{ll}
\widetilde{V}(1,\, n):=V(n) & \textrm{and}   \\
\widetilde{V}(m,\,n):=0 & \textrm{for} \ m>1.
\end{array} \right.$$
This defines an embedding of monoidal categories, that is
$\widetilde{V\circ W}= \widetilde{V} \bc \widetilde{W}$. The
category $(\Sy\textrm{-Mod}, \circ,\, I)$ is a full monoidal
subcategory of $(\Sy\textrm{-biMod}, \bc,\, I)$. Hence, an operad
is a properad.

Since the notion of properad includes the one of associative
algebras and operads, we work in this general framework throughout
the text. We resume these notions in the following table

$$\xymatrix@R=6pt{\textsc{Monoidal category :}  & {(\KK\textrm{-Mod},\,
\otk)\ \ } \ar@{>->}[r] & {(\Sy\textrm{-Mod}, \circ) \ \ }
\ar@{>->}[r]&
{(\Sy\textrm{-biMod}, \bc)\  } \\
\textsc{Monoid :} &  {\textrm{Associative algebra} \ \ }
\ar@{>->}[r] & {\textrm{Operad} \ \  }  \ar@{>->}[r]&
{\textrm{Properad}.} }$$

Remark that the first monoidal product $\otk$ is bilinear and
symmetric, the second one $\circ$ is only linear on the left and
no symmetry. The third one $\bc$ has no linear nor symmetric
properties in general.

\subsection{$\Po$-gebras}

In this section, we precise the previous analogy with multi-linear
operations and recall the notion of an (al)gebra over a properad.\\

Let $(A,\, \mu,\, u)$ be an associative algebra and let $M$ be a
$\KK$-module. Recall that a structure of \emph{module} over $A$ on
$M$ is given by a morphism of associative algebras $\phi \, : \, A
\to \End(M)$. More generally, we have the following definition.

\begin{dei}[$\Po$-gebras]
Let $\Po$ be a properad and let $M$ be a $\KK$-module. A structure
of \emph{$\Po$-gebra} on $M$ is a morphism of properads $\phi \,
:\, \Po \to \End(M)$.
\end{dei}

When $\Po$ is an operad, this corresponds to the notion of algebra
over $\Po$ or \emph{$\Po$-algebra} (see V. Ginzburg and M.
Kapranov \cite{GK}). There is an operad $\A$ such that the
category of $\A$-algebras is equal to the category of non-unital
associative algebras, an operad $\C$ such that the category of
$\C$-algebras is equal to the category of non-unital commutative
associative algebras and an operad $\Li$ such that the category of
$\Li$-algebras is equal to the
category of Lie algebras.\\

Categories of ``algebras'' defined by products and coproducts
(multiple outputs), cannot be modelled by operads, one has to use
properads. Recall from \cite{V1}, that there is a properad $\Bi$
such that the category of $\Bi$-gebras is equal to the category of
non-unital non-counital bialgebras and there exists a properad
$\BLi$ such that category of $\BLi$-gebras is equal to the
category of Lie bialgebras, for instance.

\begin{Rq}
Following the article of J-.P. Serre \cite{Serre}, we choose to
call a \emph{gebra} any algebraic structure like modules over an
associative algebra, associative algebras, Lie algebras,
commutative algebras or bialgebras, Lie bialgebras, etc ...
\end{Rq}

This point of view on algebraic structures allows us to understand
and describe general properties between different types of gebras.
Constructions on the levels of operads or properads induce general
relations between the related types of gebras.

\subsection{Free and quadratic properad}

The forgetful functor $U$ from the category of properads to the
category of $\Sy$-bimodules has a left adjoint $\F$.

$$\xymatrix{{U \, : \, \textrm{Properads}  \ } \ar@_{->}@<-1ex>[r] &
\ar@_{->}@<-1ex>[l] {\ \Sy\textrm{-biMod} \,  : \, \F .}} $$

We gave an explicit construction of it in \cite{V4} by means of a
particular colimit. For every $\Sy$-bimodule $V$, it provides the
free properad $\F(V)$ on $V$. It is given by the direct sum of
connected directed graphs with the vertices indexed by elements of
$V$. The composition product $\mu$ is simply defined by the
grafting of graphs. Therefore, the number of vertices is preserved
by $\mu$
and it induces a natural graduation denoted $\F_{(\omega)}(V)$ and called the \emph{weight}.\\

Remark that, when $V$ is a $\KK$-module, we find the tensor
algebra $T(V)$ on $V$, which is the free associative algebra on
$V$. When $V$ is an $\Sy$-module, we get the free operad in terms
of
indexed trees like in \cite{GK} Section 2.1.\\

We can generalize the notion of ideal for an associative algebra
to ideals for operads and properads (see Appendix~\ref{Categorical
algebra}). Let $V$ be an $\Sy$-bimodule and $R$ be a
sub-$\Sy$-bimodule of $\F(V)$, we consider the ideal generated by
$R$ in $\F(V)$ and we denote it by $(R)$. As usual, the quotient
$\F(V)/(R)$ has a natural structure of properads. When $R\subset
\F_{(2)}(V)$, the quotient properad is called a \emph{quadratic
properad}. When $V$ is a $\KK$-module, this definition corresponds
to the notion of \emph{quadratic algebra} (see Y. Manin \cite{M})
and when $V$ is an $\Sy$-module it corresponds to the notion of
\emph{quadratic operad} of \cite{GK}.\\

\begin{Exs}
The symmetric and the exterior algebras are natural examples of
quadratic algebras. The operads $\A$, $\C$ and $\Li$ are the most
common quadratic operads (see \cite{GK}). The properads $\BLi$ of
Lie bialgebras, $\IBi$ of infinitesimal bialgebras and
$\mathcal{F}rob$ of Frobenius bialgebras are quadratic properads
(see \cite{V1} Section 2.9).
\end{Exs}

Since $R$ is homogenous of weight $2$, the quotient properad
$\F(V)/(R)$ is graded by the weight.\\

Dually, there is a \emph{connected cofree coproperad} denoted
$\F^c(V)$ (see \cite{V1} Section 2.8).

\subsection{Hadamard tensor product}








We define another monoidal product in the category of
$\Sy$-bimodules.

\begin{dei}[Hadamard product of $\Sy$-bimodules]
Let $V$ and $W$ be two $\Sy$-bimodules. Their Hadamard product is
defined by $(V\oth W)(m,\, n):=V(m,\, n) \otk W(m,\, n)$.
\end{dei}

When $V$ and $W$ are $\Sy$-modules, the Hadamard product is equal
to $(V\oth W)(n):=V(n) \otk W(n)$. When it is obvious that in the
context we are dealing
with the Hadamard product, we simply denote it by $\ot$.\\

This monoidal product is bilinear and symmetric. The unit of the
Hadamard product is the $\Sy$-bimodule $K$ defined by
$K(m,\,n):=\KK$, with trivial action of $\Sy_n$ and $\Sy_m$, for
all $n,m$ (and $K(n)=\KK$ for $\Sy$-modules). (The properad $K$ models commutative and
cocommutative Frobenius algebras).\\

\section{Categorical algebra}
\label{Categorical algebra} The aim of this section is to define
the notion of ``ideal of a monoid'' in a modern, categorical point
of view. Working in the opposite category, we get the dual notion
for comonoids. The other purpose of this categorical treatment is
to characterize the ideal ``generated by'' and its dual notion.

\subsection{Definition of the ``ideal'' notions}

In this section we define the notions of \emph{ideal monomorphism}
and \emph{ideal subobject} of monoids. Dually, we define the
notions of \emph{coideal epimorphism} and
\emph{coideal quotient}. \\

Let us work in an abelian monoidal category, that is an abelian
category $\CA$ endowed with a monoidal product $\bt$. Consider the
subcategory $\mona$ whose objects ar monoids in $\CA$. One natural
question now is to ask whether $\mona$ is still an abelian
category. The answer is no because the class of monomorphisms
 is too wide, for instance. Recall that in an abelian category the
class of monomorphisms is equal to the class of kernels. Every
 morphism $A\xrightarrow{f} B$ in $\mona$ admits a kernel $i\, :\, K\mono
A$ is $\CA$. The following diagram is commutative
$$\xymatrix{K \bt K \ar[r]^{i\bt i} \ar@{..>}[d]^{\mu_K}  & A\bt A \ar[r]^{f\bt f}
\ar[d]^{\mu_A}
 & B \bt B  \ar[d]^{\mu_B}\\
K  \ar@{>->}[r]^{i} & A  \ar[r]^{f} & B, }$$ where $\mu_A$ and
$\mu_B$ stand for the product of the monoid $A$ and $B$
respectively. Since $(f\bt f)\circ (i \bt i) =(f\circ i)\bt(f\circ
i) =0$, the composite $f \circ \mu_A \circ (i\bt i)$ is equal to
$0$. By the universal property of the kernel $i$, there exists an
associative map $\mu_K \, :\, K\bt K \to K$ making $K$ into a
monoid and $i \, :\, K\mono A$ a morphism in $\mona$. Hence
kernels exist in $\mona$ and every kernel is a monomorphism.
Actually, $K$ has more properties than just being a submonoid of
$A$ (see \ref{ideal-classical}), which explains why not all
monomorphisms are kernels. On the other hand, let $I \mono A$ be a
monomorphism of monoids, its cokernel $A/I$ in $\CA$ is not
necessarily a monoid. Following Kummer's language, we restrict our
attention to \emph{ideal monomorphisms}, that is monomorphisms
that are kernels in $\mona$.

\begin{dei}[Ideal monomorphism]
Let $I \mono A \epi Q$ be an exact sequence in $\CA$, where $A$ is
a monoid. In other words, $I \mono A$ is the kernel of $A \epi Q$
and $A\epi Q$ is the cokernel of $I \mono A$.

The monomorphism $I\mono A$ in $\mona$ is an \emph{ideal
monomorphism} if $A \epi Q$ is a morphism in $\mona$.
\end{dei}

In this case, we say that $I$ is an \emph{ideal} (subobject) of
$A$ and $Q$ is naturally a \emph{quotient} monoid, also denoted by $A/I$.\\

Dually, recall that a comonoid in $\CA$ is a monoid in the
opposite category $\CA^{\textrm{op}}$. If we dualize the previous
arguments in the opposite category, we can see that the category
$\coma$ of comonoids in $\CA$ is not an abelian category because
the class of epimorphisms is too big. The cokernel in $\CA$ of a
morphism in $\coma$ is a morphism in $\coma$ (and even more), but
the kernel in $\CA$ of a morphism in $\coma$ is not necessarily a
morphism in $\coma$. Therefore, every epimorphism of comonoids is
not a cokernel. Once again, we call \emph{coideal epimorphisms},
the epimorphisms that are cokernels.

\begin{dei}[Coideal epimorphism]
Let $I \mono C \epi Q$ be an exact sequence in $\CA$, where $C$ is
a comonoid.  The epimorphism $C\epi Q$ in $\coma$ is a
\emph{coideal epimorphism} if $I \mono C$ is a morphism in
$\coma$.
\end{dei}

In this case, the subobject $I\mono C$ is naturally a
\emph{subcomonoid} of $C$ and the quotient $Q$ is called a
\emph{coideal quotient}.

\begin{Rq}
The term  \emph{coideal} is already used in the literature, but
stands for a (coideal) subobject $J \mono C$ (or a monomorphism)
in $\CA$ of a comonoid $C$ such that its cokernel in $\CA$ is a
morphism in $\coma$. It is equivalent to ask that the quotient
$C/J$ is a comonoid.

This notion does not correspond to the dual of the notion of
ideal, where ``dual'' means ``in the opposite category''.
\end{Rq}

\subsection{Relation with the classical definition}
\label{ideal-classical} We now relate this definition of ideal
with the classical notion. Let $\xymatrix@C=20pt{ I\
\ar@{>->}[r]^{\iota}& A \ar@{->>}[r]^{\pi}& Q}$ be an sequence in
$\mona$, exact in $\CA$. Denote by $\mu_A$ and $\mu_Q$ the
products of $A$ and $Q$ respectively. The morphism $\pi$ is a
morphism in $\mona$, means that the following diagram commutes
$$\xymatrix { A\bt A \ar[r]^{\pi \bt \pi} \ar[d]^{\mu_A}&  Q\bt Q \ar[d]^{\mu_Q} \\
A \ar[r]^{\pi} & Q. }$$ Let $\kappa \, :\, K_I \mono A\bt A$ be
the kernel of $\pi \bt \pi$ in $\CA$.

\begin{pro}
\label{ideal-kernel}
 Let $\CA$ be a monoidal category such that
the monoidal product $\bt$ preserves epimorphisms. A monomorphism
$I\mono A$ is an ideal monomorphism if and only if the composite
$\pi \circ \mu_A \circ \kappa$ is equal to $0$.
\end{pro}

\begin{proo}$ \ $

$(\Rightarrow)$ It comes from $\pi \circ \mu_A \circ \kappa= \mu_Q
\circ  (\pi \bt \pi) \circ \kappa=0$

$(\Leftarrow)$ From the hypothesis, we have that $\pi \bt \pi$ is
an epimorphism. Therefore, it is the cokernel of $\kappa$ and by
the universal property of the cokernel, there exists a morphism
$\mu_Q\, : \, Q\bt Q \to Q$ such that $\pi \circ \mu_A  $ factors
through $\pi \bt \pi$. It is then straightforward to check that
$\mu_Q$ defines an associative product on $Q$.
\end{proo}

The extra assumption, that the monoidal product has to preserve
epimorphisms, is verified in every cases studied in this paper.
For a proof for this fact, we refer to Proposition 1,  Proposition
10 and Section 5 of \cite{V4}.

The problem is now to make explicit the kernel $\kappa \, :\, K_I
\mono A\bt A$ of $\pi \bt \pi$.

\begin{dei}[Multilinear part]
The \emph{multilinear part} in $X$ of $A\bt ( X \oplus Y)\bt B$ is
defined either
\begin{itemize}
\item by the cokernel of $A\bt Y \bt B \xrightarrow{A\bt i_y \bt B
} A\bt ( X \oplus Y)\bt B$

\item or by the kernel of $A\bt (X \oplus Y) \bt B
\xrightarrow{A\bt \pi_y \bt B } A\bt  Y \bt B$,
\end{itemize}
since $i_Y$ is a section of $\pi_Y$, that is $\pi_y \circ i_Y =
\textrm{Id}_Y$, these two objects are naturally isomorphic. We
denote it $A\bt (\underline{X} \oplus Y) \bt B$.
\end{dei}

It corresponds to elements of $A\bt (X \oplus Y) \bt B$ with at
least one element of $X$ in between.

Suppose that we are working in an abelian category $\CA$ such that
every short exact sequence splits, that is $\xymatrix{{I \ }
\ar@{>->}[r]^{\iota} & A \ar@{->>}[r]_{\pi}&
\ar@{..>}@/_/[l]_{\exists} Q }$ or equivalently $A\cong I \oplus
Q$. Once again, this condition is verified in every category
studied here since they are categories of representations of
finite groups over a field of characteristic $0$.

\begin{pro}\label{Ideal-Caracterisation}
In a monoidal abelian category such that the monoidal product
preserves epimorphisms and where every short exact sequence
splits, we have
$$K_I=A\bt (A + \underline{I})+ (A +
\underline{I})\bt A,$$ where $A\bt (A+ \underline{I}):=\textrm{Im}
\, \left( A \bt (A \oplus \underline{I}) \mono A \bt (A\oplus I)
\xrightarrow{A\bt (A + \iota)} A\bt A \right)$.
\end{pro}

\begin{proo}
It is enough to prove that $A\bt (Q + \uI ) +  (Q + \uI )\bt A$ is
the kernel of $\pi \bt \pi$. We have the following commutative
diagram
$$\xymatrix@M=9pt{Q \bt (Q\oplus \uI) \ar@{>->}[r]& Q \bt A \ar@{->>}[r]^{Q \bt \pi} & Q \bt Q \\
A \bt (Q \oplus \uI)  \ar@{>->}[r]^{i_1} \ar@{>->}[d]
\ar@{..>>}[u]^{\exists !\, \theta_2} & A \bt A \ar@{->>}[r]^{A \bt
\pi} \ar@{->>}[u]_{\pi \bt A} \ar@{->>}[ur]^{\pi \bt \pi}  &
 A \bt Q \ar@{->>}[u]_{\pi \bt Q} \\
A \bt (Q \oplus \uI) \oplus  (Q\ \oplus \uI) \bt A \ar[ur]^{i_1
\oplus i_2} & \ar@{>->}[l] (Q\ \oplus \uI) \bt A
\ar@{>->}[u]_{i_2} \ar@{..>>}[r]^{\exists ! \, \theta_2} & (Q
\oplus \uI)\bt Q, \ar@{>->}[u] }
$$ where the two dotted arrows exist by the property of kernels
applied to the first line and last column. Since $A\bt (Q +\uI ) +
(Q + \uI )\bt A$ is by definition the image of the morphism
$i_1\oplus i_2$, it remains to show that $\pi \bt \pi$ is the
cokernel of $i_1\oplus i_2$. The assumption that every short exact
sequence splits implies that the maps $\theta_1$ and $\theta_2$
are epimorphisms. It is then straightforward to check that $\pi
\bt \pi$ is the cokernel of $i_1\oplus i_2$.
\end{proo}

Since $\iota \, : \, I \mono A$ is the kernel of $\pi \, : \, A
\epi Q$, there exists a morphism $\bar{\mu}$ making the following
diagram to commute
$$\xymatrix@M=6pt@W=6pt@H=10pt{A\bt(A + \underline{I}) + (A+
\underline{I})\bt A \ar@{>->}[r] \ar[d]^{\bar{\mu}} & A\bt A
\ar[d]^{\mu_A} & \\
{I \ } \ar@{>->}[r]^{\iota}& A \ar@{->>}[r]^{\pi}& Q.  }$$ Hence,
we get
 $\pi \circ \iota \circ \bar{\mu}=0$.
When the monoidal product is additive on the left and on the
right, we have $A\ot (A+ \underline{I})=A\ot I$ and $(A+
\underline{I})\ot A=I\ot A$. In this case, the notion of ideal
corresponds to the classical one.\\

Dually, let $\xymatrix@C=20pt{ I\  \ar@{>->}[r]^{\iota}& C
\ar@{->>}[r]^{\pi}& Q} $ be a sequence in $\coma$, exact in $\CA$.
Denote by $\gamma \, :\, C\bt C \epi CoK_Q$ be the cokernel of
$\iota \bt \iota$. Note that when every short exact sequence
splits, we have $\textrm{CoK}_Q\cong K_I\cong C\bt (C +
\underline{I})+ (C + \underline{I})\bt C$.

\begin{pro}
Let $\CA$ be a monoidal category such that the monoidal product
$\bt$ preserves monomorphisms. An epimorphism $C \epi Q$ is an
coideal epimorphism if and only if the composite $\gamma \circ
\Delta_C \circ \iota$ is equal to $0$.
\end{pro}

\begin{proo}
We work in the opposite category and we apply
Proposition~\ref{ideal-kernel}.
\end{proo}

In the case of a coassociative coalgebra $C=I \oplus Q$, it means
that the composite $(\pi\otimes C \oplus C\otimes \pi) \circ
\Delta_C\circ \iota \, :\, I\to Q\otimes C \oplus C\otimes Q$ is
null. In other words, we have $\Delta_C(c)\in I \otimes I$, for
$c\in I$.

\subsection{Various notions of modules}

We recall briefly the various notions of modules and relate one of
them to the notion of ideal.

\begin{dei}[Module]
An object $M$ of $\CA$ is a \emph{left module} over a monoid $A$
if there exists a map $A\bt M \to M$ compatible with the product
of $A$. Dually, there is a notion of \emph{right module}. And A
compatible left and right action defines a \emph{bimodule}.
\end{dei}

At first sight, the biadditive case could lead to the following
definition : $I$ is an ideal of $A$ if it is a bimodule over $A$ :
$A \bt I \to I$ and $I \bt A \to I$. The main problem with such a
notion is that $A/I$ is not a monoid when $\bt$ is not biadditive.
Instead of that, one has to consider a linearized version of
module.

\begin{dei}
An object $M$ of $\CA$ is called a \emph{multilinear left module
over $A$} if it is endowed with a map $A \bt (A \oplus
\underline{M}) \to M$ compatible with the product and the unit of
$A$.
\end{dei}

We have a similar notion on the right hand sight and a notion of
\emph{multilinear bimodule}. If we use this language, $I$ is an
ideal of $A$ if and only if $I$ a multilinear bimodule with the
action induced by $\mu_A$.

\begin{Rq}
The same notion arises from the work of D. Quillen on (co)homology
theories \cite{Quillen}. The coefficient for these theories are
\emph{abelian group objects}. When one wants to make explicit
Quillen (co)homology of monoids, these coefficients are exactly
linear version of modules. We refer the reader to the paper of
H.J. Baues, M. Jibladze and A. Tonks \cite{BJT} for a complete
description in the case of operads, or more generally when the
monoidal product is additive only on one side.
\end{Rq}

Dualize these definitions to get the notions of \emph{comodules}
and \emph{(multi)linear} comodules over a comonoid $C$.

\subsection{``generated by''}

\label{generatedby} Following this categorical point of view, we
define and make explicit the notions of \emph{ideal generated by}
and \emph{subcomonoid
generated by}.\\

Let $\xi\, : \, R\mono A$ be a subobject of $A$ in $\CA$, where
$A$ is a monoid. We are going to consider the ``cokernel'' $A \epi
Q $ of $\xi$ in $\mona$, that is the universal epimorphism of
monoids such that the composite  $R \mono A\epi Q$ vanishes. The
resulting quotient monoid $Q$ is the largest quotient of $A$ with relations in $R$.\\

Since our leitmotiv is to treat together ideals and quotient
monoids, we would rather use the following presentation. Consider
the category $\mathcal{S}_\xi$ of sequences $(\textbf{S})\, :I
\mono A \epi Q $ in $\mona$, exact in $\CA$ such that the
composite $R \mono A \epi Q$ is equal to $0$. Since $I\mono A$ is
the kernel of $A \epi Q$, this last condition is equivalent to the
existence of a morphism $\iota \, :\, R \mono I$ in $\CA$  such
that the following diagram commutes

$$\xymatrix@M=6pt@W=6pt@H=10pt{R \ar@/_/[drr]^{0}  \ar@{>->}[dr]_{\xi} \ar@{>..>}[d]_{\exists \, \iota}& & \\
I\ar@{>->}[r] & A \ar@{->>}[r]  &  Q.}$$

Let $(\textbf{S}')\, :J \mono A \epi O$ be another object of
$\mathcal{S}_\xi$, the morphisms between ($\textbf{S}$) and
($\textbf{S}'$) correspond to the pair of morphisms $(i,\, p)$ in
$\mona$ such that the following diagram commutes

$$\xymatrix@M=6pt@W=6pt@H=10pt{R  \ar@{>->}[dr]^{\xi} \ar@{>->}[d] \ar@/_2pc/@{>->}[dd]& & O \\
I\ar@{>->}[r] \ar@{>.>}[d]_{i}  & A \ar@{->>}[ur] \ar@{->>}[r]  &  Q \ar@{.>>}[u]_{p} \\
J \ar@{>->}[ur]&& }$$

\begin{dei}[Ideal generated by $R$]
Let $\CA$ be a category such that for every monoid $A$ and every
subobject $\xi\, : \, R\mono A$, the category $\mathcal{S}_\xi$
admits an initial object $\bar{\textbf{S}} \, : \, (R) \mono A
\epi A/(R)$.

In this case, $(R)$ is called the \emph{ideal generated by $R$}
and $A/(R)$ is the induced quotient monoid.
\end{dei}

\begin{Rq}
The terminal object of this category always exists and is given by
the sequence $A\mono A \epi 0 $.
\end{Rq}

If we dualize the previous arguments in the opposite category, we
get the same kind of diagram but with $C$ comonoid instead of $A$
monoid.\\

Let $\xi\, : \, S \Lepi C$ be a quotient of $C$ in $\CA$, where
$C$ is a comonoid. We aim to consider the largest subcomonoid of
$C$ vanishing on $S$. This notion is given by ``kernel'' $S \Lepi
C$ of $\xi$ in $\coma$, that is the universal monomorphism of
comonoids such that the composite $S \Lepi C\Lmono Q$ is equal to
$0$.\\

Consider the category $\mathcal{S}_\xi$ of sequences
$(\textbf{S})\, :I \Lepi C \Lmono Q $ in $\coma$, exact in $\CA$
such that the composite $S \Lepi C \Lmono Q$ is equal to $0$.
There exists a morphism $\pi \, :\, I \epi S$ such that the
diagram commutes

$$\xymatrix@M=6pt@W=6pt@H=10pt{S    & & \\
I \ar@{..>>}[u]^{\exists \, \pi}& C \ar@{->>}[l]
\ar@{->>}[ul]^{\xi}& \ar@{>->}[l] Q \ar@/^/[ull]_{0} .}$$

Let $(\textbf{S}')\, :J \Lepi C \Lmono O$ be another object of
$\mathcal{S}_\xi$, the morphisms between ($\textbf{S}$) and
($\textbf{S}'$) correspond to the pair of morphisms $(i,\, p)$ in
$\coma$ such that the following diagram commutes

$$\xymatrix@M=6pt@W=6pt@H=10pt{S  & &
O\ar@{>->}[dl] \ar@{>.>}[d]^{i} \\
I \ar@{->>}[u]   & C  \ar@{->>}[l] \ar@{->>}[ul]_{\xi} \ar@{->>}[dl]  & \ar@{>->}[l] Q  \\
J \ar@/^2pc/@{->>}[uu] \ar@{.>>}[u]^{p} & & }$$

\begin{dei}[Subcomonoid generated by $S$]
Let $\CA$ be a category such that for every comonoid $C$ and every
quotient $\xi\, : \,  S \Lepi C$, the category $\mathcal{S}_\xi$
admits a terminal object $\bar{\textbf{S}} \, : \, (S) \Lepi C
\Lmono C(S)$.

In this case, $C(S)$ is called the \emph{subcomonoid of $C$
generated by $S$} and $(S)$ is the induced coideal quotient.
\end{dei}

\begin{Rq}
The initial object is the sequence $C\Lepi C  \Lmono 0 $.
\end{Rq}

\subsection{Ideal generated = free multilinear bimodule}
 \label{freebimodule}
Since the notion of ideal is equivalent to the notion of
multilinear bimodule, the ideal of $A$ generated by $R$ is the
free $A$-multilinear bimodule on $R$.

\begin{pro}
The ideal generated by $R$ in $A$ is given by the image
$$\textrm{Im}\left(A\bt (A+\underline{R})\bt A \xrightarrow{\mu^2} A           \right). $$
\end{pro}

\begin{proo}
Using Proposition~\ref{Ideal-Caracterisation}, we have that it is
an ideal of $A$. It is easy to see that any ideal containing $R$
also contains $\textrm{Im}\left(A\bt (A+\underline{R})\bt A
\xrightarrow{\mu^2} A \right)$.
\end{proo}

If we dualize the arguments, we have the explicit form of the
subcomonoid of $C$ generated by $R$.

\begin{pro}
Let $S \Lepi C$ be an epimorphism in $\CA$. The subcomonoid of $C$
generated by $S$ is given by the kernel of
$$\textrm{Ker}\left(C \xrightarrow{\Delta^2} C^{\bt 3} \xrightarrow{proj} C\bt (C+\underline{S}) \bt C \right). $$
\end{pro}

\begin{proo}
Dualize.
\end{proo}

\begin{center}
\textsc{Acknowledgements}
\end{center}

I would like to thank Marcelo Aguiar to have shared with me his
knowledge on Manin's products. I am grateful to Jean-Louis Loday,
 Damien Calaque and Martin Markl for useful discussions. It is a pleasure to thank Anton
Korochkine for stimulating discussions about the counterexample of
Section~\ref{Counterexample}.\\

I am grateful to the Max Planck Institute for Mathematics in Bonn
for the excellent working conditions I found there during my stay
in May and June 2006 where I finished this paper.

\medskip

{\small \textsc{Laboratoire J.A. Dieudonn\'e, Universit\'e de
Nice, Parc Valrose, 06108 Nice
Cedex 02, France}\\
E-mail address : \texttt{brunov@math.unice.fr}\\
URL : \texttt{http://math.unice.fr/$\sim$brunov}}

\end{document}